\documentclass[11pt]{article}
\usepackage{amsfonts}
\usepackage{amsfonts,latexsym,amsmath,amscd,geometry}
\usepackage{amssymb}
\usepackage{latexsym}
\renewcommand{\theequation}{\thesection.\arabic{equation}}
\newcommand \nc{\newcommand}
\newtheorem{theorem}{Theorem}[section]
\newtheorem{lemma}[theorem]{Lemma}
\newtheorem{proposition}[theorem]{Proposition}
\newtheorem{corollary}[theorem]{Corollary}
\newtheorem{definition}[theorem]{Definition}

\newtheorem{remark}[theorem]{Remark}
\renewcommand{\thetheorem}{\thesubsection.\arabic{theorem}}
\nc{\ba}{\begin{array}}\nc{\ea}{\end{array}}
\nc{\be}{\begin{eqnarray}}\nc{\ee}{\end{eqnarray}}
\nc{\beq}{\begin{equation}}\nc{\eeq}{\end{equation}}
\nc{\bex}{\begin{eqnarray*}}\nc{\eex}{\end{eqnarray*}}
\nc{\btm}{\begin{theorem}} \nc{\etm}{\end{theorem}}
\nc{\blm}{\begin{lemma}} \nc{\elm}{\end{lemma}}
\nc{\R}{\mathbb{R}} \nc{\va}{\varepsilon} \nc{\ls}{\limits}

\def\de{\Delta}
\def\pf{\noindent{\bf Proof.\quad}}\def\endpf{\hfill$\Box$}
\def\u{\dot{u}}\def\di{\mbox{div\,}}

\begin{document}
\title{{\bf Global solutions to the three-dimensional full compressible Navier-Stokes equations with vacuum at infinity in some classes of large data}}
\author{Huanyao Wen\thanks{Department of Mathematics, South China University of Technology, Guangzhou, 510641, China.
Email: mahywen@scut.edu.cn, or huanyaowen@hotmail.com},\quad
Changjiang Zhu\thanks{Department of Mathematics, South China
University of Technology, Guangzhou, 510641, China. Email:
cjzhu@mail.ccnu.edu.cn, or machjzhu@scut.edu.cn}  }
\date{}

\maketitle

\begin{abstract}
We consider the Cauchy problem for the full compressible
Navier-Stokes equations with vanishing of density at infinity in
$\mathbb{R}^3$. Our main purpose is to prove the existence (and
uniqueness) of global strong and classical solutions and study the
large-time behavior of the solutions as well as the decay rates in
time. Our main results show that the strong solution exists globally
in time if the initial mass is small for the fixed coefficients of
viscosity and heat conduction, and can be large for the large
coefficients of viscosity and heat conduction. Moreover, large-time
behavior and a surprisingly exponential decay rate of the strong
solution are obtained. Finally, we show that the global strong
solution can become classical if the initial data is more regular.
Note that the assumptions on the initial density do not exclude that
the initial density may vanish in a subset of $\mathbb{R}^3$ and
that it can be of a nontrivially compact support.
 To our knowledge, this paper contains the first result so far for the global existence of solutions to the full compressible
  Navier-Stokes equations when density vanishes at infinity (in space). In addition, the exponential decay rate of the strong
   solution is of independent interest.
\end{abstract}

\noindent{\bf Key Words}: Full compressible Navier-Stokes equations, global classical and strong solutions, large-time behavior, vacuum.\\[0.8mm]
\noindent{\bf 2000 Mathematics Subject Classification}. 35Q30,
 76N10, 35K65.

\tableofcontents

\vspace{4mm}
\section {Introduction}
\setcounter{equation}{0} \setcounter{theorem}{0}
 The compressible Navier-Stokes
equations, describing the motion of compressible fluids, can be written in the Eulerian coordinates in
$\mathbb{R}^3$ as follows: \be\label{3d-full N-S}
\begin{cases}
\rho_t+\mathrm{div}(\rho u)=0, \\
(\rho u)_t+\mathrm{div}(\rho u\otimes u)+\nabla
P=\mathrm{div}(\mathcal {T}),\\
(\rho E)_t+\mathrm{div}(\rho E u+P u)=\mathrm{div}(\mathcal {T}
u)+\mathrm{div}(\kappa\nabla\theta).
\end{cases}
\ee Here $\mathcal {T}$ is the stress tensor given by
$$\mathcal
{T}=\mu\left(\nabla u+(\nabla u)^\prime\right)+\lambda \di u I_3,$$
where $I_3$ is a $3\times 3$ unit matrix;
 $\rho=\rho(x,t)$, $u=u(x, t)=(u_1, u_2, u_3)(x,t)$ and $\theta=\theta(x,t)$
are unknown functions denoting  the density, velocity and absolute
temperature, respectively; $P=P(\rho,\theta)$, $E$ and $\kappa$
denote
 pressure, total energy, and coefficient
of heat conduction, respectively, where $\displaystyle
E=e+\frac{|u|^2}{2}$ ($e$ is the internal energy);
 $\mu$ and $\lambda$ are coefficients of viscosity, satisfying the following physical restrictions
\be\label{viscosity assumption} \mu>0,\ \lambda+\frac{2\mu}{3}\ge0.
\ee Assume that \bex P=R\rho\theta,\ e=C_\nu\theta, \eex for some
 constants $R>0$ and $C_\nu>0$. Note that the temperature function can be written as $\theta=\mathcal{A}\exp\{S/C_\nu\}\rho^{\gamma-1}$ where $S=S(x,t)$ is the entropy, $ \mathcal{A}>0$ and $\gamma>1$ are constants.

When the density has a positive lower bound (i.e., no vacuum at any
point), the momentum equation and the energy equation are parabolic.
In this case, there have been a lot of works so far on the
well-posedness of solutions to the Cauchy problem and the
initial-boundary-value problem for (\ref{3d-full N-S}). Refer, for
instance, to these elegant works \cite{Itaya, Jiang1, Jiang:Math
Nachr, Kawohl, Kazhikhov-Shelukhi, Matsumura-Nishida: Kyoto Un,
Matsumura-Nishida: CMP, Tani} for local and global existence of
classical solutions from one dimension to high dimensions. In
particular, Matsumura and Nishida in \cite{Matsumura-Nishida: Kyoto
Un, Matsumura-Nishida: CMP} showed that the global classical
solution in three dimensions exists provided that the initial data
is small in $H^3$. In one dimension and high dimensions with
spherically symmetric solutions, the global existence of classical
solutions with large initial data has been obtained (see
\cite{Kazhikhov-Shelukhi, Kawohl, Jiang1}). On the existence,
asymptotic behavior of the weak solutions to the system
(\ref{3d-full N-S}), please refer, for instance, to \cite{ Jiang2,
Jiang3, Jiang-Zhang:weak solutions, Hoff, Hoff-arma, Feireisl1} from
one dimension to high dimensions.

When vacuum is allowed, some new challenging difficulties arise,
such as degeneracy of the equations. In spite of this, some
important progress on global existence of weak solutions, local
existence of strong solutions and global existence of classical
solutions has been achieved by Feireisl, Bresch-Desjardins, Cho-Kim
and Huang-Li \cite{Feireisl-book, Bresch-Desjardins, cho-Kim:
perfect gas, Huang-Li}. More precisely, Feireisl in his pioneering
work \cite{Feireisl-book} got the global existence of so-called {\em
variational} solutions to (\ref{3d-full N-S}) with
temperature-dependent coefficient of heat conduction in a bounded
domain $\Omega\subseteq\mathbb{R}^N$ for $N\ge 2$. The temperature
equation in \cite{Feireisl-book} is satisfied only as an inequality
in sense of distribution. This work is the very first attempt
towards the existence of weak solutions to the full compressible
Navier-Stokes equations in high dimensions. With viscosity
coefficients which are only density-dependent, the existence of
global weak solutions to the Navier-Stokes equations for viscous
compressible and heat conducting fluids in $\mathbb{T}^3$ or
$\mathbb{R}^3$ was obtained by Bresch and Desjardins
\cite{Bresch-Desjardins}. One of the important estimates in
\cite{Bresch-Desjardins} is the Bresch-Desjardins entropy estimate
which gives more regularity of the density function with the help of
the density-dependent viscosity. As pointed out in
\cite{Bresch-Desjardins}, the assumption on the initial density that
$\rho_0-\rho_\infty\in L^1$ is necessary for some positive constant
$\rho_\infty$. Huang and Li \cite{Huang-Li} obtained the existence
and uniqueness of global classical solutions to the Cauchy problem
for (\ref{3d-full N-S}) in $\mathbb{R}^3$ provided that the initial
energy is small. Vacuum is allowed interiorly but not at infinity in
\cite{Huang-Li}. When the domain is bounded, with vacuum and large
initial data, we got existence and uniqueness of global classical
solutions to the initial-boundary-value problem of (\ref{3d-full
N-S}) in one dimension and in high dimensions with spherically and
cylindrically symmetric initial data \cite{Wen-Zhu, Wen-Zhu 3}.
Temperature-dependent coefficient of heat conduction plays a crucial
role in the proof.

When the initial density and initial temperature vanish at
infinity\footnote{In this case, initial vacuum state at infinity is
allowed.}, there is no useful basic energy equality (or inequality).
To our knowledge, there is only one result so far in this direction,
please refer to \cite{cho-Kim: perfect gas} for the perfect gas with
constant coefficients of viscosity and heat conduction, where the
authors obtained the local existence and uniqueness of strong
solutions to (\ref{3d-full N-S}) with vacuum at infinity in
$\mathbb{R}^3$. Two natural questions are: do there exist some
global strong and more regular solutions to (\ref{3d-full N-S}) when
vacuum state is allowed at infinity$?$ If do, then how do they
behave when time goes to infinity? We will answer the questions in
the paper.

We would like to introduce our main ideas in the paper. To prove the
global existence of the strong solution, we establish a sharper
blow-up criterion than that we have obtained in \cite{Wen-Zhu 4} for
strong solution if the strong solution blows up in finite time. Then
we get a crucial proposition (Proposition \ref{prop 3.1}) which
implies that the terms in the criterion will never blow up in finite
time when the initial mass is small in some sense (refer to the
proof of Corollary \ref{3d-cor 3.2} for more details). This together
with the contradiction arguments indicates that the strong solution
exists globally in time. This is the main ingredient of the proof.
Moreover, our result shows that the initial mass can be large if the
coefficients of viscosity and heat conduction are taken to be large,
which implies that large viscosity and heat conduction mean large
solution. Furthermore, large-time behavior of the strong solution is
considered, and a surprisingly exponential decay rate of the strong
solution is obtained. Finally, we show that the global strong
solution can become classical if the initial data is more regular.

The main challenges in studying the global well-posedness of
solution are summarized as follows:\\

$(D_1):$ No useful basic energy equality (or inequality).\\

When $\rho\rightarrow\tilde{\rho}>0,\
\theta\rightarrow\tilde{\theta}>0$ as $|x|\rightarrow\infty$, the
following classical basic energy equality holds: \bex\begin{split} \mathcal
{C}(t)+\int_0^t\int_{\mathbb{R}^3}\left(\frac{\lambda(\mathrm{div}u)^2
+\frac{\mu|\nabla u+(\nabla
u)^\prime|^2}{2}}{\theta}+\frac{\kappa|\nabla\theta|^2}{\theta^2}\right)=
\mathcal {C}(0),
\end{split}
\eex where
$$\mathcal
{C}(t)=\int\left(\rho
(\frac{\theta}{\tilde{\theta}}-\log\frac{\theta}{\tilde{\theta}}-1)+\frac{\rho
|u|^2}{2\tilde{\theta}}+(\tilde{\rho}-\rho+\rho\log\frac{\rho}{\tilde{\rho}})\right).
$$
It is easy to verify that $\mathcal {C}(0)\ge0$. This equality plays an important role in the proof of the main
theorems for instance in \cite{Huang-Li,Jiang1,Jiang-Zhang:weak solutions}. If $(\tilde{\rho},
\tilde{\theta})=(0,0)$, the energy equality (or inequality) is unavailable.\\

$(D_2):$ Zlotnik inequality (see  Appendix A in Section 7) which was
used in \cite{Huang-Li-Xin} for isentropic flow to get the upper
bound of the density does not work here. In \cite{Huang-Li-Xin},
$g(\rho)$ is defined as $g(\rho)=-\frac{\rho P(\rho)}{2\mu+\lambda}$
for the case that $\tilde{\rho}=0$, where $P(\rho)=a\rho^\gamma$ for
$a>0$ and $\gamma>1$. However, in the present paper,
$P=R\rho\theta$. Thus $g(\infty)\not=-\infty$ due to the possible
vanishing of
$\theta$. \\

Our strategies on handling $(D_1)$ and $(D_2)$ are as follows.
Firstly, for $(D_1)$, we define a new $B(T)$ (see (\ref{BT})) in
Proposition \ref{prop 3.1} and do not need time-weighted terms with
more regularity like those in \cite{Huang-Li}. Besides, we prove
that the mass is conserved for all time with the regularity of the
strong solution, i.e.,
$$
\int_{\mathbb{R}^3}\rho\,dx=\int_{\mathbb{R}^3}\rho_0\,dx.
$$ This gives that the mass is small for all time if we assume the initial mass is small.
 The ``smallness" of the mass and the {\it a priori} assumption of $B(T)$ and density make us get the estimate of $\nabla u$ in $L^2_{xt}$
  norm which is the very important starting point. Secondly, for
$(D_2)$, we use the idea of Lions (\cite{Lions2}) and Desjardins
(\cite{Desjardin}) to construct a ``$\log\rho$" equation. Then we
define a different function $g$ from the isentropic case, i.e.,
$g(\rho)=-\frac{R\rho\theta}{2\mu+\lambda}$. With the help of the
``$\log\rho$" equation, the ``smallness" of mass and the {\it a
priori} assumptions of $A(T)$ (see (\ref{AT})), $B(T)$ and $\rho$,
we have $N_1=0$ in Appendix A. This suggests that $g\le0$ is enough
in stead of $g(\infty)=-\infty$ in Appendix A. This is the main
ingredient in the proof of
the upper bound of the density. See Lemma \ref{3d-le:rho} for more details.\\

Before we state our main results, we would like to give some
notation which will be used throughout
this paper.\\

 \noindent{\bf Notation:}\\

(i)\ $\displaystyle\int_{\mathbb{R}^3} f =\int_{\mathbb{R}^3} f \,dx$.\\

(ii)\ For $1\le l\le \infty$, denote the $L^l$ spaces and the
standard Sobolev spaces as follows:
$$L^l=L^l(\Sigma),  \ D^{k,l}=\left\{ u\in L^1_{\rm{loc}}(\Sigma): \|\nabla^k u \|_{L^l}<\infty\right\},$$
$$W^{k,l}=L^l\cap D^{k,l},  \ H^k=W^{k,2}, \ D^k=D^{k,2},$$
$$D_0^1=\Big\{u\in L^6: \ \|\nabla u\|_{L^2}<\infty\},$$
$$\|u\|_{D^{k,l}}=\|\nabla^k u\|_{L^l}.$$

(iii)\ $G=(2\mu+\lambda)\mathrm{div}u-P$ is the effective viscous flux.\\

(iv)\ $\dot{h}=h_t+u\cdot\nabla h$ denotes the material
derivative.\\

(v)\ $\displaystyle m_0=\int_{\mathbb{R}^3}\rho_0(x)\,dx$.\\

The rest of the paper is organized as follows. In Section 2, we
present our main results. In Section \ref{3d-sec 3}, we establish a
sharp blow-up criterion for strong solution. In Section \ref{3d-sec
4}, motivated by the blow-up criterion established in Section
\ref{3d-sec 3}, we prove the global existence of strong solution provided that the initial mass is small in some sense. In
Section \ref{3d-sec 5}, we study the large-time behavior of the
solution and get the exponential decay estimate. In
Section \ref{3d-sec 7}, based on the local well-posedness of the
classical solution, we get the existence and uniqueness of global
classical solution by establishing some higher-order {\it a priori} estimates globally
in time. In Section \ref{3d-sec 8}, we give the proof of the local
well-posedness of classical solution.
\section{ Main results}
\setcounter{equation}{0} \setcounter{theorem}{0}
\renewcommand{\thetheorem}{\thesection.\arabic{theorem}}
 Assume that $\mu$,
$\lambda$ and $\kappa$ are constants. We assume $R=C_\nu=1$
henceforth, since the constants $R$ in the pressure function and
$C_\nu$ in the internal energy play no role in the analysis. In this
case, if the solutions are regular enough (such as strong solutions
and classical solutions), (\ref{3d-full N-S}) is equivalent to the
following system
 \be\label{full N-S+1}
\begin{cases}
\rho_t+\nabla \cdot (\rho u)=0, \\
\rho  u_t+\rho u\cdot\nabla u+\nabla
P(\rho,\theta)=\mu\Delta u+(\mu+\lambda)\nabla\mathrm{div}u,\\
\rho \theta_t+\rho
u\cdot\nabla\theta+\rho\theta\mathrm{div}u=\frac{\mu}{2}\left|\nabla
u+(\nabla
u)^\prime\right|^2+\lambda(\mathrm{div}u)^2+\kappa\Delta\theta,\
\mathrm{in}\ \mathbb{R}^3\times(0,\infty).
\end{cases}
\ee System (\ref{full N-S+1}) is supplemented with initial
conditions \be\label{3d-initial} (\rho, u, \theta)|_{t=0}=(\rho_0,
u_0, \theta_0),\ x\in\mathbb{R}^3, \ee with \be\label{3d-boundary}
\rho(x,t)\rightarrow0,\ u(x,t)\rightarrow0,\
\theta(x,t)\rightarrow0,\ \mathrm{as}\ |x|\rightarrow\infty,\
\mathrm{for}\ t\ge0. \ee

We give the definition of the strong solution to (\ref{full
N-S+1})-(\ref{3d-boundary}) throughout this paper, which is similar
to \cite{cho-Kim: perfect gas}.
\begin{definition}(Strong solution) For $T>0$, $(\rho, u, \theta)$ is called a strong solution to the compressible Navier-Stokes equations
(\ref{full N-S+1})-(\ref{3d-boundary}) in $\mathbb{R}^3\times
[0,T]$, if for some $q\in (3, 6)$, \be\label{global strong}
\begin{split}0\le \rho\in
C([0,T];W^{1,q}\cap H^1),\ \rho_t\in C([0,T];L^2\cap L^q),&\\
 (u, \theta)\in C([0,T];D^2\cap D^1_0)\cap L^2([0,T];D^{2,q}),\
(u_t, \theta_t)\in L^2([0,T];D^1_0),&\\
  (\sqrt{\rho} u_t, \sqrt{\rho} \theta_t)\in
L^\infty([0,T];L^2),\
\theta\ge0,\quad\quad\quad\quad\quad\quad\quad&\end{split}\ee and
$(\rho,u,\theta)$ satisfies (\ref{full N-S+1}) a.e. in
$\mathbb{R}^3\times (0,T]$. In particular, the strong solution
$(\rho, u, \theta)$ of (\ref{full N-S+1})-(\ref{3d-boundary}) is
called global strong solution, if the strong solution satisfies
(\ref{global strong}) for any $T>0$, and satisfies (\ref{full
N-S+1}) a.e. in $\mathbb{R}^3\times (0,\infty)$.
\end{definition}

\subsection{A blow-up criterion}\label{3d-sec 2.1}
 \setcounter{equation}{0}
\setcounter{theorem}{0}
\renewcommand{\theequation}{\thesubsection.\arabic{equation}}
\renewcommand{\thetheorem}{\thesubsection.\arabic{theorem}}

We state our main theorem in Section \ref{3d-sec 2.1}, which is on a
blow-up criterion for strong solutions to (\ref{full
N-S+1})-(\ref{3d-boundary}), as follows:
\begin{theorem}\label{blowup-th:1.1}
Assume $\rho_0\geq0$,
 $\rho_0\in H^1\cap W^{1,q}\cap L^1$, for some $q\in(3,6)$, $(u_0, \theta_0)\in D^2\cap D_0^1$, and the
following compatibility conditions are satisfied:
\beq\label{3d-compatibility}
\begin{cases}
\mu\Delta u_0+(\mu+\lambda)\nabla\mathrm{div} u_0-\nabla P(\rho_0, \theta_0)=\sqrt{\rho}_0g_1,\\
\kappa\Delta\theta_0+\frac{\mu}{2}\left|\nabla u_0+(\nabla
u_0)^\prime\right|^2+\lambda(\mathrm{div}u_0)^2=\sqrt{\rho_0}g_2,\
x\in\mathbb{R}^3,
\end{cases}
\eeq for some $g_i\in L^2$, $i=1,2$. Let $(\rho, u, \theta)$ be a
strong solution to (\ref{full N-S+1})-(\ref{3d-boundary}) in
$\mathbb{R}^3\times[0, T]$. If $0<T^*<+\infty$ is the maximal existence time of the strong solution, then \be\label{3d-result}
\lim\sup\limits_{T\nearrow T^*}\left(\|\rho\|_{L^\infty(0,T;
L^\infty)}+\|\rho\theta\|_{L^4(0,T; L^\frac{12}{5})}\right)=\infty,
\ee provided that $3\mu>\lambda$.
\end{theorem}

\begin{remark}
Theorem \ref{blowup-th:1.1} is an extension of our previous
result in \cite{Wen-Zhu 4} where a blow-up criterion in terms of
the upper bounds of $\rho$ and $\theta$ was established.
\end{remark}

\begin{remark}\label{3d-remark blowup}
The additional restriction on the viscosity, i.e., $3\mu>\lambda$, is only used in Lemma \ref{blowup-le:2.2}. Thus, this additional restriction $3\mu>\lambda$ can be removed if (\ref{3d-result}) is replaced by
\bex
\lim\sup\limits_{T\nearrow T^*}\left(\|\rho\|_{L^\infty(0,T;
L^\infty)}+\|\rho\theta\|_{L^4(0,T; L^\frac{12}{5})}+\|\rho^\frac{1}{4}u\|_{L^\infty(0,T;L^4)}+\big\| |u||\nabla
u| \big\|_{L^2(0,T;L^2)}\right)=\infty.
\eex
\end{remark}

\begin{remark}
Theorem \ref{blowup-th:1.1} gives a necessary condition that the
strong solution blows up in finite time. Thus, to prove the global
existence of strong solution to (\ref{full
N-S+1})-(\ref{3d-boundary}), it suffices to find some suitable
initial data and explore some global {\it a priori} estimates to
make the necessary condition fail. Please see Theorem
\ref{3d-th:1.1} and its proof in Section \ref{3d-sec 4} for more
details.
\end{remark}

\begin{remark}  Under the conditions of Theorem \ref{blowup-th:1.1},
the local existence and uniqueness of the strong solutions was
obtained by Cho and Kim in \cite{cho-Kim: perfect gas}. Thus, the
assumption $T^*>0$ makes sense.
\end{remark}
\subsection{Global strong solution}\label{3d-sec 2.2}
 \setcounter{equation}{0}
\setcounter{theorem}{0}
\renewcommand{\theequation}{\thesubsection.\arabic{equation}}

The second result is concerned with global existence and uniqueness
of strong solution.
\begin{theorem}\label{3d-th:1.1} (Global strong solution) For any given $K_i>0$ ($i=1,2$) and $\bar{\rho}>0$, assume that the initial data
$(\rho_0,u_0,\theta_0)$ satisfies \be\label{initial
data}\rho_0\geq0,\  \theta_0\ge0,\ in\ \mathbb{R}^3,\ \rho_0\in
H^1\cap W^{1,q}\cap L^1,\  (u_0, \theta_0)\in D^2\cap D_0^1, \ee for
some $q\in(3,6)$, and
 \be\label{3d-initial assumption}\begin{cases}
0\le\rho_0\le\bar{\rho},\ in\ \mathbb{R}^3,&\\
\|\nabla u_0\|_{L^2}\le \sqrt{K_1},\
\|\sqrt{\rho_0}\theta_0\|_{L^2}\le\sqrt{K_2},&\end{cases} \ee and
the compatibility conditions (\ref{3d-compatibility}). Then there
exists a unique global strong solution $(\rho,u,\theta)$ in
$\mathbb{R}^3\times[0,T]$ for any $T>0$, provided that
{\small\bex
\begin{split}m_0\le\varepsilon_0\triangleq\min\left\{C_3,\ \frac{\check{\mathcal
{C}}(2\mu+\lambda)^{6}}{\tilde{E}^3},\
\check{\mathcal
{C}}(2\mu+\lambda)^\frac{12}{7},\
\frac{\check{\mathcal
{C}}\mu^{12}(2\mu+\lambda)^{12}}{\tilde{E}^{12}},\ \check{\mathcal
{C}}(2\mu+\lambda)^{36}\kappa^{12}\right\},\end{split}\eex }
where $\displaystyle m_0=\int_{\mathbb{R}^3}\rho_0(x)\,dx$, $\tilde{E}=\frac{(14\mu+9\lambda)}{2\mu}+\frac{6\bar{\rho}K_2}{\mu(\mu+\lambda)K_1}
+\frac{8\bar{\rho}\kappa K_2}{\mu(\mu+\lambda)^2K_1}+1,$ and
{\small\bex\begin{split}C_3=\min\left\{\frac{\check{\mathcal
{C}}\kappa^6(\mu+\lambda)^6\mu^6}{\Big(\kappa(\mu+\lambda)+1\Big)^6},\
\check{\mathcal {C}}\mu^3\kappa^3(2\mu+\lambda)^6,\
\frac{\check{\mathcal {C}}\mu^6}{\tilde{E}^2},\ \check{\mathcal
{C}}\mu^2(2\mu+\lambda)^8, \frac{\check{\mathcal
{C}}\kappa^4\mu^2}{(2\mu+\lambda)^8\tilde{E}^6},\
\frac{\check{\mathcal {C}}\kappa^4}{\tilde{E}^4\mu^2},\
\frac{\check{\mathcal
{C}}\kappa^6}{\tilde{E}^3}\right\},\end{split}\eex } for some constant
$\check{\mathcal {C}}>0$ depending on $\bar{\rho}, K_1, K_2,$ and
some other known constants but independent of $\mu,\lambda,\kappa,$
and $t$.
\end{theorem}


\begin{remark}
For the fixed viscosity and heat conduction, we need the ``smallness" of initial mass. But the velocity and temperature can be large.
Even the density can also be large in some small regions.

\end{remark}

\begin{remark}
In Theorem \ref{3d-th:1.1}, for any fixed $\lambda$ satisfying
(\ref{viscosity assumption}), as $\mu$ and $\kappa$ are large
enough, we have
$$\tilde{E}\thicksim\frac{1}{\mu^2}
+\frac{\kappa}{\mu^3}+1,$$
\bex\begin{split}C_3\thicksim&\min\left\{\mu^6,\ \mu^9\kappa^3,\
\frac{\mu^6}{\tilde{E}^2},\ \mu^{10},
\frac{\kappa^4}{\mu^6\tilde{E}^6},\
\frac{\kappa^4}{\tilde{E}^4\mu^2},\
\frac{\kappa^6}{\tilde{E}^3}\right\}\\ \thicksim&\min\left\{\mu^6,\
\frac{\mu^{12}}{\kappa^2+\mu^6},\
\frac{\kappa^4\mu^{12}}{\kappa^6+\mu^{18}},\
\frac{\kappa^4\mu^{10}}{\kappa^4+\mu^{12}},\
\frac{\kappa^6\mu^9}{\kappa^3+\mu^9}\right\}.\end{split}\eex  Thus,
\bex\begin{split} \varepsilon_0\thicksim&\min\left\{
\frac{\mu^{12}}{\kappa^2+\mu^6},\
\frac{\kappa^4\mu^{12}}{\kappa^6+\mu^{18}},\
\frac{\kappa^4\mu^{10}}{\kappa^4+\mu^{12}},\
\frac{\kappa^6\mu^9}{\kappa^3+\mu^9},\
\frac{\mu^{15}}{\kappa^3+\mu^9},\
\mu^\frac{12}{7},\
\frac{\mu^{60}}{\kappa^{12}+\mu^{36}},\ \mu^{36}\kappa^{12}\right\}\\ \thicksim&\mu^{\alpha},
\end{split}\eex for some $\alpha=\alpha(r_1)>0$, provided that $\kappa\thicksim\mu^{r_1}$, for
any $r_1\in(\frac{3}{2},5)$. In this case, the initial mass (i.e.,
$C_0$) could be large if $\mu$ are sufficiently large. In fact, for isentropic flow (no temperature equation), there have been some works on the
global large regular solutions with vacuum for the initial-boundary-value problem in one dimension and in high dimensions (symmetric initial data), Cauchy problem and periodic problem in two dimensions \cite{Ding-Wen-Zhu,
Ding-Wen-Yao-Zhu, Jiu-Wang-XinJDE, Jiu-Wang-Xin}. For the full system, please refer to our previous works \cite{Wen-Zhu, Wen-Zhu 3} for the initial-boundary-value problem in one dimension and in high dimensions (symmetric initial data).
\end{remark}

\begin{remark}
For the coefficient of heat conduction $\kappa=0$, two important works by Xin et al. \cite{Xin, Xin-Yan} indicate that there are no global
smooth (classical) solutions to (\ref{full N-S+1}) with initial density of nontrivial compactly support or with initial density satisfying $$\begin{cases}V\subset\bar{V}\subset U\subset\Omega,&\\
\rho_0=0,\ \mathrm{in}\ U-V,& \end{cases}$$ and $\rho_0$ is not
identically equal to zero on $V$, where $U$ is a bounded and
connected open set.
\end{remark}

\begin{remark} For $\kappa>0$, with small initial mass, one can not expect generally that the global solutions as in Theorems \ref{3d-th:1.1} and \ref{3d-th1.3} are highly decreasing at infinity (in space) due to \cite{Rozanova} even if they are initially, or that the entropy $S$ has better regularity due to \cite{Tan-Wang}.
\end{remark}

\begin{remark}
We would like to mention the Ref. \cite{cho-jin} which shows that there is no global strong solution in $\mathbb{R}^3$ if the initial density is of nontrivial compactly support. While the authors in \cite{cho-jin} essentially impose on the strong solution an assumption for the entropy function $S$ (i.e., $S=S(x,t)<\infty$ for all $(x,t)\in\mathbb{R}^3\times[0,T]$) in order that the temperature function vanishes in the vacuum region. Thus here we find a global strong solution in a bigger space than that in \cite{cho-jin}. As a byproduct, when the initial mass is small in some sense, our result shows that the entropy function $S$ of the strong solution is not always less than infinity in $\mathbb{R}^3\times[0,T]$ even if it is initially.
\end{remark}

\subsection{Asymptotic behavior}\label{3d-sec 2.3}
 \setcounter{equation}{0}
\setcounter{theorem}{0}
\renewcommand{\theequation}{\thesubsection.\arabic{equation}}

The third result is concerned with the large-time behavior of the
solution as well as its decay rate.
\begin{theorem}\label{3d-th:1.2} (Asymptotic behavior)
Under the conditions of Theorem \ref{3d-th:1.1}, it holds that
\be\label{3d-thm longtime}
\int_{\mathbb{R}^3}\left(\rho|\theta|^2+|\nabla
u|^2+|\nabla\theta|^2\right)\rightarrow0, \ee as
$t\rightarrow\infty$, provided that
$$
m_0\le\varepsilon_0.
$$
Furthermore, the following decay rate  \be\label{3d-thmdecay}
\int_{\mathbb{R}^3}\left(\rho|\theta|^2+|\nabla u|^2\right)\le
\bar{C}\exp\{-\bar{C}_1t\}, \ee holds for any $t\in[0,\infty)$,
provided that
$$
m_0\le\min\{\varepsilon_0,\ \tilde{\varepsilon}_0\},
$$ for some positive constants $\tilde{\varepsilon}_0$, $\bar{C}$, and $\bar{C}_1$ depending on $\mu,\lambda,\kappa, K_1, K_2,
\bar{\rho}$, and some other known constants but independent of $t$.
\end{theorem}
\begin{remark}
Some large-time behavior of the solutions to Cauchy problem for
(\ref{3d-full N-S}) with non-vacuum state at infinity have been
studied before, see for instance \cite{Hoff-arma, Huang-Li} and
references therein. While, there seem few results on the decay rate.
Here we get the exponential decay estimate (\ref{3d-thmdecay}) which
seems surprising. The main ingredient is that here we have the
integrability and the uniform upper bound of $\rho$ in
$\mathbb{R}^3$ such that the inequality obtained from
(\ref{3d-asytheta u+1}) (assume all the constant coefficients are 1
w.l.o.g.) \bex\begin{split} \int_{\mathbb{R}^3}\left(|\nabla
u|^2+|\nabla\theta|^2+\rho
|\u|^2\right)+\frac{d}{dt}\int_{\mathbb{R}^3}\left(\rho|u|^2+\rho|\theta|^2+|\mathrm{curl}
u|^2+G^2\right) \le0
\end{split}
\eex is able to give the desired inequality (assume the constant coefficients are 1) which implies the exponential decay estimate, i.e.,
\bex\begin{split} \int_{\mathbb{R}^3}\left(\rho|u|^2+\rho|\theta|^2+|\mathrm{curl}
u|^2+G^2\right)+\frac{d}{dt}\int_{\mathbb{R}^3}\left(\rho|u|^2+\rho|\theta|^2+|\mathrm{curl}
u|^2+G^2\right) \le0.
\end{split}
\eex Refer to Section \ref{3d-sec 5.2} for more details.
\end{remark}

\begin{remark}
From (\ref{full N-S+1})$_3$, Lemmas \ref{3d-le:low uniform}, \ref{3d-le: H 2 of u and H 1of theta}, \ref{3d-le: H 2 of theta}, Corollary \ref{corollary-long time}, and the standard interpolation inequality, it is easy to verify that the large-time behavior (\ref{3d-thm longtime}) and the decay estimate (\ref{3d-thmdecay}) can be improved as follows:
\bex
\int_{\mathbb{R}^3}\left(\rho|\theta|^{r_1}+|\nabla
u|^r+|\nabla\theta|^r\right)\rightarrow0 \eex for any $r\in[2,6)$ and any $r_1\in[2,\infty)$, as
$t\rightarrow\infty$,
and
\bex
\int_{\mathbb{R}^3}\left(\rho|\theta|^{r_1}+|\nabla u|^r\right)\le
\bar{\bar{C}}\exp\{-\bar{\bar{C}}_1t\} \eex  for any $t\in[0,\infty)$ and some positive constants $\bar{\bar{C}}$ and $\bar{\bar{C}}_1$ independent of t.

\end{remark}

\subsection{Global classical solution}\label{3d-sec 2.5}
 \setcounter{equation}{0}
\setcounter{theorem}{0}
\renewcommand{\theequation}{\thesubsection.\arabic{equation}}
The global strong solution as in Theorem \ref{3d-th:1.1} can become
classical, if there is more regularity of the initial density(c.f.
\cite{Huang-Li, Jiu-Wang-Xin}). Before presenting the main result in
this part, we would like to give a definition of classical solution
in the paper.
\begin{definition}(classical solution) For $T>0$, $(\rho, u, \theta)$ is called a classical solution to the compressible Navier-Stokes equations
(\ref{full N-S+1})-(\ref{3d-boundary}) in $\mathbb{R}^3\times
[0,T]$, if for some $q\in (3,
6)$,\be\label{global-regularities}\begin{split} &\rho\in C([0,T];
H^2\cap W^{2,q}),\ \rho_t\in C([0,T]; H^1),\ \rho\ge0,\
\theta\ge0,&\\& (u,\theta)\in C([0,T];D^2\cap D^1_0)\cap
L^2([0,T];D^3),\ (u_t,\theta_t)\in L^2([0,T]; D_0^1),&\\&
  (\sqrt{\rho} u_t, \sqrt{\rho} \theta_t)\in
L^\infty([0,T];L^2),\ \sqrt{t}\sqrt{\rho} u_{tt}\in L^2([0,T];
L^2),\ t\sqrt{\rho} u_{tt}\in L^\infty([0,T];L^2),&\\& \sqrt{t}u\in
L^\infty([0,T];D^3),\ \sqrt{t} u_t\in L^\infty([0,T];D_0^1)\cap
L^2([0,T]; D^2),&\\& tu\in L^\infty([0,T];D^{3,q}),\ tu_t\in
L^\infty([0,T]; D^2),\ tu_{tt}\in L^2([0,T];D_0^1),&\\& t\theta\in
L^\infty([0,T];D^3)\cap L^2([0,T];D^4),\ t\theta_t\in
L^\infty([0,T]; D_0^1)\cap L^2([0,T];D^2),&\\&
t^\frac{3}{2}\theta\in L^\infty([0,T]; D^4),\
t^\frac{3}{2}\theta_t\in L^\infty([0,T]; D^2),\
t\sqrt{\rho}\theta_{tt}\in
L^2([0,T];L^2),&\\&t^\frac{3}{2}\sqrt{\rho}\theta_{tt}\in
L^\infty([0,T];L^2),\  t^\frac{3}{2}\theta_{tt}\in L^2([0,T];D_0^1),
\end{split}
\ee and $(\rho,u,\theta)$ satisfies (\ref{full N-S+1}) in
$\mathbb{R}^3\times (0,T]$. In particular, the classical solution
$(\rho, u, \theta)$ of (\ref{full N-S+1})-(\ref{3d-boundary}) is
called global classical solution, if the classical solution
satisfies (\ref{global-regularities}) for any $T>0$, and satisfies
(\ref{full N-S+1}) in $\mathbb{R}^3\times (0,\infty)$.
\end{definition}
Now we are in a position to state our main result in this part.
\begin{theorem}\label{3d-th1.3} (Global classical solution) Under
the conditions of Theorem \ref{3d-th:1.1}, if in addition that
 \be\label{glocal-initial data}
\rho_0\in H^2\cap W^{2,q}, \ee for some $q\in(3,6)$, then there
exists a unique global classical solution $(\rho,u,\theta)$ of
(\ref{full N-S+1})-(\ref{3d-boundary}).
\end{theorem}
\begin{remark}
Though the initial data can be large if the coefficients of
viscosity and heat conduction are large, it is still unknown whether
the global classical solution exists when the initial data are large
for the fixed coefficients of viscosity and heat conduction. It
should be noted that a similar question of whether there exists a
global smooth solution of the three-dimensional incompressible
Navier-Stokes equations with smooth initial data is one of the most
outstanding mathematical open problems (\cite{fefferman}).
\end{remark}

\section{A blow-up criterion}\label{3d-sec 3}

 \setcounter{equation}{0}\setcounter{theorem}{0}
\renewcommand{\theequation}{\thesection.\arabic{equation}}
\renewcommand{\thetheorem}{\thesection.\arabic{theorem}}

Let $0<T^*<\infty$ be the maximal existence time of the strong
solution $(\rho, u, \theta)$ to (\ref{full
N-S+1})-(\ref{3d-boundary}). Namely, $(\rho, u, \theta)$ is a strong
solution to (\ref{full N-S+1})-(\ref{3d-boundary}) in
$\mathbb{R}^3\times [0, T]$ for any $0<T<T^*$, but not a strong
solution in $\mathbb{R}^3\times [0, T^*]$. We shall prove Theorem
\ref{blowup-th:1.1} by using a contradiction argument. Suppose that
(\ref{3d-result}) is false, i.e. \beq\label{blowup-2.1}
M:=\lim\sup\limits_{t\nearrow
T^*}(\|\rho(t)\|_{L^\infty}+\int_0^t\|\rho\theta(s)\|_{L^\frac{12}{5}}^4ds)<\infty.
\eeq Our aim is to show that under the assumption
(\ref{blowup-2.1}), there is a bound $C>0$ depending only on $M,
\rho_0, u_0, \theta_0, \mu,\lambda, \kappa$, and $T^*$ such that
\beq\label{non-uniform_est1} \sup_{0\le
t<T^*}\|\theta(t)\|_{L^\infty}\le C. \eeq With
(\ref{non-uniform_est1}) and (\ref{blowup-2.1}), we showed in our
previous paper \cite{Wen-Zhu 4} that $T^*$ is not the maximal time,
which is the desired contradiction.

Throughout the rest of the section, we denote by $C$ a generic
constant depending only on $\rho_0$, $u_0$, $\theta_0$, $T^*$, $M$,
$\lambda$, $\mu$, $\kappa$.
\begin{lemma}\label{3d-le:3.1}
Under the conditions of Theorem \ref{blowup-th:1.1}, it holds that
\be\label{3d-mass conservation} \int_{\mathbb{R}^3}
\rho=\int_{\mathbb{R}^3}\rho_0\triangleq m_0,\ee for any
$t\in[0,T^*)$.
\end{lemma}
\pf For any $r>1$, let $\phi_r$ be the classical cut-off function
satisfying
$$
\phi_r\in C_0^\infty(\mathbb{R}^3),\ \mathrm{supp}\phi_r\in
B_{r}(0),\  \phi_r\equiv1\ \mathrm{in}\ B_{\frac{r}{2}}(0),\
0\le\phi_r\le1,\ \mathrm{and}\ |\nabla \phi_r|\le\frac{C}{r}\
\mathrm{in}\ \mathbb{R}^3.
$$
Multiplying (\ref{full N-S+1})$_1$ by $\phi_r$, we have
\be\label{new mass equation} (\phi_r\rho)_t+\phi_r\nabla\cdot(\rho
u)=0. \ee Integrating (\ref{new mass equation}) over
$\mathbb{R}^3\times[0,t]$ for $0\le t< T^*$, and using integration
by parts, we have \be\label{rho phi r}\begin{split}
\int_{\mathbb{R}^3}\phi_r\rho=&\int_{\mathbb{R}^3}\phi_r\rho_0+\int_0^t\int_{\mathbb{R}^3}\rho
u\cdot\nabla\phi_r\\
=&\int_{\mathbb{R}^3}\phi_r\rho_0+\int_0^t\int_{B_r(0)/B_{\frac{r}{2}}(0)}\rho
u\cdot\nabla\phi_r\\=&I_1+I_2.
\end{split}
\ee Since $\rho_0\in L^1$, we have \be\label{3d-I 1}
I_1\rightarrow\int_{\mathbb{R}^3}\rho_0 \ee as $r\rightarrow\infty$.

For $I_2$, we have \bex\begin{split}
|I_2|\le&C\int_0^t\int_{B_r(0)/B_{\frac{r}{2}}(0)}\rho \frac{|u|}{r}\\
\le&C\int_0^t\big\|u\big\|_{L^6(\mathbb{R}^3)}\big\|\rho\big\|_{L^2\left(B_r(0)/B_{\frac{r}{2}}(0)\right)}
\big\|\frac{1}{r}\big\|_{L^3\left(B_r(0)/B_{\frac{r}{2}}(0)\right)}\\
\le&C\int_0^t\big\|\nabla
u\big\|_{L^2(\mathbb{R}^3)}\big\|\rho\big\|_{L^2\left(B_r(0)/B_{\frac{r}{2}}(0)\right)}.
\end{split}
\eex Since $\nabla u\in C([0,t];L^2(\mathbb{R}^3))$ and $\rho\in
C([0,t];L^2(\mathbb{R}^3))$ for $t\in[0,T^*)$, let $r$ go to
$\infty$, one has
$$ I_2\rightarrow0.$$ This together with (\ref{rho phi r}) and
(\ref{3d-I 1}) deduces
$$
\rho(\cdot,t)\in L^1,
$$  and
$$
\int_{\mathbb{R}^3}\rho=\int_{\mathbb{R}^3}\rho_0
$$ for any $t\in[0,T]$.
\endpf
\begin{lemma}\label{blowup-le:2.2}
Under the conditions of Theorem \ref{blowup-th:1.1} and
(\ref{blowup-2.1}), if $3\mu>\lambda$, it holds that \be\label{u
2nabla u 2} \sup\limits_{0\leq t\leq
T}\int_{\mathbb{R}^3}\rho|u|^4\,dx+\int_0^T\int_{\mathbb{R}^3}|u|^2|\nabla
u|^2\,dx\,dt\leq C, \ee for any $T\in(0,T^*)$.
\end{lemma}
\pf  The detailed proof of Lemma \ref{blowup-le:2.2} could be found
in \cite{Wen-Zhu 4} (Lemma 4.2 therein), which might be slightly
modified (only for the pressure term).
\endpf
\begin{lemma}\label{blowup-le: int nabla u}Under the conditions of Theorem \ref{blowup-th:1.1} and (\ref{blowup-2.1}), it holds that for any $T\in[0,T^*)$
\be\label{H 1 of u}\sup\limits_{0\le t\le
T}\int_{\mathbb{R}^3}(\rho|\theta|^2+|\nabla u|^2)\,
dx+\int_0^T\int_{\mathbb{R}^3}(|\nabla\theta|^2+\rho
|u_t|^2)\,dxdt\le C.\ee
\end{lemma}
\pf Multiplying (\ref{full N-S+1})$_2$ by $u_t$, and integrating by
parts over $\mathbb{R}^3$, we have

\beq\label{dt nabla u 2-1}\begin{split} &\int_{\mathbb{R}^3}\rho
|u_t|^2+\frac{1}{2}\frac{d}{dt}\int_{\mathbb{R}^3}\left(\mu|\nabla
u|^2+(\mu+\lambda)|\mathrm{div}u|^2\right)\\=&\frac{d}{dt}\int_{\mathbb{R}^3}P\mathrm{div}u-\frac{1}{2(2\mu+\lambda)}
\frac{d}{dt}\int_{\mathbb{R}^3}P^2
-\frac{1}{2\mu+\lambda}\int_{\mathbb{R}^3}P_tG-\int_{\mathbb{R}^3}\rho
u\cdot\nabla u \cdot u_t\\=&\sum\limits_{i=1}^4(I)_i,
\end{split}
\eeq where $G=(2\mu+\lambda)\mathrm{div}u-P$.

Recalling $P=\rho\theta$, we obtain from (\ref{full N-S+1})$_1$ and
(\ref{full N-S+1})$_3$ \beq\label{equation of P t}\begin{split}
P_t=-\mathrm{div}(P u)-\rho\theta\mathrm{div}u+\mu\left(\nabla
u+(\nabla u)^\prime\right):\nabla
u+\lambda\mathrm{div}u\mathrm{div}u+\kappa\Delta\theta.
\end{split}
\eeq Substituting (\ref{equation of P t}) into $(I)_3$, and using
integration by parts and the H\"older inequality, we have
 \beq\label{blowup-(I) 3}\begin{split}
(I)_3=&-\frac{1}{2\mu+\lambda}\int_{\mathbb{R}^3}P u\cdot\nabla
G+\frac{1}{2\mu+\lambda}\int_{\mathbb{R}^3}\rho\theta\mathrm{div}u
G\\&+\frac{\mu}{2\mu+\lambda}\int_{\mathbb{R}^3}\left(\nabla
u+(\nabla u)^\prime\right):\left(\nabla G\otimes
u\right)+\frac{\lambda}{2\mu+\lambda}\int_{\mathbb{R}^3}\mathrm{div}u
u\cdot\nabla
G\\&+\frac{1}{2\mu+\lambda}\int_{\mathbb{R}^3}\big(\mu\Delta
u+(\mu+\lambda)\nabla \mathrm{div}u\big)\cdot uG
+\frac{\kappa}{2\mu+\lambda}\int_{\mathbb{R}^3}\nabla\theta\cdot\nabla G\\
\le& C\|\rho u\theta\|_{L^2}\|\nabla
G\|_{L^2}+\frac{1}{2\mu+\lambda}\int_{\mathbb{R}^3}\rho\theta\mathrm{div}u
G+C\|\nabla G\|_{L^2}\big\|u|\nabla u|\big\|_{L^2}\\&+C\|\nabla
G\|_{L^2}\|\nabla\theta\|_{L^2}+\frac{1}{2\mu+\lambda}\int_{\mathbb{R}^3}\big(\mu\Delta
u+(\mu+\lambda)\nabla \mathrm{div}u\big)\cdot uG.
\end{split}
\eeq Substituting (\ref{full N-S+1})$_2$ into (\ref{blowup-(I) 3}),
and using the Sobolev inequality, (\ref{blowup-2.1}) and integration
by parts, we have \beq\label{blowup-(I) 3+1}\begin{split} (I)_3\le&
C\|\nabla G\|_{L^2}\Big(\|\rho u\theta\|_{L^2}+\big\|u|\nabla
u|\big\|_{L^2}+\|\nabla\theta\|_{L^2}\Big)+\frac{1}{2\mu+\lambda}\int_{\mathbb{R}^3}\rho
u_t\cdot uG\\&+\frac{1}{2\mu+\lambda}\int_{\mathbb{R}^3}\rho
u\cdot\nabla u\cdot uG-\frac{1}{2\mu+\lambda}\int_{\mathbb{R}^3} Pu\cdot\nabla G\\
\le&C\|\nabla G\|_{L^2}\Big(\|\rho u\theta\|_{L^2}+\big\|u|\nabla
u|\big\|_{L^2}+\|\nabla\theta\|_{L^2}\Big)+\frac{1}{6}\int_{\mathbb{R}^3}\rho
|u_t|^2\\&+C\int_{\mathbb{R}^3}\rho|u|^2|G|^2+C\big\|u|\nabla u|\big\|_{L^2}^2\\
\le&C\|\nabla G\|_{L^2}\Big(\|\rho u\theta\|_{L^2}+\big\|u|\nabla
u|\big\|_{L^2}+\|\nabla\theta\|_{L^2}\Big)+\frac{1}{6}\int_{\mathbb{R}^3}\rho
|u_t|^2\\&+C\int_{\mathbb{R}^3}|u|^2|\nabla
u|^2+C\int_{\mathbb{R}^3}\rho|u|^2|\rho\theta|^2.
\end{split}
\eeq
 Taking $div$ on both side of
(\ref{full N-S+1})$_2$, we get \be\label{equation of G} \Delta
G=\mathrm{div}(\rho u_t+\rho u\cdot\nabla u). \ee
 By (\ref{equation of G}) and the standard
$L^2$-estimates together with (\ref{blowup-2.1}), we get
\beq\label{H 1 of G}
\begin{split}\|\nabla G\|_{L^{2}}\le C\|\rho u_t\|_{L^{2}}+C\|\rho u\cdot\nabla u\|_{L^2}\le
 C\|\sqrt{\rho} u_t\|_{L^{2}}+C\big\| |u||\nabla u|\big\|_{L^2}.
\end{split}
\eeq
 Substituting (\ref{H 1 of G}) into (\ref{blowup-(I) 3+1}), and using
the Cauchy inequality, we have \beq\label{blowup-(I)
3+2-1}\begin{split} (I)_3 \le&C\|\rho
u\theta\|_{L^2}^2+C\big\|u|\nabla
u|\big\|_{L^2}^2+C\|\nabla\theta\|_{L^2}^2+\frac{1}{3}\int_{\mathbb{R}^3}\rho
|u_t|^2.
\end{split}
\eeq For the first term of the right hand side of (\ref{blowup-(I)
3+2-1}), using the H\"older inequality, the Sobolev inequality and
the Cauchy inequality, we have \be\label{rho u theta}\begin{split}
\|\rho
u\theta\|_{L^2}^2\le&\big\||u|^2\big\|_{L^6}\|\rho\theta\|_{L^\frac{12}{5}}^2\le
C\big\|u|\nabla u|\big\|_{L^2}\|\rho\theta\|_{L^\frac{12}{5}}^2\\
\le&C\big\|u|\nabla
u|\big\|_{L^2}^2+C\|\rho\theta\|_{L^\frac{12}{5}}^4.
\end{split}
\ee Substituting (\ref{rho u theta}) into (\ref{blowup-(I) 3+2-1}),
we have
 \beq\label{blowup-(I) 3+2}\begin{split} (I)_3
\le&C\|\rho\theta\|_{L^\frac{12}{5}}^4+C\big\|u|\nabla
u|\big\|_{L^2}^2+C\|\nabla\theta\|_{L^2}^2+\frac{1}{3}\int_{\mathbb{R}^3}\rho
|u_t|^2.
\end{split}
\eeq
 For $(I)_4$, using Cauchy inequality and (\ref{blowup-2.1}), we have
 \be\label{blowup-(I) 4}\begin{split} (I)_4\le&\frac{1}{6}\int_{\mathbb{R}^3}\rho
|u_t|^2+C\int_{\mathbb{R}^3} |u|^2 |\nabla u|^2.
\end{split}
\ee Putting (\ref{blowup-(I) 3+2}) and (\ref{blowup-(I) 4}) into
(\ref{dt nabla u 2-1}), and integrating it over $[0,t]$, for
$t<T^*$, we have \bex\begin{split} &\int_0^t\int_{\mathbb{R}^3}\rho
|u_t|^2+\int_{\mathbb{R}^3}\left(\mu|\nabla
u|^2+(\mu+\lambda)|\mathrm{div}u|^2\right)\\
\le&2\int_{\mathbb{R}^3}P\mathrm{div}u+C\int_0^t\|\nabla\theta\|_{L^2}^2+C\\
\le&(\mu+\lambda)\int_{\mathbb{R}^3}|\mathrm{div}u|^2+C\big(\int_{\mathbb{R}^3}\rho\theta^2+\int_0^t\int_{\mathbb{R}^3}|\nabla\theta|^2\big)+C,
\end{split}
\eex where we have used Cauchy inequality, (\ref{blowup-2.1}) and
(\ref{u 2nabla u 2}). Therefore, \be\label{zy}\begin{split}
\int_0^t\int_{\mathbb{R}^3}\rho|u_t|^2+\int_{\mathbb{R}^3}|\nabla
u|^2\leq
C\big(\int_{\mathbb{R}^3}\rho\theta^2+\int_0^t\int_{\mathbb{R}^3}|\nabla\theta|^2\big)+C.
\end{split}
\ee Multiplying (\ref{full N-S+1})$_3$ by $\theta$ and integrating
by parts over $\mathbb{R}^3$, we have \beq\label{zyy}\begin{split}
&\kappa\int_{\mathbb{R}^3}|\nabla\theta|^2+\frac{1}{2}\frac{d}{dt}\int_{\mathbb{R}^3}\rho|\theta|^2\\
=&-\int_{\mathbb{R}^3}\rho\theta^2\mathrm{div}u+\int_{\mathbb{R}^3}\frac{\mu}{2}|\nabla u+(\nabla u)'|^2\theta+\int_{\mathbb{R}^3}\lambda|\mathrm{div}u|^2\theta\\
=&\sum\limits_{i=1}^3(II)_i.
\end{split}
\eeq  For $(II)_2$ and $(II)_3$, we have
\beq\label{blowup-(II)_2,3}\begin{split}
(II)_2+(II)_3=&\int_{\mathbb{R}^3}\mu\left(\nabla u+(\nabla u)'\right):\nabla u\theta+\int_{\mathbb{R}^3}\lambda|\mathrm{div}u|^2\theta\\
=&-\int_{\mathbb{R}^3}\mu(\triangle u+\nabla\mathrm{div}u )\cdot
u\theta-\int_{\mathbb{R}^3}\mu\left(\nabla u+(\nabla
u)'\right):(\nabla\theta\otimes u)
\\&-\int_{\mathbb{R}^3}\lambda u\cdot\nabla\mathrm{div}u\theta-\int_{\mathbb{R}^3}\lambda\mathrm{div}u u\cdot\nabla \theta\\
=&-\int_{\mathbb{R}^3}(\rho u_t+\rho u\cdot \nabla u+\nabla P)\cdot
u\theta-\int_{\mathbb{R}^3}\mu\left(\nabla u+(\nabla
u)'\right):(\nabla\theta\otimes u)
\\&-\int_{\mathbb{R}^3}\lambda\mathrm{div}u u\cdot\nabla \theta\\
=&-\int_{\mathbb{R}^3}\rho u_t\cdot
u\theta-\int_{\mathbb{R}^3}\rho(u\cdot\nabla)u\cdot
u\theta+\int_{\mathbb{R}^3}\rho\theta^2\mathrm{div}u+\int_{\mathbb{R}^3}\rho\theta
u\cdot\nabla\theta
\\&-\int_{\mathbb{R}^3}\mu\left(\nabla u+(\nabla u)'\right):(\nabla\theta\otimes u)-\int_{\mathbb{R}^3}\lambda\mathrm{div}u u\cdot\nabla \theta,
\end{split}
\eeq where we have used integration by parts and (\ref{full N-S+1})$_2$.\\
Using the H\"older inequality, the Cauchy inequality,
(\ref{blowup-2.1}) and (\ref{u 2nabla u 2}), we have
\beq\label{blowup-(II) 2,3+1}\begin{split} &(II)_2+(II)_3\\
\leq&
\|\sqrt{\rho}u_t\|_{L^2}\|\sqrt[4]{\rho}u\|_{L^4}\|\sqrt[4]{\rho}\theta\|_{L^4}
+ \big\|u|\nabla u|\big\|_{L^2}\|\rho u\|_{L^3}\|\theta\|_{L^6}
+\int_{\mathbb{R}^3}\rho\theta^2\mathrm{div}u\\&+\|\nabla\theta\|_{L^2}\|\sqrt{\rho}u\|_{L^4}\|\sqrt{\rho}\theta\|_{L^4}
+C\big\|u|\nabla u|\big\|_{L^2}\|\nabla\theta\|_{L^2}\\
\leq&C\|\sqrt{\rho}u_t\|_{L^2}\|\sqrt[4]{\rho}\theta\|_{L^4}+\frac{\kappa}{2}\|\nabla\theta\|_{L^2}^2+C\big\|u|\nabla
u|\big\|_{L^2}^2+\int_{\mathbb{R}^3}\rho\theta^2\mathrm{div}u+C\|\sqrt{\rho}\theta\|_{L^4}^2.
\end{split}
\eeq Substituting (\ref{blowup-(II) 2,3+1}) into (\ref{zyy}), we
have \beq\label{zyyy}\begin{split}
\kappa\int_{\mathbb{R}^3}|\nabla\theta|^2+\frac{d}{dt}\int_{\mathbb{R}^3}\rho|\theta|^2\leq
C\|\sqrt{\rho}u_t\|_{L^2}\|\sqrt[4]{\rho}\theta\|_{L^4}+C\big\|u|\nabla
u|\big\|_{L^2}^2+C\|\sqrt{\rho}\theta\|_{L^4}^2.
\end{split}
\eeq Integrating (\ref{zyyy}) over $[0,t]$ ($t<T^*$), and using
(\ref{u 2nabla u 2}), we have \beq\label{sy}\begin{split}
\int_0^t\int_{\mathbb{R}^3}|\nabla\theta|^2+\int_{\mathbb{R}^3}\rho|\theta|^2\leq
C\int_0^t\|\sqrt{\rho}u_t\|_{L^2}\|\sqrt[4]{\rho}\theta\|_{L^4}+C\int_0^t\|\sqrt{\rho}\theta\|_{L^4}^2+C.
\end{split}
\eeq Multiplying (\ref{sy}) by $2C$, and adding the resulting
inequality into (\ref{zy}), we have \bex\begin{split}
&\int_0^t\int_{\mathbb{R}^3}\rho|u_t|^2+\int_{\mathbb{R}^3}|\nabla u|^2+\int_0^t\int_{\mathbb{R}^3}|\nabla\theta|^2+\int_{\mathbb{R}^3}\rho|\theta|^2\\
\le&
C\int_0^t\|\sqrt{\rho}u_t\|_{L^2}\|\sqrt[4]{\rho}\theta\|_{L^4}+C\int_0^t\|\sqrt{\rho}\theta\|_{L^4}^2+C\\
\le&
\frac{1}{2}\int_0^t\int_{\mathbb{R}^3}\rho|u_t|^2+C\int_0^t\left(\int_{\mathbb{R}^3}\sqrt{\rho}|\theta|\sqrt{\rho}|\theta|^3\right)^\frac{1}{2}+C\\
\le&\frac{1}{2}\int_0^t\int_{\mathbb{R}^3}\rho|u_t|^2+C\int_0^t\|\sqrt{\rho}\theta\|_{L^2}^\frac{1}{2}\|\sqrt[6]{\rho}\theta\|_{L^6}^\frac{3}{2}+C\\
\le&\frac{1}{2}\int_0^t\int_{\mathbb{R}^3}\rho|u_t|^2+\frac{1}{2}\int_0^t\int_{\mathbb{R}^3}|\nabla\theta|^2+C\int_0^t\int_{\mathbb{R}^3}\rho|\theta|^2+C,
\end{split}
\eex where we have used the Young inequality, the H\"older
inequality, the Sobolev inequality and (\ref{blowup-2.1}).
 This, together with the Gronwall inequality, gives (\ref{H 1 of
u}).
\endpf

\begin{lemma}\label{blowup-le: int rho u t}Under the conditions of Theorem \ref{blowup-th:1.1} and (\ref{blowup-2.1}), it holds that for any $t\in(0,T^*)$
\be\label{H 2 of u}\int_{\mathbb{R}^3}(|\nabla
\theta|^2+\rho|\u|^2)\,dx+\int_0^t\int_{\mathbb{R}^3}(\rho
|\dot{\theta}|^2+|\nabla\u|^2)\,dx\,ds\le C.\ee
\end{lemma}
\pf By the definition of $\dot{u}$, (\ref{full N-S+1})$_2$ can be
reformulated as follows: \be\label{N-1}\rho\dot{u}+\nabla
P=\mu\Delta u+(\mu+\lambda)\nabla\mathrm{div}u.\ee Denote
$$
f_\epsilon\doteq f_{\epsilon}(x,\cdot)=\eta_\epsilon\star
f(x,\cdot),\quad f_\delta\doteq f_{\delta}(\cdot,t)=\phi_\delta\star
f(\cdot,t),\quad
f_{\epsilon,\delta}=\phi_\delta\star(\eta_\epsilon\star f),
$$
where
$$
\eta_\epsilon(\cdot)=\frac{1}{\epsilon}\eta(\frac{\cdot}{\epsilon}),\quad
\phi_\delta(\cdot)=\frac{1}{\delta^3}\phi(\frac{\cdot}{\delta}).
$$ Here $\eta$ and $\phi$ are the standard mollifiers in one
dimension and in three dimensions respectively. For any given
$\tau>0$, let $\epsilon\in(0, \tau]$. Taking convolutions of both
sides of (\ref{N-1}) with $\eta$ and $\phi$, we have
\beq\label{N-0}(\rho\dot{u})_{\epsilon,\delta}+\nabla
P_{\epsilon,\delta}=\mu\Delta
u_{\epsilon,\delta}+(\mu+\lambda)\nabla\mathrm{div}u_{\epsilon,\delta}\eeq
in $\mathbb{R}^3\times(\tau,T-\tau)$ where $T<T^*$.

Differentiating (\ref{N-0}) with respect to $t$, we have
\be\label{N1}
\begin{split}
\frac{\partial}{\partial t}(\rho\dot{u})_{\epsilon,\delta}+\nabla
(P_t)_{\epsilon,\delta}
=&\mu\Delta(\u)_{\epsilon,\delta}+(\mu+\lambda)\nabla\mathrm{div}(\u)_{\epsilon,\delta}-\mu\Delta[(u\cdot\nabla)
u]_{\epsilon,\delta}\\&-(\mu+\lambda)\nabla\mathrm{div}[(u\cdot\nabla)
u]_{\epsilon,\delta}.
\end{split}
\ee Multiplying (\ref{N1}) by $(\u)_{\epsilon,\delta}$, and
integrating by parts over $\mathbb{R}^3$, we have 
\bex
\begin{split}
&\int_{\mathbb{R}^3}\frac{\partial}{\partial
t}(\rho\dot{u})_{\epsilon,\delta}(\u)_{\epsilon,\delta}+\int_{\mathbb{R}^3}\left[\mu|\nabla(\u)_{\epsilon,\delta}|^2
+(\mu+\lambda)|\di(\u)_{\epsilon,\delta}|^2\right]\\
=&\int_{\mathbb{R}^3}
(P_t)_{\epsilon,\delta}\di(\u)_{\epsilon,\delta}
+\mu\int_{\mathbb{R}^3}\nabla[(u\cdot\nabla) u]_{\epsilon,\delta}
\cdot\nabla(\u)_{\epsilon,\delta} \\
&+(\mu+\lambda) \int_{\mathbb{R}^3}\di[(u\cdot\nabla)
u]_{\epsilon,\delta}\di(\u)_{\epsilon,\delta}.
\end{split}
\eex
 Let $\delta\rightarrow0^+$, we have \beq\label{N2}
\begin{split}
&\int_{\mathbb{R}^3}\frac{\partial}{\partial
t}(\rho\dot{u})_{\epsilon}(\u)_{\epsilon}+\int_{\mathbb{R}^3}\left[\mu|\nabla(\u)_{\epsilon}|^2
+(\mu+\lambda)|\di(\u)_{\epsilon}|^2\right]\\
=&\int_{\mathbb{R}^3} (P_t)_{\epsilon}\di(\u)_{\epsilon}
+\mu\int_{\mathbb{R}^3}\nabla[(u\cdot\nabla) u]_{\epsilon}
\cdot\nabla(\u)_{\epsilon} \\
&+(\mu+\lambda) \int_{\mathbb{R}^3}\di[(u\cdot\nabla)
u]_{\epsilon}\di(\u)_{\epsilon}.
\end{split}
\eeq Note that \bex\begin{split} \frac{\partial}{\partial
t}(\rho\dot{u})_{\epsilon}=&\frac{1}{\epsilon^2}\int_0^T\eta^\prime(\frac{t-s}{\epsilon})\rho(\cdot,s)\dot{u}(\cdot,s)\,ds\\
=&\frac{1}{\epsilon}\int_{-1}^1\eta^\prime(s)\rho(\cdot,t-\epsilon
s)\dot{u}(\cdot,t-\epsilon s)\,ds\\
=&\frac{1}{\epsilon}\int_{-1}^1\eta^\prime(s)[\rho(\cdot,t-\epsilon
s)-\rho(\cdot,t)]\dot{u}(\cdot,t-\epsilon
s)\,ds+\rho(\cdot,t)\frac{1}{\epsilon}\int_{-1}^1\eta^\prime(s)\dot{u}(\cdot,t-\epsilon
s)\,ds\\
=&\frac{1}{\epsilon}\int_{-1}^1\eta^\prime(s)[\rho(\cdot,t-\epsilon
s)-\rho(\cdot,t)][\dot{u}(\cdot,t-\epsilon
s)-\dot{u}(\cdot,t)]\,ds\\&+\frac{1}{\epsilon}\dot{u}(\cdot,t)\int_{-1}^1\eta^\prime(s)[\rho(\cdot,t-\epsilon
s)-\rho(\cdot,t)]\,ds+\rho(\cdot,t)\frac{\partial}{\partial
t}(\dot{u})_\epsilon
\end{split}
\eex Thus, \be\label{N3}\begin{split}
 \int_{\mathbb{R}^3}\frac{\partial}{\partial
t}(\rho\dot{u})_{\epsilon}(\u)_{\epsilon}\,dx=&\int_{\mathbb{R}^3}(\u)_{\epsilon}(x,t)\int_{-1}^1\eta^\prime(s)\frac{\rho(x,t-\epsilon
s)-\rho(x,t)}{\epsilon}[\dot{u}(x,t-\epsilon
s)-\dot{u}(x,t)]\,ds\,dx\\&+\int_{\mathbb{R}^3}(\u)_{\epsilon}(x,t)\dot{u}(x,t)\int_{-1}^1\eta^\prime(s)\frac{\rho(x,t-\epsilon
s)}{\epsilon}\,ds\,dx\\&+
\frac{1}{2}\frac{d}{dt}\int_{\mathbb{R}^3}\rho|(\dot{u})_\epsilon|^2\,dx-\frac{1}{2}\int_{\mathbb{R}^3}\rho_t|(\dot{u})_\epsilon|^2\,dx\\=&
G_\epsilon(t)+\int_{\mathbb{R}^3}(\u)_{\epsilon}\dot{u}(\rho_t)_\epsilon\,dx+
\frac{1}{2}\frac{d}{dt}\int_{\mathbb{R}^3}\rho|(\dot{u})_\epsilon|^2\,dx-\frac{1}{2}\int_{\mathbb{R}^3}\rho_t|(\dot{u})_\epsilon|^2\,dx,
\end{split}
\ee where we have used the conclusions that
$$\displaystyle\int_{-1}^1\eta^\prime(s)\,ds=0,\quad \int_{-1}^1\eta^\prime(s)\frac{\rho(x,t-\epsilon
s)}{\epsilon}\,ds=(\rho_t)_\epsilon.$$ Here
$$G_\epsilon(t)=\int_{\mathbb{R}^3}(\u)_{\epsilon}(x,t)\int_{-1}^1\eta^\prime(s)\frac{\rho(x,t-\epsilon
s)-\rho(x,t)}{\epsilon}[\dot{u}(x,t-\epsilon
s)-\dot{u}(x,t)]\,ds\,dx.$$ By virtue of the mass equation, it is
easy to check that $\rho_t\in L^\infty([0,T];L^\frac{3}{2})$ for
$T<T^*$. This, combined with the H\"older inequality, the Sobolev
inequality and the fact that
$$\lim\limits_{\epsilon\rightarrow0^+}\|\nabla\dot{u}(\cdot,\cdot-\epsilon
s)-\nabla\dot{u}(\cdot,\cdot)\|_{L^2(\mathbb{R}^3\times[\tau,T-\tau])}=0,$$
gives \be\label{G}
\lim\limits_{\epsilon\rightarrow0^+}\|G_\epsilon\|_{L^1([\tau,T-\tau])}=0.
\ee Putting (\ref{N3}) into (\ref{N2}), and integrating the result
over $(\tau,t)$ for $t\le T-\tau$, we have
 \beq\label{N4}
\begin{split}
&\int_\tau^tG_\epsilon(s)\,ds+\int_\tau^t\int_{\mathbb{R}^3}(\u)_{\epsilon}\dot{u}(\rho_s)_\epsilon\,dx\,ds+
\frac{1}{2}\int_{\mathbb{R}^3}\rho(x,t)|(\dot{u})_\epsilon(x,t)|^2\,dx\\&
-\frac{1}{2}\int_\tau^t\int_{\mathbb{R}^3}\rho_s|(\dot{u})_\epsilon|^2\,dx\,ds
+\int_\tau^t\int_{\mathbb{R}^3}\left[\mu|\nabla(\u)_{\epsilon}|^2
+(\mu+\lambda)|\di(\u)_{\epsilon}|^2\right]\,dx\,ds\\
=&\frac{1}{2}\int_{\mathbb{R}^3}\rho(x,\tau)|(\dot{u})_\epsilon(x,\tau)|^2\,dx+\int_\tau^t\int_{\mathbb{R}^3}
(P_s)_{\epsilon}\di(\u)_{\epsilon}\,dx\,ds  \\
&+\mu \int_\tau^t\int_{\mathbb{R}^3}\nabla[(u\cdot\nabla)
u]_{\epsilon}
\cdot\nabla(\u)_{\epsilon}\,dx\,ds+(\mu+\lambda)\int_\tau^t\int_{\mathbb{R}^3}\di[(u\cdot\nabla)
u]_{\epsilon}\di(\u)_{\epsilon}\,dx\,ds.
\end{split}
\eeq 
Let $\epsilon$ go to zero, and use (\ref{G}) to conclude that
 \beq\label{N6}
\begin{split}
&\frac{1}{2}\int_\tau^t\int_{\mathbb{R}^3}\rho_s|\dot{u}|^2\,dx\,ds+
\frac{1}{2}\int_{\mathbb{R}^3}\rho(x,t)|\dot{u}(x,t)|^2\,dx\\&
+\int_\tau^t\int_{\mathbb{R}^3}\left(\mu|\nabla\u|^2
+(\mu+\lambda)|\di\u|^2\right)\,dx\,ds\\
=&\frac{1}{2}\int_{\mathbb{R}^3}\rho(x,\tau)|\dot{u}(x,\tau)|^2\,dx+\int_\tau^t\int_{\mathbb{R}^3}
P_s\di\u\,dx\,ds +\mu
\int_\tau^t\int_{\mathbb{R}^3}\nabla[(u\cdot\nabla) u]
\cdot\nabla\u\,dx\,ds \\
&+(\mu+\lambda)\int_\tau^t\int_{\mathbb{R}^3}\di[(u\cdot\nabla)
u]\di\u\,dx\,ds.
\end{split}
\eeq For the first term on the left hand side of (\ref{N6}), we make
use of the mass equation and integration by parts. Then we arrive at
\be\label{N7}
\begin{split}
\frac{1}{2}\int_\tau^t\int_{\mathbb{R}^3}\rho_s|\dot{u}|^2\,dx\,ds=\int_\tau^t\int_{\mathbb{R}^3}\rho\dot{u}\cdot
[(u\cdot\nabla)\dot{u}] \,dx\,ds.
\end{split}
\ee (\ref{N-1}) implies that
$$
\rho \dot{u}=-\nabla P+\mu\Delta u+(\mu+\lambda)\nabla\mathrm{div}u.
$$
Putting this into (\ref{N7}), we have \be\label{N8}
\begin{split}
\frac{1}{2}\int_\tau^t\int_{\mathbb{R}^3}\rho_s|\dot{u}|^2\,dx\,ds=&-\int_\tau^t\int_{\mathbb{R}^3}u\otimes\nabla
P:\nabla\u \,dx\,ds+\mu\int_\tau^t\int_{\mathbb{R}^3}u\otimes\Delta
u: \nabla\dot{u}
\,dx\,ds\\&+(\mu+\lambda)\int_\tau^t\int_{\mathbb{R}^3}u\otimes\nabla\mathrm{div}u:
\nabla\dot{u} \,dx\,ds.
\end{split}
\ee The combination of (\ref{N6}) and (\ref{N8}) gives
 \bex
\begin{split}
& \frac{1}{2}\int_{\mathbb{R}^3}\rho(x,t)|\dot{u}(x,t)|^2\,dx
+\int_\tau^t\int_{\mathbb{R}^3}\left(\mu|\nabla\u|^2
+(\mu+\lambda)|\di\u|^2\right)\,dx\,ds\\
=&\frac{1}{2}\int_{\mathbb{R}^3}\rho(x,\tau)|\dot{u}(x,\tau)|^2\,dx+\int_\tau^t\int_{\mathbb{R}^3}
(P_s\di\u + u\otimes\nabla P:\nabla\u)\,dx\,ds
\\
&+\mu \int_\tau^t\int_{\mathbb{R}^3}\Big[\nabla[(u\cdot\nabla)u]
\cdot\nabla\u-u\otimes\Delta u: \nabla\dot{u}\Big]\,dx\,ds
\\&+(\mu+\lambda)\int_\tau^t\int_{\mathbb{R}^3}\Big[\di[(u\cdot\nabla)
u]\di\u-u\otimes\nabla\mathrm{div}u: \nabla\dot{u}\Big]\,dx\,ds.
\end{split}
\eex Let $\tau\rightarrow0^+$ (take subsequence if necessary), we
have\beq\label{dt rho u t-1}
\begin{split}
& \frac{1}{2}\int_{\mathbb{R}^3}\rho(x,t)|\dot{u}(x,t)|^2\,dx
+\int_0^t\int_{\mathbb{R}^3}\left(\mu|\nabla\u|^2
+(\mu+\lambda)|\di\u|^2\right)\,dx\,ds\\
\le&\int_0^t\int_{\mathbb{R}^3} (P_s\di\u + u\otimes\nabla
P:\nabla\u)\,dx\,ds
\\
&+\mu \int_0^t\int_{\mathbb{R}^3}\Big[\nabla[(u\cdot\nabla)u]
\cdot\nabla\u-u\otimes\Delta u: \nabla\dot{u}\Big]\,dx\,ds
\\&+(\mu+\lambda)\int_0^t\int_{\mathbb{R}^3}\Big[\di[(u\cdot\nabla)
u]\di\u-u\otimes\nabla\mathrm{div}u:
\nabla\dot{u}\Big]\,dx\,ds+C\\=& \sum\ls_{i=1}^3(III)_i+C.
\end{split}
\eeq
 For $(III)_1$, using (\ref{full N-S+1})$_1$, integration by
parts, (\ref{blowup-2.1}) and H\"older inequality, we have
\beq\label{blowup-(III) 1}
\begin{split}
(III)_1=&\lim\limits_{\delta\rightarrow0^+}\int_0^t\int_{\mathbb{R}^3}
[P_s\di\u + u\otimes\nabla P:\nabla(\u)_\delta]
\\=&
\lim\limits_{\delta\rightarrow0^+}\int_0^t\int_{\mathbb{R}^3}\Big[
(\rho\theta)_s\di\u-\rho\theta(\nabla u)^\prime:\nabla(\u)_\delta
-\rho\theta u\cdot\nabla\di(\u)_\delta \Big]\\
=&\lim\limits_{\delta\rightarrow0^+}\int_0^t\int_{\mathbb{R}^3}\Big[(\rho\theta)_s\di\u-\rho\theta(\nabla u)^\prime:\nabla(\u)_\delta
+\di(\rho\theta u)\di(\u)_\delta  \Big]\\
=&\int_0^t\int_{\mathbb{R}^3}\Big[\rho\dot{\theta}\di\u-\rho\theta(\nabla
u)^\prime:\nabla\u  \Big]\\
 \le&C\int_0^t\left(\|\sqrt{\rho}\dot{\theta}\|_{L^2}\|\nabla\u\|_{L^2}+\|\sqrt[4]{\rho}\theta\|_{L^4}\|\nabla u\|_{L^4}\|\nabla \u\|_{L^2}\right).
\end{split}
\eeq
 For $(III)_2$, we have
\bex
\begin{split}
(III)_2=&\mu
\lim\limits_{\delta\rightarrow0^+}\int_0^t\int_{\mathbb{R}^3}\Big[\nabla[(u_\delta\cdot\nabla)u_\delta]
\cdot\nabla(\u)_\delta-u_\delta\otimes\Delta u_\delta:
\nabla(\dot{u})_\delta\Big]\\=& \mu
\lim\limits_{\delta\rightarrow0^+}\int_0^t\int_{\mathbb{R}^3}\Big[\di(\Delta
u_\delta\otimes u_\delta)-\Delta[(u_\delta\cdot\nabla)u_\delta]
\Big]\cdot(\u)_\delta.
\end{split}
\eex It is not difficult to check that
$$\di(\de u_\delta\otimes u_\delta)-\de[(u_\delta\cdot\nabla)
u_\delta] =\nabla_k(\di
u_\delta\nabla_ku_\delta)-\nabla_k(\nabla_ku_\delta^j\nabla_ju_\delta)-\nabla_j(\nabla_ku_\delta^j\nabla_ku_\delta),$$
where $\nabla_k=\frac{\partial}{\partial x_k}$. Thus we have
\be\label{blowupIII2}
\begin{split}
(III)_2=& -\mu
\lim\limits_{\delta\rightarrow0^+}\int_0^t\int_{\mathbb{R}^3}(\di
u_\delta\nabla_ku_\delta)\cdot\nabla_k(\u)_\delta+\mu
\lim\limits_{\delta\rightarrow0^+}\int_0^t\int_{\mathbb{R}^3}
(\nabla_ku_\delta^j\nabla_ju_\delta) \cdot\nabla_k(\u)_\delta\\&+\mu
\lim\limits_{\delta\rightarrow0^+}\int_0^t\int_{\mathbb{R}^3}
(\nabla_ku_\delta^j\nabla_ku_\delta) \cdot\nabla_j(\u)_\delta \\=&
-\mu  \int_0^t\int_{\mathbb{R}^3}(\di
u\nabla_ku)\cdot\nabla_k(\u)+\mu
 \int_0^t\int_{\mathbb{R}^3} (\nabla_ku^j\nabla_ju)
\cdot\nabla_k(\u)\\&+\mu  \int_0^t\int_{\mathbb{R}^3}
(\nabla_ku^j\nabla_ku) \cdot\nabla_j(\u)\\ \le&  C\int_0^t\|\nabla
\u\|_{L^2}\|\nabla u\|_{L^4}^2,
\end{split}
\ee where we have used integration by parts and the H\"older
inequality. Note that
$$
\di(\nabla\di u_\delta\otimes
u_\delta)-\nabla\di[(u_\delta\cdot\nabla)
u_\delta]=\nabla(\nabla_ju_\delta^j\nabla_iu_\delta^i)-\nabla(\nabla_ju_\delta^i\nabla_iu_\delta^j)-\nabla_i(\nabla
u_\delta^i\nabla_ju_\delta^j).
$$ Similar to the arguments for $(III)_2$, we have
\be\label{blowupIII3} (III)_3\le C\int_0^t\|\nabla
\u\|_{L^2}\|\nabla u\|_{L^4}^2.\ee
 Substituting (\ref{blowup-(III) 1}), (\ref{blowupIII2}) and (\ref{blowupIII3}) into (\ref{dt rho u t-1}), and using the Cauchy inequality and
(\ref{blowup-2.1}), we have \bex
\begin{split}
&\frac{1}{2}\int_{\mathbb{R}^3}\rho|\u|^2+\int_0^t\int_{\mathbb{R}^3}\left(\mu|\nabla\u|^2+(\mu+\lambda)|\di\u|^2\right)
\\ \le&\frac{\mu}{2}\int_0^t\|\nabla
\u\|_{L^2}^2+C\int_0^t\|\sqrt{\rho}\dot{\theta}\|_{L^2}^2+C\int_0^t\|\sqrt[4]{\rho}\theta\|_{L^4}^4+C\int_0^t\|\nabla
u\|_{L^4}^4+C.
\end{split}
\eex The first term on the right hand side can be absorbed by the
left. Thus we have \beq\label{dt rho u t-2}
\begin{split}
\int_{\mathbb{R}^3}\rho|\u|^2+\int_0^t\int_{\mathbb{R}^3}|\nabla\u|^2\le
C\int_0^t\|\sqrt{\rho}\dot{\theta}\|_{L^2}^2+C\int_0^t\left(\|\nabla
u\|_{L^4}^4+\|\nabla \theta\|_{L^2}^4\right)+C,
\end{split}
\eeq where we have used \bex\begin{split}
\int_0^t\|\sqrt[4]{\rho}\theta\|_{L^4}^4=&\int_0^t\int_{\mathbb{R}^3}\rho\theta^4\le\int_0^t\|\rho\|_{L^3}\|\theta\|_{L^6}^4\\
\le &C\int_0^t\|\nabla\theta\|_{L^2}^4.
\end{split}
\eex
 The next
step is to get some estimates for $\theta$. We rewrite (\ref{full
N-S+1})$_3$ as follows: \be\label{rho theta+ =} \rho
\dot{\theta}+\rho\theta\mathrm{div}u=\frac{\mu}{2}\left|\nabla
u+(\nabla
u)^\prime\right|^2+\lambda(\mathrm{div}u)^2+\kappa\Delta\theta. \ee
Multiplying (\ref{rho theta+ =}) by $\dot{\theta}$ and integrating
by parts over $\mathbb{R}^3$, we have \be\label{dt nabla
theta}\begin{split} &\int_{\mathbb{R}^3}\rho
|\dot{\theta}|^2+\frac{\kappa}{2}\frac{d}{dt}\int_{\mathbb{R}^3}|\nabla
\theta|^2\\
=&-\int_{\mathbb{R}^3}\rho\theta\mathrm{div}u\dot{\theta}+\int_{\mathbb{R}^3}\left(\frac{\mu}{2}\left|\nabla
u+(\nabla
u)^\prime\right|^2+\lambda(\mathrm{div}u)^2\right)\theta_t\\&+\int_{\mathbb{R}^3}\left(\frac{\mu}{2}\left|\nabla
u+(\nabla
u)^\prime\right|^2+\lambda(\mathrm{div}u)^2\right)u\cdot\nabla\theta+\kappa\int_{\mathbb{R}^3}\Delta\theta
u\cdot\nabla\theta\\=& \sum\limits_{i=1}^4(IV)_i.
\end{split}
\ee For $(IV)_1$, using H\"older inequality, (\ref{blowup-2.1}),
(\ref{3d-mass conservation}) and Cauchy inequality, we have
\be\label{blowup-(IV) 1}\begin{split} (IV)_1\leq
C\|\sqrt{\rho}\dot{\theta}\|_{L^2}\|\nabla
u\|_{L^4}\|\sqrt{\rho}\|_{L^{12}}\|\theta\|_{L^6}
\le\frac{1}{8}\|\sqrt{\rho}\dot{\theta}\|_{L^2}^2+C\|\nabla
u\|_{L^4}^4+C\|\nabla\theta\|_{L^2}^4.
\end{split}
\ee For $(IV)_2$, we have \bex\begin{split}
(IV)_2=&\frac{d}{dt}\int_{\mathbb{R}^3}\left(\frac{\mu}{2}\left|\nabla u+(\nabla u)^\prime\right|^2+\lambda(\mathrm{div}u)^2\right)\theta-\mu\int_{\mathbb{R}^3}\left(\nabla u+(\nabla u)^\prime\right):\left(\nabla u_t+(\nabla u_t)^\prime\right)\theta\\
&-2\lambda\int_{\mathbb{R}^3}\mathrm{div}u\mathrm{div}u_t\theta\\
=&\frac{d}{dt}\int_{\mathbb{R}^3}\left(\frac{\mu}{2}\left|\nabla
u+(\nabla u)^\prime\right|^2+\lambda(\mathrm{div}u)^2\right)\theta
-\mu\int_{\mathbb{R}^3}\left(\nabla u+(\nabla u)^\prime\right):\left(\nabla \u+(\nabla \u)^\prime\right)\theta\\
&+\mu\int_{\mathbb{R}^3}\left(\nabla u+(\nabla u)^\prime\right):\left(\nabla u\cdot\nabla u+(\nabla u\cdot\nabla u)^\prime\right)\theta\\
&+\mu\int_{\mathbb{R}^3}\left(\nabla u+(\nabla u)^\prime\right)\cdot
(u\cdot\nabla)\left(\nabla u+(\nabla u)^\prime\right)\theta
-2\lambda\int_{\mathbb{R}^3}\mathrm{div}u\mathrm{div}\dot{u}\theta\\
&+2\lambda\int_{\mathbb{R}^3}\mathrm{div}u(\nabla u)^\prime:\nabla
u\theta+2\lambda\int_{\mathbb{R}^3}u\cdot\nabla
\mathrm{div}u\mathrm{div}u\theta.
\end{split}
\eex Using integration by parts, we have \be\label{blowup-(IV)
2-1}\begin{split}
(IV)_2=&\frac{d}{dt}\int_{\mathbb{R}^3}\left(\frac{\mu}{2}\left|\nabla
u+(\nabla u)^\prime\right|^2+\lambda(\mathrm{div}u)^2\right)\theta
-\mu\int_{\mathbb{R}^3}\left(\nabla u+(\nabla u)^\prime\right):\left(\nabla \u+(\nabla \u)^\prime\right)\theta\\
&+\mu\int_{\mathbb{R}^3}\left(\nabla u+(\nabla
u)^\prime\right):\left(\nabla u\cdot\nabla u+(\nabla u\cdot\nabla
u)^\prime\right)\theta-\mu\int_{\mathbb{R}^3}\frac{|\nabla u+(\nabla
u)^\prime|^2}{2}\mathrm{div}u\theta\\
&-\mu\int_{\mathbb{R}^3}\frac{|\nabla u+(\nabla
u)^\prime|^2}{2}u\cdot\nabla\theta
-2\lambda\int_{\mathbb{R}^3}\mathrm{div}u\mathrm{div}\dot{u}\theta+2\lambda\int_{\mathbb{R}^3}\mathrm{div}u(\nabla
u)^\prime:\nabla
u\theta\\
&-\lambda\int_{\mathbb{R}^3} (\mathrm{div}u)^3\theta-\lambda\int_{\mathbb{R}^3}|\mathrm{div}u|^2u\cdot\nabla\theta\\
=&\sum\limits_{i=1}^9(IV)_{2,i}.
\end{split}
\ee\\
For $(IV)_{2,2}$ and $(IV)_{2,6}$, using the H\"older inequality,
the Sobolev inequality, we have \be\label{blowup-(IV) 2,2 and
2,6}\begin{split} (IV)_{2,2}+(IV)_{2,6}\le C\|\nabla
\u\|_{L^2}\|\nabla u\|_{L^3}\|\theta\|_{L^6}\le C\|\nabla
\u\|_{L^2}\|\nabla u\|_{L^3}\|\nabla \theta\|_{L^2}.
\end{split}
\ee Since $\nabla
u=\nabla\Delta^{-1}\left(\nabla\mathrm{div}u-\nabla\times\mathrm{curl}u\right)$,
we apply Calderon-Zygmund inequality to
get\be\label{Calderon-Zygmund}\begin{split} \|\nabla
u\|_{L^3}\le&C\|\mathrm{curl}u\|_{L^3}+C\|\mathrm{div}u\|_{L^3}.
\end{split}
\ee Taking $\mathrm{curl}$ on both sides of (\ref{full N-S+1})$_2$,
we have \be\label{equation of curlu}
\mu\Delta(\mathrm{curl}u)=\mathrm{curl}(\rho \dot{u}). \ee By
(\ref{equation of curlu}), the $L^2$-estimates of the elliptic
equations and (\ref{blowup-2.1}), we have \beq\label{H 1 of curl u}
\begin{split}\|\nabla \mathrm{curl}u\|_{L^{2}} \le C\|\rho \dot{u}\|_{L^{2}}\le
 C\|\sqrt{\rho} \dot{u}\|_{L^{2}}.
\end{split}
\eeq By (\ref{equation of G}), (\ref{Calderon-Zygmund}) and (\ref{H
1 of curl u}), together with the Sobolev inequality, we have
 \be\label{tu3}\begin{split} \|\nabla
u\|_{L^3}\le&C\|\mathrm{curl}u\|_{L^3}+C\|G\|_{L^3}+C\|\rho\theta\|_{L^3}
\le C\|\mathrm{curl}u\|_{H^1}+C\|G\|_{H^1}+C\|\rho\|_{L^6}\|\theta\|_{L^6}\\
\le&C\|\mathrm{curl}u\|_{L^2}+C\|\mathrm{div}u\|_{L^2}+C\|\nabla(\mathrm{curl}u)\|_{L^2}+C\|\nabla G\|_{L^2}+C\|\nabla\theta\|_{L^2}\\
\le&C\|\nabla
u\|_{L^2}+C\|\sqrt{\rho}\u\|_{L^2}+C\|\nabla\theta\|_{L^2},
\end{split}
\ee where we have used (\ref{blowup-2.1}), (\ref{3d-mass
conservation}) and (\ref{H 1 of u}). Substituting (\ref{tu3}) into
(\ref{blowup-(IV) 2,2 and 2,6}), we obtain \be\label{blowup-(IV) 2,2
and 2,6+1}\begin{split} (IV)_{2,2}+(IV)_{2,6}\le C\|\nabla
\u\|_{L^2}\left(\|\nabla
u\|_{L^2}+\|\sqrt{\rho}\u\|_{L^2}+\|\nabla\theta\|_{L^2}\right)\|\nabla
\theta\|_{L^2}.
\end{split}
\ee  For $(IV)_{2,3}$, $(IV)_{2,4}$, $(IV)_{2,7}$ and $(IV)_{2,8}$,
using the H\"older inequality, the Sobolev inequality and the
Calderon-Zygmund inequality, we have \be\label{blowup-(IV)
3,4,7,8}\begin{split}
&(IV)_{2,3}+(IV)_{2,4}+(IV)_{2,7}+(IV)_{2,8}\\
\le&C\int_{\mathbb{R}^3} |\nabla u|^3|\theta|\le C\|\nabla
u\|_{L^{\frac{18}{5}}}^3\|\theta\|_{L^6}\le C\|\nabla
u\|_{L^{\frac{18}{5}}}^3\|\nabla \theta\|_{L^2}\\
\le&C\|\mathrm{curl} u\|_{L^{\frac{18}{5}}}^3\|\nabla
\theta\|_{L^2}+C\|\mathrm{div} u\|_{L^{\frac{18}{5}}}^3\|\nabla
\theta\|_{L^2}\\
\le&C\|\mathrm{curl}u\|_{L^{\frac{18}{5}}}^3\|\nabla
\theta\|_{L^2}+C\|G\|_{L^{\frac{18}{5}}}^3\|\nabla
\theta\|_{L^2}+C\|\rho\theta\|_{L^{\frac{18}{5}}}^3\|\nabla
\theta\|_{L^2}.
\end{split}
\ee Using the H\"older inequality again, together with
(\ref{blowup-2.1}), (\ref{3d-mass conservation}), the Sobolev
inequality, the Gagliardo-Nirenberg inequality, (\ref{H 1 of u}),
(\ref{equation of G}) and (\ref{H 1 of curl u}), we get
\be\label{blowup-(IV) 3,4,7,8+1}\begin{split}
&(IV)_{2,3}+(IV)_{2,4}+(IV)_{2,7}+(IV)_{2,8}\\
\le&C\|\mathrm{curl}u\|_{L^2}\|\nabla\mathrm{curl}u\|_{L^2}^2\|\nabla
\theta\|_{L^2}+C\|G\|_{L^2}\|\nabla G\|_{L^2}^2\|\nabla
\theta\|_{L^2}+C\|\rho\|_{L^9}^3\|\theta\|_{L^6}^3\|\nabla
\theta\|_{L^2}\\
\le& C\|\sqrt{\rho} \dot{u}\|_{L^{2}}^2\|\nabla
\theta\|_{L^2}+C\|\nabla \theta\|_{L^2}^4.
\end{split}
\ee  For $(IV)_{2,5}$ and $(IV)_{2,9}$, using the H\"older
inequality, the Cauchy inequality, the Sobolev inequality, the
Gagliardo-Nirenberg inequality and (\ref{H 1 of u}), we have
\be\label{blowup-(IV) 5,9}\begin{split} (IV)_{2,5}+(IV)_{2,9}\le
&C\int_{\mathbb{R}^3}|\nabla u|^2|u||\nabla\theta|\le C\|\nabla
u\|_{L^4}^2\|u\|_{L^6}\|\nabla\theta\|_{L^3}\\ \le&
C\|\nabla u\|_{L^4}^4+C\|\nabla u\|_{L^2}^2\|\nabla\theta\|_{L^2}\|\nabla^2\theta\|_{L^2}\\
\le& C\|\nabla
u\|_{L^4}^4+C\|\nabla\theta\|_{L^2}\|\nabla^2\theta\|_{L^2}.
\end{split}
\ee From the standard elliptic estimates and (\ref{rho theta+ =}),
we have \be\label{H 2 of theta}\begin{split}
\|\nabla^2\theta\|_{L^2}\le&C\|\rho \dot{\theta}\|_{L^2}+C\|\rho\theta\mathrm{div}u\|_{L^2}+C\|\nabla u\|_{L^4}^2\\
\le&C\|\sqrt{\rho} \dot{\theta}\|_{L^2}+C\|\rho\|_{L^{12}}\|\theta\|_{L^6}\|\nabla u\|_{L^4}+C\|\nabla u\|_{L^4}^2\\
\le& C\|\sqrt{\rho}\dot{\theta}\|_{L^2}+C\|\nabla
u\|_{L^4}^2+C\|\nabla\theta\|_{L^2}^2,
\end{split}
\ee where we have used the H\"older inequality, (\ref{blowup-2.1}),
(\ref{3d-mass conservation}), the Sobolev inequality and the Cauchy
inequality. Substituting (\ref{H 2 of theta}) into (\ref{blowup-(IV)
5,9}), and using the Cauchy inequality, we have
\be\label{blowup-(IV) 5,9+1}\begin{split} (IV)_{2,5}+(IV)_{2,9}\le
\frac{1}{8}\|\sqrt{\rho}\dot{\theta}\|_{L^2}^2+C\|\nabla
u\|_{L^4}^4+C\|\nabla\theta\|_{L^2}^4+C\|\nabla\theta\|_{L^2}^2.
\end{split}
\ee Substituting (\ref{blowup-(IV) 2,2 and 2,6+1}),
(\ref{blowup-(IV) 3,4,7,8+1}) and (\ref{blowup-(IV) 5,9+1}) into
(\ref{blowup-(IV) 2-1}), and using the Cauchy inequality, we have
\be\label{blowup-(IV) 2}\begin{split} (IV)_2\le&
\frac{d}{dt}\int_{\mathbb{R}^3}\left(\frac{\mu}{2}\left|\nabla
u+(\nabla u)^\prime\right|^2+\lambda(\mathrm{div}u)^2\right)\theta
\\ &+C\|\nabla
\u\|_{L^2}\left(\|\nabla
u\|_{L^2}+\|\sqrt{\rho}\u\|_{L^2}+\|\nabla\theta\|_{L^2}\right)\|\nabla
\theta\|_{L^2}+C\|\sqrt{\rho} \dot{u}\|_{L^{2}}^2\|\nabla
\theta\|_{L^2}\\ &+C\|\nabla
\theta\|_{L^2}^4+\frac{1}{8}\|\sqrt{\rho}\dot{\theta}\|_{L^2}^2+C\|\nabla
u\|_{L^4}^4+C\|\nabla\theta\|_{L^2}^2.
\end{split}
\ee For $(IV)_3$, using (\ref{blowup-(IV) 5,9}) and
(\ref{blowup-(IV) 5,9+1}), we have \be\label{blowup-(IV)
3}\begin{split} (IV)_3\le&C\int_{\mathbb{R}^3}|\nabla
u|^2|u||\nabla\theta|\le\frac{1}{8}\|\sqrt{\rho}\dot{\theta}\|_{L^2}^2+C\|\nabla
u\|_{L^4}^4+C\|\nabla\theta\|_{L^2}^4+C\|\nabla\theta\|_{L^2}^2.
\end{split}
\ee For $(IV)_4$, using the H\"older inequality, the Sobolev
inequality, the Gagliardo-Nirenberg inequality, (\ref{H 1 of u}),
(\ref{H 2 of theta}) and the Young inequality, we have
\be\label{blowup-(IV) 4}\begin{split} (IV)_4\le&
C\|\Delta\theta\|_{L^2}\|u\|_{L^6}\|\nabla\theta\|_{L^3}\le
C\|\Delta\theta\|_{L^2}\|\nabla
u\|_{L^2}\|\nabla\theta\|_{L^2}^\frac{1}{2}\|\nabla^2\theta\|_{L^2}^\frac{1}{2}\\
\le&C
\|\nabla\theta\|_{L^2}^\frac{1}{2}\|\nabla^2\theta\|_{L^2}^\frac{3}{2}\le
\frac{1}{8}\|\sqrt{\rho}\dot{\theta}\|_{L^2}^2+C\|\nabla\theta\|_{L^2}^4
+C\|\nabla u\|_{L^4}^4+C\|\nabla\theta\|_{L^2}^2.
\end{split}
\ee Putting (\ref{blowup-(IV) 1}), (\ref{blowup-(IV) 2}),
(\ref{blowup-(IV) 3}) and (\ref{blowup-(IV) 4}) into (\ref{dt nabla
theta}), integrating the resulting inequality over $[0,t]$ for
$t\in(0,T^*)$, and using the Cauchy inequality, (\ref{blowup-2.1}),
(\ref{u 2nabla u 2}) and (\ref{H 1 of u}), we have \be\label{dt
nabla theta+1a}\begin{split} \int_0^t\int_{\mathbb{R}^3}\rho
|\dot{\theta}|^2+\int_{\mathbb{R}^3}|\nabla \theta|^2
\le&C\int_0^t\|\nabla u\|_{L^4}^4+C\int_{\mathbb{R}^3}|\nabla
u|^2|\theta|+\varepsilon\int_0^t\|\nabla
\u\|_{L^2}^2\\&+C_\varepsilon\int_0^t \left(\|\nabla
u\|_{L^2}^2+\|\sqrt{\rho}\u\|_{L^2}^2+\|\nabla
\theta\|_{L^2}^2\right)\|\nabla \theta\|_{L^2}^2+C.
\end{split}
\ee For the second term of the right hand side of (\ref{dt nabla
theta+1a}), we have \be\label{D}\begin{split}
C\int_{\mathbb{R}^3}|\nabla u|^2|\theta|\le&C\|\nabla
u\|_{L^\frac{12}{5}}^2\|\theta\|_{L^6}\le C\|\mathrm{curl}
u\|_{L^\frac{12}{5}}^2\|\nabla\theta\|_{L^2}+
C\|\mathrm{div} u\|_{L^\frac{12}{5}}^2\|\nabla\theta\|_{L^2}\\
\le&C\left(\|\mathrm{curl}u\|_{L^2}^{\frac{3}{2}}\|\nabla\mathrm{curl}u\|_{L^2}^{\frac{1}{2}}+
\|G\|_{L^2}^{\frac{3}{2}}\|\nabla
G\|_{L^2}^{\frac{1}{2}}\right)\|\nabla\theta\|_{L^2}+C\|\rho\theta\|_{L^\frac{12}{5}}^2\|\nabla\theta\|_{L^2}\\
\le&
C\|\sqrt{\rho}\u\|_{L^2}^{\frac{1}{2}}\|\nabla\theta\|_{L^2}+C\|\rho\theta\|_{L^2}^{\frac{3}{2}}\|\rho\theta\|_{L^6}^{\frac{1}{2}}\|\nabla\theta\|_{L^2}
\\
\le&C\|\sqrt{\rho}\u\|_{L^2}^{\frac{1}{2}}\|\nabla\theta\|_{L^2}+C\|\nabla\theta\|_{L^2}^{\frac{3}{2}}
\le\frac{1}{2}\|\nabla\theta\|_{L^2}^2+C\|\sqrt{\rho}\u\|_{L^2}+C,
\end{split}
\ee where we have used the H\"older inequality, the Calderon-Zygmund
inequality, the Gagliardo-Nirenberg inequality, (\ref{blowup-2.1}),
(\ref{H 1 of u}), (\ref{equation of G}), (\ref{H 1 of curl u}), the
Sobolev inequality and the Young inequality. Substituting (\ref{D})
into (\ref{dt nabla theta+1a}), we have \be\label{dt nabla
theta+1}\begin{split} \int_0^t\int_{\mathbb{R}^3}\rho
|\dot{\theta}|^2+\int_{\mathbb{R}^3}|\nabla \theta|^2
\le&C\int_0^t\|\nabla u\|_{L^4}^4+C
\|\sqrt{\rho}\u\|_{L^2}+\varepsilon\int_0^t\|\nabla \u\|_{L^2}^2
\\&+C_\varepsilon\int_0^t
\left(\|\nabla u\|_{L^2}^2+\|\sqrt{\rho}\u\|_{L^2}^2+\|\nabla
\theta\|_{L^2}^2\right)\|\nabla \theta\|_{L^2}^2+C.
\end{split}
\ee Multiplying (\ref{dt nabla theta+1}) by $2C$ and adding the
resulting inequality into (\ref{dt rho u t-2}), we have \bex
\begin{split}
&C\int_0^t\int_{\mathbb{R}^3}\rho |\dot{\theta}|^2+2C\int_{\mathbb{R}^3}|\nabla \theta|^2+\int_{\mathbb{R}^3}\rho|\u|^2+\int_0^t\int_{\mathbb{R}^3}|\nabla\u|^2\\
\leq&2C^2\int_0^t\|\nabla u\|_{L^4}^4+2C^2\|\sqrt{\rho}\u\|_{L^2}+2\varepsilon C\int_0^t\|\nabla \u\|_{L^2}^2\\
&+2CC_\varepsilon\int_0^t \left(\|\nabla
u\|_{L^2}^2+\|\sqrt{\rho}\u\|_{L^2}^2+\|\nabla
\theta\|_{L^2}^2\right)\|\nabla \theta\|_{L^2}^2+2C^2.
\end{split}
\eex Taking $\varepsilon$ sufficiently small, together with the
Cauchy inequality, we have \be\label{non-sum}
\begin{split}
&\int_{\mathbb{R}^3}(|\nabla
\theta|^2+\rho|\u|^2)+\int_0^t\int_{\mathbb{R}^3}(\rho
|\dot{\theta}|^2+|\nabla\u|^2) \\ \le&C\int_0^t\|\nabla u\|_{L^4}^4+
C\int_0^t \left(\|\nabla
u\|_{L^2}^2+\|\sqrt{\rho}\u\|_{L^2}^2+\|\nabla
\theta\|_{L^2}^2\right)\|\nabla \theta\|_{L^2}^2+C.
\end{split}
\ee For the first term of the right hand side of (\ref{non-sum}),
similar to (\ref{tu3}), we have \be\label{non-sum+1}
\begin{split}
\int_0^t\|\nabla u\|_{L^4}^4
\le&C\int_0^t\|\mathrm{curl}u\|_{L^4}^4+C\int_0^t\|G\|_{L^4}^4+C\int_0^t\|\nabla\theta\|_{L^2}^4\\
\le&C\int_0^t\|\mathrm{curl}u\|_{L^2}\|\nabla\mathrm{curl}u\|_{L^2}^3+C\int_0^t\|G\|_{L^2}\|\nabla
G\|_{L^2}^3+C\int_0^t\|\nabla\theta\|_{L^2}^4\\
\le&C\int_0^t\|\sqrt{\rho}\u\|_{L^2}^3+C\int_0^t\|\nabla\theta\|_{L^2}^4.
\end{split}
\ee
 By (\ref{non-sum}), (\ref{non-sum+1}) and the Cauchy inequality, we have \be\label{non-sum+2}
\begin{split}
&\int_{\mathbb{R}^3}(|\nabla
\theta|^2+\rho|\u|^2)+\int_0^t\int_{\mathbb{R}^3}(\rho
|\dot{\theta}|^2+|\nabla\u|^2)\\
\le&C\int_0^t\|\sqrt{\rho}\u\|_{L^2}^3+C \int_0^t \left(\|\nabla
u\|_{L^2}^2+\|\sqrt{\rho}\u\|_{L^2}^2+\|\nabla
\theta\|_{L^2}^2\right)\|\nabla \theta\|_{L^2}^2+C\\ \le&C\int_0^t
\left(\|\nabla u\|_{L^2}^2+\|\sqrt{\rho}\u\|_{L^2}^2+\|\nabla
\theta\|_{L^2}^2\right)\left(\|\nabla
\theta\|_{L^2}^2+\|\sqrt{\rho}\u\|_{L^2}^2\right)+C.
\end{split}
\ee By (\ref{blowup-2.1}), (\ref{u 2nabla u 2}), (\ref{H 1 of u}),
we have
$$
\int_0^t \left(\|\nabla
u\|_{L^2}^2+\|\sqrt{\rho}\u\|_{L^2}^2+\|\nabla
\theta\|_{L^2}^2\right)\le C,
$$ for any $t\in(0,T^*)$. This, together with (\ref{non-sum+2}) and
the Gronwall inequality, deduces (\ref{H 2 of u}).
\endpf

\begin{corollary}\label{blowup-cor 3.5}
Under the conditions of Theorem \ref{blowup-th:1.1} and
(\ref{blowup-2.1}), it holds that for any $t\in(0,T^*)$ \be\label{u
infty} \|u\|_{L^\infty}+\|\nabla
u\|_{L^r}+\int_0^t\int_{\mathbb{R}^3}\left(\rho|\theta_t|^2+|\nabla^2\theta|^2\right)\,dx\,ds\le
C, \ee for any $r\in[2,6]$, and any $t\in[0,\infty)$.
\end{corollary}
\pf Similar to (\ref{tu3}), using (\ref{blowup-2.1}), the Sobolev
inequality and (\ref{H 2 of u}), we have \be\label{nabla L 6}
\begin{split}
\|\nabla u\|_{L^6}\le&
C\|\mathrm{culr}u\|_{L^6}+C\|G\|_{L^6}+C\|\rho\theta\|_{L^6}\\
\le&C\|\nabla\mathrm{culr}u\|_{L^2}+C\|\nabla
G\|_{L^2}+C\|\nabla\theta\|_{L^2}\\
\le&C\|\sqrt{\rho}\u\|_{L^2}+C\|\nabla\theta\|_{L^2}\le C.
\end{split} \ee
By (\ref{H 1 of u}), we have \bex \|\nabla u\|_{L^2}\le C. \eex
This, together with (\ref{nabla L 6}) and the interpolation
inequality and the Sobolev inequality, implies \be\label{3d-u infty
and nabla L r}\begin{split} \|u\|_{L^\infty}+\|\nabla u\|_{L^r}\le
C,
\end{split}
\ee for any $r\in[2,6]$.

By (\ref{H 2 of theta}), we have
 \be\label{3d-L 2 H 2 of theta}\begin{split}
\|\nabla^2\theta\|_{L^2}^2\le&
C\|\sqrt{\rho}\dot{\theta}\|_{L^2}^2+C\|\nabla
u\|_{L^4}^4+C\|\nabla\theta\|_{L^2}^4\\
 \le&
C\|\sqrt{\rho}\dot{\theta}\|_{L^2}^2+C\|\nabla
u\|_{L^2}\|\nabla u\|_{L^6}^3+C\|\nabla\theta\|_{L^2}^2\\
\le&C\|\sqrt{\rho}\dot{\theta}\|_{L^2}^2+C\|\sqrt{\rho}\u\|_{L^2}^2+C\|\nabla\theta\|_{L^2}^2,
\end{split}
\ee where we have used the Gagliardo-Nirenberg inequality, (\ref{H 2
of u}), (\ref{nabla L 6}) and (\ref{3d-u infty and nabla L r}).

By (\ref{3d-L 2 H 2 of theta}), (\ref{H 2 of u}) and (\ref{H 1 of
u}), we get
$$
\int_0^t\int_{\mathbb{R}^3}|\nabla^2\theta|^2\le C.
$$
Recall $\dot{\theta}=\theta_t+u\cdot\nabla\theta$, we have
\bex\begin{split} \int_0^t\int_{\mathbb{R}^3}\rho|\theta_t|^2\le&
C\int_0^t\int_{\mathbb{R}^3}\rho|\dot{\theta}|^2+C\int_0^t\int_{\mathbb{R}^3}\rho|u\cdot\nabla\theta|^2\\
\le&C\int_0^t\int_{\mathbb{R}^3}\rho|\dot{\theta}|^2+C\int_0^t\int_{\mathbb{R}^3}|\nabla\theta|^2\le
C,
\end{split}\eex where we have used (\ref{blowup-2.1}), (\ref{H 1 of u}), (\ref{H 2 of u}), and (\ref{3d-u infty and nabla L r}).
\endpf

\begin{lemma}\label{blowup-le: H 2 of theta}Under the conditions of Theorem \ref{blowup-th:1.1} and (\ref{blowup-2.1}), it holds that for any $t\in(0,T^*)$
\be\label{non-H 2 of theta}
\int_{\mathbb{R}^3}\rho|\theta_t|^2\,dx+\int_0^t\int_{\mathbb{R}^3}|\nabla\theta_s|^2\,dx\,ds\le
C. \ee
\end{lemma}
\pf To get the estimate (\ref{non-H 2 of theta}), we use the
arguments similar to the proof of Lemma \ref{blowup-le: int rho u
t}, i.e., mollifying each term in (\ref{full N-S+1})$_3$,
differentiating the result w.r.t. t, multiplying by
$(\theta_{\epsilon,\delta})_t$, integrating by parts, and passing to
the limit. Then we arrive at \be\label{dt H 2 of theta}\begin{split}
&\frac{1}{2}\int_{\mathbb{R}^3}\rho|\theta_t|^2+\kappa\int_0^t\int_{\mathbb{R}^3}|\nabla\theta_t|^2\\
\le&-
\int_0^t\int_{\mathbb{R}^3}\rho_t\left(\frac{\theta_t}{2}+u\cdot\nabla\theta+\theta\mathrm{div}u\right)\theta_t-\int_0^t\int_{\mathbb{R}^3}\rho(
u_t\cdot\nabla\theta+u\cdot\nabla\theta_t+\theta_t\mathrm{div}u)\theta_t\\&-\int_0^t
\int_{\mathbb{R}^3}\rho\theta\mathrm{div}u_t\theta_t+\mu\int_0^t\int_{\mathbb{R}^3}\left(\nabla
u+(\nabla u)^\prime\right):\left(\nabla u_t+(\nabla
u_t)^\prime\right)\theta_t\\&+2\lambda\int_0^t\int_{\mathbb{R}^3}\mathrm{div}u\mathrm{div}u_t\theta_t+C\\
=&\sum\limits_{i=1}^5(V)_i+C,
\end{split}
\ee where \be\label{blowup-(V) 1}\begin{split} (V)_1=&
-\int_0^t\int_{\mathbb{R}^3}\rho
u\cdot\nabla\theta_t\left(\frac{\theta_t}{2}+u\cdot\nabla\theta+\theta\mathrm{div}u\right)-\int_0^t\int_{\mathbb{R}^3}\rho
u\cdot\frac{\nabla\theta_t}{2}\theta_t\\&-\int_0^t\int_{\mathbb{R}^3}\rho
u\cdot\left(\nabla
(u\cdot\nabla)\theta+u\cdot\nabla\nabla\theta\right)\theta_t-\int_0^t\int_{\mathbb{R}^3}\rho
u\cdot\left(\nabla\theta\mathrm{div}u+\theta\nabla\mathrm{div}u\right)\theta_t\\=&\sum\limits_{i=1}^4(V)_{1,i}.
\end{split}
\ee For $(V)_{1,1}$, we have \be\label{blowup-(V) 1,1}\begin{split}
(V)_{1,1} \le& \frac{\kappa}{24}\int_0^t\|\nabla
\theta_t\|_{L^2}^2+C\int_0^t\|u\|_{L^\infty}^2\|\sqrt{\rho}\theta_t\|_{L^2}^2
+C\int_0^t\|u\|_{L^\infty}^4\|\nabla\theta\|_{L^2}^2\\&+C\int_0^t\|u\|_{L^\infty}^2\|\theta\|_{L^6}^2\|\nabla u\|_{L^3}^2\\
\le&\frac{\kappa}{24}\int_0^t\int_{\mathbb{R}^3}|\nabla\theta_t|^2+C\int_0^t\int_{\mathbb{R}^3}\rho|\theta_t|^2
+C\int_0^t\int_{\mathbb{R}^3}|\nabla\theta|^2,
\end{split}
\ee where we have used the Cauchy inequality, (\ref{blowup-2.1}) and
(\ref{u infty}).\\ For $(V)_{1,2}$ and $(V)_{1,3}$, using the Cauchy
inequality, (\ref{blowup-2.1}) and (\ref{u infty}) again, we have
\be\label{blowup-(V) 1,2}\begin{split}
(V)_{1,2}\le&\frac{\kappa}{24}\int_0^t\int_{\mathbb{R}^3}|\nabla\theta_t|^2+C\int_0^t\int_{\mathbb{R}^3}\rho
|\theta_t|^2,
\end{split}
\ee and \be\label{blowup-(V) 1,3}\begin{split} (V)_{1,3}\le&C
\int_0^t\|\sqrt{\rho}\theta_t\|_{L^2}\left(\|\nabla\theta\|_{L^6}\|\nabla
u\|_{L^3}+\|\nabla^2\theta\|_{L^2}\right)\\
\le&C\int_0^t\|\sqrt{\rho}\theta_t\|_{L^2}^2+C\int_0^t\|\nabla^2\theta\|_{L^2}^2.
\end{split}
\ee For $(V)_{1,4}$, integrating by parts, we have

\be\label{blowup-(V) 1,4a}\begin{split} (V)_{1,4}
=&-\int_0^t\int_{\mathbb{R}^3}\rho
u\cdot\nabla\theta\mathrm{div}u\theta_t-\displaystyle\frac{1}{2\mu+\lambda}\int_0^t\int_{\mathbb{R}^3}\rho\theta
u\cdot\nabla G\theta_t
\\&-\displaystyle\frac{1}{2\mu+\lambda}\int_0^t\int_{\mathbb{R}^3}\rho\theta u\cdot\nabla (\rho\theta)\theta_t\\
=&-\int_0^t\int_{\mathbb{R}^3}\rho
u\cdot\nabla\theta\mathrm{div}u\theta_t-\displaystyle\frac{1}{2\mu+\lambda}\int_0^t\int_{\mathbb{R}^3}\rho\theta
u\cdot\nabla G\theta_t
\\
&+\displaystyle\frac{1}{2(2\mu+\lambda)}\int_0^t\int_{\mathbb{R}^3}\rho^2\theta^2
u\cdot\nabla\theta_t+\displaystyle\frac{1}{2(2\mu+\lambda)}\int_0^t\int_{\mathbb{R}^3}\rho^2\theta^2\mathrm{div}u\theta_t
.
\end{split}
\ee Furthermore, we get \be\label{blowup-(V) 1,4}\begin{split}
(V)_{1,4} \le&C\int_0^t\|\sqrt{\rho}\theta_t\|_{L^2}\|\nabla
u\|_{L^3}\|\nabla \theta\|_{L^6}+C\int_0^t\|\nabla
G\|_{L^2}\|\theta_t\|_{L^6}\|\theta\|_{L^6}\|\rho\|_{L^{6}}\\&
+C\int_0^t\|\theta\|_{L^6}^2\|\rho\|_{L^{12}}^2\|\nabla
\theta_t\|_{L^2}+C\int_0^t\|\nabla
u\|_{L^3}\|\theta_t\|_{L^6}\|\theta\|_{L^6}^2\|\rho\|_{L^{12}}^2\\
\le&\frac{\kappa}{24}\int_0^t\int_{\mathbb{R}^3}|\nabla\theta_t|^2+C\int_0^t\|\sqrt{\rho}\theta_t\|_{L^2}^2+C\int_0^t\|\nabla^2
\theta\|_{L^2}^2\\&+C\int_0^t\|\sqrt{\rho}\u\|_{L^2}^2
+C\int_0^t\|\nabla\theta\|_{L^2}^2,
\end{split}
\ee where we have used the H\"older inequality, the Sobolev
inequality, (\ref{blowup-2.1}), (\ref{3d-mass conservation}),
(\ref{H 1 of G}), (\ref{H 2 of u}), (\ref{u infty}) and the Cauchy
inequality. Substituting (\ref{blowup-(V) 1,1}), (\ref{blowup-(V)
1,2}), (\ref{blowup-(V) 1,3}) and (\ref{blowup-(V) 1,4}) into
(\ref{blowup-(V) 1}), we have \be\label{blowup-(V) 1+1}\begin{split}
(V)_1\le&\frac{\kappa}{8}\int_0^t\int_{\mathbb{R}^3}|\nabla\theta_t|^2+C\int_0^t\int_{\mathbb{R}^3}\rho|\theta_t|^2
+C\int_0^t\int_{\mathbb{R}^3}|\nabla\theta|^2\\&+C\int_0^t\|\nabla^2
\theta\|_{L^2}^2+C\int_0^t\|\sqrt{\rho}\u\|_{L^2}^2.
\end{split}
\ee For $(V)_2$, \be\label{blowup-(V) 2}\begin{split} (V)_2
=&-\int_0^t\int_{\mathbb{R}^3}\rho
\u\cdot\nabla\theta\theta_t+\int_0^t\int_{\mathbb{R}^3}\rho
(u\cdot\nabla)u\cdot\nabla\theta\theta_t
-\int_0^t\int_{\mathbb{R}^3}\rho u\cdot\nabla\theta_t\theta_t\\&-\int_0^t\int_{\mathbb{R}^3}\rho|\theta_t|^2\mathrm{div}u\\
\le&C\int_0^t\|\sqrt{\rho}\u\|_{L^2}\|\theta_t\|_{L^6}\|\nabla
\theta\|_{L^3} +C\int_0^t\|\sqrt{\rho}\theta_t\|_{L^2}\|\nabla
u\|_{L^6}\|\nabla
\theta\|_{L^3}\\&+C\int_0^t\|\sqrt{\rho}\theta_t\|_{L^2}\|\nabla
\theta_t\|_{L^2}
+C\int_0^t\|\sqrt{\rho}\theta_t\|_{L^2}\|\theta_t\|_{L^6}\|\nabla u\|_{L^3}\\
\leq&\frac{\kappa}{8}\int_0^t\int_{\mathbb{R}^3}|\nabla\theta_t|^2
+C\int_0^t\|\nabla
\theta\|_{L^3}^2+C\int_0^t\|\sqrt{\rho}\theta_t\|_{L^2}^2\\
\leq&\frac{\kappa}{8}\int_0^t\int_{\mathbb{R}^3}|\nabla\theta_t|^2
+C\int_0^t\|\nabla \theta\|_{L^2}^2+C\int_0^t\|\nabla^2
\theta\|_{L^2}^2+C\int_0^t\|\sqrt{\rho}\theta_t\|_{L^2}^2,
\end{split}
\ee
where we have used the H\"older inequality, the Sobolev inequality, (\ref{blowup-2.1}), (\ref{H 2 of u}) and (\ref{u infty}).\\
For $(V)_3$, integrating by parts, we have
 \be\label{blowup-(V)
3}\begin{split}
(V)_3=&-\int_0^t\int_{\mathbb{R}^3}\rho\theta\mathrm{div}\u\theta_t+\int_0^t\int_{\mathbb{R}^3}\rho\theta\mathrm{div}(u\cdot\nabla u)\theta_t\\
=&-\int_0^t\int_{\mathbb{R}^3}\rho\theta\mathrm{div}\u\theta_t+\int_0^t\int_{\mathbb{R}^3}\rho\theta\nabla
u:(\nabla u)'\theta_t
+\int_0^t\int_{\mathbb{R}^3}\rho\theta u\cdot\nabla \mathrm{div}u\theta_t\\
=&-\int_0^t\int_{\mathbb{R}^3}\rho\theta\mathrm{div}\u\theta_t+\int_0^t\int_{\mathbb{R}^3}\rho\theta\nabla
u:(\nabla u)'\theta_t
+\displaystyle\frac{1}{2\mu+\lambda}\int_0^t\int_{\mathbb{R}^3}\rho\theta
u\cdot\nabla G\theta_t
\\
&-\displaystyle\frac{1}{2\mu+\lambda}\int_0^t\int_{\mathbb{R}^3}\frac{\rho^2}{2}\theta^2\mathrm{div}u\theta_t-\displaystyle\frac{1}{2\mu+\lambda}\int_0^t
\int_{\mathbb{R}^3}\frac{\rho^2}{2}\theta^2u\cdot\nabla\theta_t.
\end{split}
\ee Furthermore, using the H\"older inequality, the Sobolev
inequality, (\ref{blowup-2.1}), (\ref{3d-mass conservation}),
(\ref{H 1 of G}), (\ref{H 2 of u}), (\ref{u infty}) and the Young
inequality, we have \be\label{blowup-(V) 3+1}\begin{split} (V)_3
\le&C \int_0^t\|\nabla
\u\|_{L^2}\|\theta_t\|_{L^6}\|\theta\|_{L^6}\|\rho\|_{L^{6}}+C\int_0^t\|\rho\|_{L^\infty}
\|\theta\|_{L^6}\|\nabla u\|_{L^2}\|\nabla
u\|_{L^6}\|\theta_t\|_{L^6}\\&+ C\int_0^t\|u\|_{L^\infty}\|\nabla
G\|_{L^2}\|\theta_t\|_{L^6}\|\theta\|_{L^6}\|\rho\|_{L^{6}}
+C\int_0^t\|\nabla
u\|_{L^3}\|\theta_t\|_{L^6}\|\theta\|_{L^6}^2\|\rho\|_{L^{12}}^2\\&+C\int_0^t\|u\|_{L^\infty}\|\theta\|_{L^6}^2\|\rho\|_{L^{12}}^2\|\nabla
\theta_t\|_{L^2}\\
\le&
\frac{\kappa}{8}\int_0^t\int_{\mathbb{R}^3}|\nabla\theta_t|^2+C\int_0^t\|\nabla
\u\|_{L^2}^2+C\int_0^t\|\nabla\theta\|_{L^2}^2+
C\int_0^t\|\sqrt{\rho}\u\|_{L^2}^2.
\end{split}
\ee Similar to $(V)_2$ and $(V)_3$, for $(V)_4$ and $(V)_5$, we
deduce \be\label{blowup-(V) 4 and 5}
\begin{split}
(V)_4+(V)_5\le&C\int_0^t\|\nabla\u\|_{L^2}\|\nabla
u\|_{L^3}\|\theta_t\|_{L^6}+C\int_0^t\int_{\mathbb{R}^3}|\nabla
u|^3|\theta_t|\\&+C\int_0^t\int_{\mathbb{R}^3}|u|^2|\nabla
u|^4+\frac{\kappa}{16}\int_0^t\int_{\mathbb{R}^3}|\nabla\theta_t|^2\\
\le&\frac{3\kappa}{32}\int_0^t\int_{\mathbb{R}^3}|\nabla\theta_t|^2+C\int_0^t\|\nabla\u\|_{L^2}^2+\bar{C}\int_0^t\|\nabla
u\|_{L^2}\|\nabla
u\|_{L^6}^2\|\theta_t\|_{L^6}\\&+C\int_0^t\|u\|_{L^6}^2\|\nabla
u\|_{L^6}^4\\
\le&\frac{\kappa}{8}\int_0^t\int_{\mathbb{R}^3}|\nabla\theta_t|^2+C\int_0^t\|\nabla\u\|_{L^2}^2+C\int_0^t\|\nabla
u\|_{L^2}^2,
\end{split}
\ee
where we have used the H\"older inequality, integration by parts, the Cauchy inequality, (\ref{u infty}), the interpolation inequality and the Sobolev inequality.\\
Putting (\ref{blowup-(V) 1+1}), (\ref{blowup-(V) 2}),
(\ref{blowup-(V) 3+1}) and (\ref{blowup-(V) 4 and 5}) into (\ref{dt
H 2 of theta}), we have \be\label{dt H 2 of theta+1}\begin{split}
&\frac{1}{2}\int_{\mathbb{R}^3}\rho|\theta_t|^2+\kappa\int_0^t\int_{\mathbb{R}^3}|\nabla\theta_t|^2\\
\le& C\int_0^t\int_{\mathbb{R}^3}\rho|\theta_t|^2
+C\int_0^t\int_{\mathbb{R}^3}|\nabla\theta|^2+C\int_0^t\|\nabla^2
\theta\|_{L^2}^2+C\int_0^t\|\sqrt{\rho}\u\|_{L^2}^2\\&+C\int_0^t\|\nabla
\u\|_{L^2}^2+C\int_0^t\|\nabla u\|_{L^2}^2+C.
\end{split}
\ee This combined with (\ref{H 1 of u}), (\ref{H 2 of u}) and
(\ref{u infty}) completes the proof of Lemma \ref{blowup-le: H 2 of
theta}.
\endpf
\begin{corollary}Under the conditions of Theorem
\ref{blowup-th:1.1} and (\ref{blowup-2.1}), it holds that for any
$t\in( 0,T^*)$ \be\label{cor:H 2 of theta}
\int_{\mathbb{R}^3}|\nabla^2\theta|^2\,dx\le C. \ee
\end{corollary}
\pf It follows from (\ref{full N-S+1})$_3$, (\ref{blowup-2.1}),
(\ref{3d-mass conservation}), (\ref{H 2 of u}), (\ref{u infty}),
(\ref{non-H 2 of theta}) and the interpolation inequality that
\bex\begin{split} \|\nabla^2\theta\|_{L^2}\le& C\|\rho
\theta_t\|_{L^2}+ C\|\rho
u\cdot\nabla\theta\|_{L^2}+C\|\rho\theta\mathrm{div}u\|_{L^2}+C\|\nabla
u\|_{L^4}^2\\ \le&C\|\sqrt{\rho} \theta_t\|_{L^2}+
C\|\nabla\theta\|_{L^2}+C\|\rho\|_{L^6}\|\theta\|_{L^6}\|\mathrm{div}u\|_{L^6}+C\\
\le&C.
\end{split}\eex
\endpf\\

By (\ref{H 2 of u}), (\ref{cor:H 2 of theta}) and the Sobolev
inequality, we get the following corollary which is the desired one,
i.e., (\ref{non-uniform_est1}).
\begin{corollary}Under the conditions of Theorem \ref{blowup-th:1.1} and (\ref{blowup-2.1}), it holds that for any $t\in(0,T^*)$
\be\label{cor:H 2 of u} \|\theta\|_{L^\infty(0,t;L^\infty)}\le C.
\ee
\end{corollary}

\section{Global strong solution}\label{3d-sec 4}

 \setcounter{equation}{0} \setcounter{theorem}{0}

 In this section,
we shall prove the global existence and uniqueness of the strong
solution. Since the local existence and uniqueness of the strong
solution has been obtained in \cite{cho-Kim: perfect gas} under the
conditions of Theorem \ref{3d-th:1.1}, we assume that $T^*>0$ is the
maximal existence time of the strong solution. We shall
prove $T^*=\infty$ by using contradiction arguments.

 Remark \ref{3d-remark blowup} says that if $T^*<\infty$, then
{\small\be\label{3d-blow up criterion+1}
\lim\sup\limits_{T\nearrow T^*}\left(\|\rho\|_{L^\infty(0,T;
L^\infty)}+\|\rho\theta\|_{L^4(0,T; L^\frac{12}{5})}+\|\rho^\frac{1}{4}u\|_{L^\infty(0,T;L^4)}+\big\| |u||\nabla
u| \big\|_{L^2(0,T;L^2)}\right)=\infty
\ee} for all $\mu$ and $\lambda$ satisfying only the physical restriction
(\ref{viscosity assumption}).

If $T^*<\infty$, our aim is to prove that (\ref{3d-blow up criterion+1})
is not true under the conditions of Theorem \ref{3d-th:1.1}, which is the desired contradiction.

To do this, we define \be\label{AT} A(T)=\sup\limits_{0\le t\le
T}\int_{\mathbb{R}^3}|\nabla
u|^2+\int_0^T\int_{\mathbb{R}^3}\frac{\rho|\u|^2}{\mu}, \ee and
\be\label{BT} B(T)=\sup\limits_{0\le t\le T}\int_{\mathbb{R}^3}\rho
\theta^2+\int_0^T\int_{\mathbb{R}^3}\kappa|\nabla\theta|^2. \ee The
following proposition plays a crucial role in the section.
\begin{proposition}\label{prop 3.1}
Assume that the initial data satisfies (\ref{initial data}),
(\ref{3d-initial assumption}), and (\ref{3d-compatibility}). If the
strong solution $(\rho, u, \theta)$ satisfies \be\label{a priori
assumption} A(T)\le 2\tilde{E} K_1,\ B(T)\le 2K_2,\ 0\le \rho\le
2\bar{\rho},\ (x,t)\in \mathbb{R}^3\times[0,T], \ee then
\be\label{3d-result from the a p} A(T)\le \frac{3}{2}\tilde{E} K_1,\
B(T)\le\frac{3}{2}K_2,\ 0\le \rho\le \frac{3}{2}\bar{\rho},\
(x,t)\in \mathbb{R}^3\times[0,T], \ee provided that
$m_0\le\varepsilon_0$. Here $\displaystyle
m_0=\int_{\mathbb{R}^3}\rho_0(x)\,dx$,
$\tilde{E}=\frac{(14\mu+9\lambda)}{2\mu}+\frac{6\bar{\rho}K_2}{\mu(\mu+\lambda)K_1}
+\frac{8\bar{\rho}\kappa K_2}{\mu(\mu+\lambda)^2K_1}+1$, and
{\small\bex
\begin{split}\varepsilon_0=&\min\left\{C_3,\ \frac{\check{\mathcal
{C}}(2\mu+\lambda)^{6}}{\tilde{E}^3},\
\check{\mathcal
{C}}(2\mu+\lambda)^\frac{12}{7},\
\frac{\check{\mathcal
{C}}\mu^{12}(2\mu+\lambda)^{12}}{\tilde{E}^{12}},\ \check{\mathcal
{C}}(2\mu+\lambda)^{36}\kappa^{12}\right\},\end{split}\eex }
 where
{\small\bex\begin{split}C_3=\min\left\{\frac{\check{\mathcal
{C}}\kappa^6(\mu+\lambda)^6\mu^6}{\Big(\kappa(\mu+\lambda)+1\Big)^6},\
\check{\mathcal {C}}\mu^3\kappa^3(2\mu+\lambda)^6,\
\frac{\check{\mathcal {C}}\mu^6}{\tilde{E}^2},\ \check{\mathcal
{C}}\mu^2(2\mu+\lambda)^8, \frac{\check{\mathcal
{C}}\kappa^4\mu^2}{(2\mu+\lambda)^8\tilde{E}^6},\
\frac{\check{\mathcal {C}}\kappa^4}{\tilde{E}^4\mu^2},\
\frac{\check{\mathcal
{C}}\kappa^6}{\tilde{E}^3}\right\},\end{split}\eex} for some constant
$\check{\mathcal {C}}>0$ depending on $\bar{\rho}, K_1, K_2,$ and
some other known constants but independent of $\mu,\lambda,\kappa,$
and $t$.
\end{proposition}

With Proposition \ref{prop 3.1}, we shall get $T^*=\infty$. More
precisely, we obtain the following corollary.
\begin{corollary}\label{3d-cor 3.2}
With Proposition \ref{prop 3.1}, it holds that $T^*=\infty$ with
(\ref{a priori assumption}) valid for any $0\le T<\infty$.
\end{corollary}
\pf
 If $T_1>0$ is the
maximal time such that (\ref{a priori assumption}) is valid, then
$T_1=T^*$. For otherwise, (\ref{3d-result from the a p}) implies
that $T_1$ is not the maximal time.

With $T_1=T^*$, (\ref{a priori assumption}) and the $L^1$-bound of
$\rho$ (see Lemma \ref{3d-le:3.1}), one can easily get
\be\label{contradiction} \|\rho\|_{L^\infty(0,t;
L^\infty)}+\|\rho\theta\|_{L^4(0,t;
L^\frac{12}{5})}+\|\rho^\frac{1}{4}u\|_{L^\infty(0,t;L^4)}\le
\breve{C}(1+t^\frac{1}{8}) \ee for all $t\in[0,T^*)$, where
$\breve{C}$ is a positive constant independent of t. Using the
H\"older inequality, the inequality $\|u\|_{L^6}\le
\breve{C}\|\nabla u\|_{L^2}$ and (\ref{a priori assumption}), we
have \bex\begin{split} \big\| |u||\nabla u|
\big\|_{L^2}\le\|u\|_{L^6}\|\nabla u\|_{L^3} \le \breve{C}\|\nabla
u\|_{L^3}.
\end{split}
\eex  This together with (\ref{a priori assumption}), the
$L^1$-bound of $\rho$ and the estimate for $\|\nabla u\|_{L^3}$ (see
(\ref{3d-nabla u 3}) for the detail) gives
 \be\label{contradiction+1} \big\|
|u||\nabla u| \big\|_{L^2(0,t;L^2)}\le \breve{C}(1+t^\frac{1}{4}).
\ee

Therefore, if $T^*<\infty$, (\ref{contradiction}) and (\ref{contradiction+1}) will contradict with (\ref{3d-blow up criterion+1}). Thus, $T^*$ must be $\infty$. Then
we get $T_1=T^*=\infty$ which implies that the strong solution exists globally in time and that (\ref{a priori assumption}) is
valid for any $T\in[0,\infty)$.
\begin{remark}
Corollary \ref{3d-cor 3.2} means that if Proposition \ref{prop 3.1}
is valid, the global existence of strong solutions will be got. The
uniqueness of the solutions can be referred to \cite{cho-Kim:
perfect gas}. The proof of Theorem \ref{3d-th:1.1} is
complete.
\end{remark}
Let's come back to prove Proposition \ref{prop 3.1}. Throughout the
rest of the paper, we denote generic constants by $C$ depending on
$\bar{\rho}, K_1, K_2,$ and some other known constants but independent of $\mu,\lambda,\kappa,$ and $t$.\\

{\noindent\bf \underline{Proof of Proposition \ref{prop 3.1}}}:

\begin{lemma}\label{3d-le:3.2}
Under the conditions of Proposition \ref{prop 3.1}, it holds that
\be\label{3d-energy1} \int_0^T\int_{\mathbb{R}^3}|\nabla
u|^2\,dx\leq m_0^\frac{1}{2},\ee provided $$ m_0\le
\frac{\kappa^6(\mu+\lambda)^6\mu^6}{C^6\Big(\kappa(\mu+\lambda)+1\Big)^6}\triangleq
C_1.
$$
\end{lemma}
\pf Multiplying (\ref{full N-S+1})$_2$ by $u$, integrating by parts
over $\mathbb{R}^3$, and using the Cauchy inequality, we have
\bex\begin{split}
&\frac{1}{2}\frac{d}{dt}\int_{\mathbb{R}^3}\rho|u|^2+\int_{\mathbb{R}^3}\left(\mu|\nabla
u|^2+(\mu+\lambda)|\mathrm{div}u|^2\right)\\=&\int_{\mathbb{R}^3}\rho\theta\mathrm{div}u\le
(\mu+\lambda)\int_{\mathbb{R}^3}|\mathrm{div}u|^2+\frac{1}{4(\mu+\lambda)}\int_{\mathbb{R}^3}\rho^2\theta^2.
\end{split}
\eex This, together with the H\"older inequality, the Sobolev
inequality, (\ref{a priori assumption}) and (\ref{3d-mass
conservation}), deduces \be\label{3d-dt rho u}\begin{split}
\frac{1}{2}\frac{d}{dt}\int_{\mathbb{R}^3}\rho|u|^2+\mu\int_{\mathbb{R}^3}|\nabla
u|^2
\le\frac{1}{4(\mu+\lambda)}\|\rho\|_{L^3}^2\|\theta\|_{L^6}^2\le
\frac{Cm_0^\frac{2}{3}}{\mu+\lambda}\|\nabla\theta\|_{L^2}^2.
\end{split}
\ee Integrating (\ref{3d-dt rho u}) over $[0,T]$, and using (\ref{a
priori assumption}) again, we have \bex\begin{split}
\int_0^T\int_{\mathbb{R}^3}|\nabla u|^2 \le
\frac{C}{\mu}\|\rho_0\|_{L^\frac{3}{2}}\|\nabla u_0\|_{L^2}^2+
\frac{Cm_0^\frac{2}{3}}{\mu(\mu+\lambda)}\int_0^T\|\nabla\theta\|_{L^2}^2
\le
\left(1+\frac{1}{\kappa(\mu+\lambda)}\right)\frac{Cm_0^\frac{2}{3}}{\mu}\le
m_0^\frac{1}{2},
\end{split}
\eex provided that
$$
m_0\le
\frac{\kappa^6(\mu+\lambda)^6\mu^6}{C^6\Big(\kappa(\mu+\lambda)+1\Big)^6}\triangleq
C_1.
$$

\endpf
\begin{lemma}\label{3d-le:3.3}Under the conditions of Proposition \ref{prop 3.1}, it holds that
\be\label{3d-A}A(T)\le \frac{3 \tilde{E}K_1}{2},\ee provided that $$
m_0\le \min\left\{C_1,\
\frac{\mu^3\kappa^3(2\mu+\lambda)^6}{216C^3},\
\frac{\mu^6}{36C^2\tilde{E}^2},\
\frac{\mu^2(2\mu+\lambda)^8}{36C^2}\right\}\triangleq C_2,
$$ where $\tilde{E}=\frac{(14\mu+9\lambda)}{2\mu}+\frac{6\bar{\rho}K_2}{\mu(\mu+\lambda)K_1}
+\frac{8\bar{\rho}\kappa K_2}{\mu(\mu+\lambda)^2K_1}+1$.
\end{lemma}
\pf Recall from (\ref{dt nabla u 2-1}) \beq\label{3d-dt nabla u
2-1}\begin{split} &\int_{\mathbb{R}^3}\rho
|\u|^2+\frac{1}{2}\frac{d}{dt}\int_{\mathbb{R}^3}\left(\mu|\nabla
u|^2+(\mu+\lambda)|\mathrm{div}u|^2\right)\\=&\frac{d}{dt}\int_{\mathbb{R}^3}P\mathrm{div}u-\frac{1}{2(2\mu+\lambda)}
\frac{d}{dt}\int_{\mathbb{R}^3}P^2
-\frac{1}{2\mu+\lambda}\int_{\mathbb{R}^3}P_tG+\int_{\mathbb{R}^3}\rho
(u\cdot\nabla) u \cdot \u\\=&\sum\limits_{i=1}^4II_i,
\end{split}
\eeq where $G=(2\mu+\lambda)\mathrm{div}u-P$.

 Substituting
(\ref{equation of P t}) into $II_3$, and using integration by parts,
the H\"older inequality and the Sobolev inequality, we have
 \beq\label{3d-II 3}\begin{split}
II_3 \le&
\frac{1}{2\mu+\lambda}\|\rho\theta\|_{L^3}\|u\|_{L^6}\|\nabla
G\|_{L^2}+\frac{1}{2\mu+\lambda}\|\rho\theta\|_{L^3}\|\mathrm{div}u\|_{L^2}
\|G\|_{L^6}\\&+C\|G\|_{L^6}\|\nabla u\|_{L^2}\|\nabla
u\|_{L^3}+\frac{\kappa}{2\mu+\lambda}\|\nabla
G\|_{L^2}\|\nabla\theta\|_{L^2}\\
\le& \frac{C}{2\mu+\lambda}\|\rho\|_{L^6}\|\theta\|_{L^6}\|\nabla
u\|_{L^2}\|\nabla G\|_{L^2}+C\|\nabla G\|_{L^2}\|\nabla
u\|_{L^2}\|\nabla u\|_{L^3}\\&+\frac{\kappa}{2\mu+\lambda}\|\nabla
G\|_{L^2}\|\nabla\theta\|_{L^2}.
\end{split}
\eeq
 By (\ref{equation of G}), (\ref{equation of curlu}), the standard
$L^2$-estimates, and (\ref{a priori assumption}), we get
\beq\label{3d-H 1 of G}
\begin{split}\|\nabla G\|_{L^{2}}\le \|\rho \u\|_{L^{2}}\le
 \sqrt{2\bar{\rho}}\|\sqrt{\rho} \u\|_{L^{2}},
\end{split}
\eeq and\beq\label{3d-H 1 of curl u}
\begin{split}\|\nabla \mathrm{curl}u\|_{L^{2}} \le\frac{1}{\mu} \|\rho \dot{u}\|_{L^{2}}\le\frac{ \sqrt{2\bar{\rho}}}{\mu}
 \|\sqrt{\rho}\u\|_{L^{2}}.
\end{split}
\eeq
 Similar to (\ref{tu3}), using (\ref{3d-H 1 of G}), (\ref{3d-H 1 of curl
u}), the Sobolev inequality, the H\"older inequality, and the
Gagliardo-Nirenberg inequality, we have
 \be\label{3d-nabla u 3}\begin{split} \|\nabla
u\|_{L^3}\le&C\|\mathrm{curl}u\|_{L^3}+\frac{C}{2\mu+\lambda}\|G\|_{L^3}+\frac{C}{2\mu+\lambda}\|\rho\theta\|_{L^3}\\
\le&
C\|\mathrm{curl}u\|_{L^2}^\frac{1}{2}\|\nabla\mathrm{curl}u\|_{L^2}^\frac{1}{2}+\frac{C}{2\mu+\lambda}\|G\|_{L^2}^\frac{1}{2}\|\nabla
G\|_{L^2}^\frac{1}{2}
+\frac{C}{2\mu+\lambda}\|\rho\|_{L^6}\|\theta\|_{L^6}\\
\le&\frac{C}{\sqrt{\mu}}\|\mathrm{curl}u\|_{L^2}^\frac{1}{2}\|\sqrt{\rho}\u\|_{L^2}^\frac{1}{2}+\frac{C}{2\mu+\lambda}\|G\|_{L^2}^\frac{1}{2}\|\sqrt{\rho}\u\|_{L^2}^\frac{1}{2}
+\frac{C}{2\mu+\lambda}\|\rho\|_{L^6}\|\nabla\theta\|_{L^2}.
\end{split}
\ee
 Substituting (\ref{3d-H 1 of G}) into (\ref{3d-II 3}), we have
\beq\label{3d-II 3+1}\begin{split} II_3 \le&
\frac{C}{2\mu+\lambda}\|\rho\|_{L^6}\|\nabla\theta\|_{L^2}\|\nabla
u\|_{L^2}\|\sqrt{\rho}\u\|_{L^2}+C\|\sqrt{\rho}\u\|_{L^2}\|\nabla
u\|_{L^2}\|\nabla u\|_{L^3}\\&+\frac{\kappa
 \sqrt{2\bar{\rho}}}{2\mu+\lambda}\|\sqrt{\rho}\u\|_{L^2}\|\nabla\theta\|_{L^2}.
\end{split}
\eeq For $II_4$, using the H\"older inequality, (\ref{a priori
assumption}), and the Sobolev inequality, we have \be\label{3d-II
4}II_4\le C\|\sqrt{\rho} \u\|_{L^2}\|u\|_{L^6}\|\nabla u\|_{L^3}\le
C\|\sqrt{\rho} \u\|_{L^2}\|\nabla u\|_{L^2}\|\nabla u\|_{L^3}.\ee
Putting (\ref{3d-II 3+1}) and (\ref{3d-II 4}) together, and using
(\ref{3d-nabla u 3}), and the Young inequality, we have

\beq\label{3d-II 3+II 4}\begin{split} II_3+II_4 \le&
\frac{C}{2\mu+\lambda}\|\rho\|_{L^6}\|\nabla\theta\|_{L^2}\|\nabla
u\|_{L^2}\|\sqrt{\rho}\u\|_{L^2}+\frac{C}{\sqrt{\mu}}\|\sqrt{\rho}\u\|_{L^2}^\frac{3}{2}\|\nabla
u\|_{L^2}^\frac{3}{2}\\&+\frac{C}{2\mu+\lambda}\|\sqrt{\rho}\u\|_{L^2}^\frac{3}{2}\|\nabla
u\|_{L^2}\|\rho\theta\|_{L^2}^\frac{1}{2} +\frac{\kappa
 \sqrt{2\bar{\rho}}}{2\mu+\lambda}\|\sqrt{\rho}\u\|_{L^2}\|\nabla\theta\|_{L^2}\\ \le&
\frac{1}{2}\|\sqrt{\rho}\u\|_{L^2}^2+\frac{C}{(2\mu+\lambda)^2}\|\rho\|_{L^6}^2\|\nabla\theta\|_{L^2}^2\|\nabla
u\|_{L^2}^2+\frac{C}{\mu^2}\|\nabla
u\|_{L^2}^6\\&+\frac{C}{(2\mu+\lambda)^4}\|\nabla
u\|_{L^2}^4\|\rho\theta\|_{L^2}^2 +\frac{2\bar{\rho}\kappa^2
}{(2\mu+\lambda)^2}\|\nabla\theta\|_{L^2}^2.
\end{split}
\eeq Substituting (\ref{3d-II 3+II 4}) into (\ref{3d-dt nabla u
2-1}),
 we have
\beq\label{3d-dt nabla u 2}\begin{split}
&\frac{1}{2}\int_{\mathbb{R}^3}\rho
|\u|^2+\frac{1}{2}\frac{d}{dt}\int_{\mathbb{R}^3}\left(\mu|\nabla
u|^2+(\mu+\lambda)|\mathrm{div}u|^2\right)\\
\le&\frac{d}{dt}\int_{\mathbb{R}^3}P\mathrm{div}u-\frac{1}{2(2\mu+\lambda)}
\frac{d}{dt}\int_{\mathbb{R}^3}P^2+
\frac{C}{(2\mu+\lambda)^2}\|\rho\|_{L^6}^2\|\nabla\theta\|_{L^2}^2\|\nabla
u\|_{L^2}^2\\&+\frac{C}{\mu^2}\|\nabla
u\|_{L^2}^6+\frac{C}{(2\mu+\lambda)^4}\|\nabla
u\|_{L^2}^4\|\rho\theta\|_{L^2}^2 +\frac{2\bar{\rho}\kappa^2
}{(2\mu+\lambda)^2}\|\nabla\theta\|_{L^2}^2.
\end{split}
\eeq Integrating (\ref{3d-dt nabla u 2}) over $[0,t]$, and using the
Cauchy inequality, we have \bex\begin{split}
&\int_0^t\int_{\mathbb{R}^3}\rho
|\u|^2+\int_{\mathbb{R}^3}\left(\mu|\nabla
u|^2+(\mu+\lambda)|\mathrm{div}u|^2\right)\\
\le&\int_{\mathbb{R}^3}\left(\mu|\nabla
u_0|^2+(\mu+\lambda)|\mathrm{div}u_0|^2\right)-2\int_{\mathbb{R}^3}\rho_0\theta_0\mathrm{div}u_0+\frac{1}{2\mu+\lambda}
\int_{\mathbb{R}^3}\rho_0^2\theta_0^2\\&+(\mu+\lambda)\int_{\mathbb{R}^3}|\mathrm{div}u|^2+\frac{1}{\mu+\lambda}\int_{\mathbb{R}^3}\rho^2\theta^2+
\frac{C}{(2\mu+\lambda)^2}\int_0^t\|\rho\|_{L^6}^2\|\nabla\theta\|_{L^2}^2\|\nabla
u\|_{L^2}^2\\&+\frac{C}{\mu^2}\int_0^t\|\nabla
u\|_{L^2}^6+\frac{C}{(2\mu+\lambda)^4}\int_0^t\|\nabla
u\|_{L^2}^4\|\rho\theta\|_{L^2}^2 +\frac{4\bar{\rho}\kappa^2
}{(2\mu+\lambda)^2}\int_0^t\|\nabla\theta\|_{L^2}^2,
\end{split}
\eex which, together with (\ref{a priori assumption}), (\ref{3d-mass
conservation}), and (\ref{3d-energy1}), gives
\beq\label{3d-dt nabla
u 2+1}\begin{split} &\int_0^t\int_{\mathbb{R}^3}\rho
|\u|^2+\mu\int_{\mathbb{R}^3}|\nabla
u|^2\\
\le&\mu
K_1+3(\mu+\lambda)K_1+\frac{3(2\mu+\lambda)K_1}{2}+\frac{2\bar{\rho}K_2}{2\mu+\lambda}+\frac{4\bar{\rho}K_2}{\mu+\lambda}+
\frac{Cm_0^\frac{1}{3}A(T)}{(2\mu+\lambda)^2}\frac{K_2}{\kappa}\\&+\frac{CA(T)^2}{\mu^2}\int_0^t\|\nabla
u\|_{L^2}^2+\frac{CK_2A(T)}{(2\mu+\lambda)^4}\int_0^t\|\nabla
u\|_{L^2}^2 +\frac{8\bar{\rho}\kappa K_2}{(2\mu+\lambda)^2}\\
\le&\frac{(14\mu+9\lambda)}{2}K_1+\frac{6\bar{\rho}K_2}{\mu+\lambda}
+\frac{8\bar{\rho}\kappa K_2}{(\mu+\lambda)^2}+
\frac{Cm_0^\frac{1}{3}\tilde{E}K_1}{\kappa(2\mu+\lambda)^2}+\frac{C\tilde{E}^2K_1}{\mu^2}m_0^\frac{1}{2}
+\frac{C\tilde{E}K_1}{(2\mu+\lambda)^4} m_0^\frac{1}{2}
\\ \le&\mu \tilde{E}K_1+\frac{\mu \tilde{E}K_1}{2}=\frac{3\mu \tilde{E}K_1}{2},
\end{split}
\eeq provided that
$$
m_0\le \min\left\{C_1,\
\frac{\mu^3\kappa^3(2\mu+\lambda)^6}{216C^3},\
\frac{\mu^6}{36C^2\tilde{E}^2},\
\frac{\mu^2(2\mu+\lambda)^8}{36C^2}\right\}\triangleq C_2,
$$ where $\tilde{E}=\frac{(14\mu+9\lambda)}{2\mu}+\frac{6\bar{\rho}K_2}{\mu(\mu+\lambda)K_1}
+\frac{8\bar{\rho}\kappa K_2}{\mu(\mu+\lambda)^2K_1}+1$.

 By
(\ref{3d-dt nabla u 2+1}), we get (\ref{3d-A}).
\endpf

\begin{lemma}\label{3d-le:3.3}Under the conditions of Proposition \ref{prop 3.1}, it holds that
\be\label{3d-B(T)}B(T)\le \frac{3K_2}{2},\ee provided that $$ m_0\le
\min\left\{C_2,
\frac{\kappa^4\mu^2}{6^4C^4(2\mu+\lambda)^8\tilde{E}^6},\
\frac{\kappa^4}{6^4C^4\tilde{E}^4\mu^2},\
\frac{\kappa^6}{6^6C^6\tilde{E}^3}\right\}\triangleq C_3.
$$
\end{lemma}
\pf Recall from (\ref{zyy}) \beq\label{3d-dt rho theta}\begin{split}
&\kappa\int_{\mathbb{R}^3}|\nabla\theta|^2+\frac{1}{2}\frac{d}{dt}\int_{\mathbb{R}^3}\rho|\theta|^2\\
=&-\int_{\mathbb{R}^3}\rho\theta^2\mathrm{div}u+\int_{\mathbb{R}^3}\frac{\mu}{2}|\nabla u+(\nabla u)'|^2\theta+\int_{\mathbb{R}^3}\lambda(\mathrm{div}u)^2\theta\\
=&\sum\limits_{i=1}^3III_i.
\end{split}
\eeq For $III_1$, using the H\"older inequality, and the Sobolev
inequality, we have \be\label{3d-III 1}\begin{split}
III_1\le\|\mathrm{div}u\|_{L^2}\|\theta\|_{L^6}^2\|\rho\|_{L^6}\le
C\|\mathrm{div}u\|_{L^2}\|\nabla\theta\|_{L^2}^2\|\rho\|_{L^6}.
\end{split}
\ee For $III_2$ and $III_3$, using the H\"older inequality, and the
Sobolev inequality again, together with (\ref{3d-nabla u 3}), we
have \beq\label{3d-III_2, 3}\begin{split}
&III_2+III_3\\ \le&C(2\mu+\lambda)\|\nabla u\|_{L^2}\|\nabla u\|_{L^3}\|\theta\|_{L^6}\le C(2\mu+\lambda)\|\nabla u\|_{L^2}\|\nabla u\|_{L^3}\|\nabla\theta\|_{L^2}\\
\le& \frac{C(2\mu+\lambda)}{\sqrt{\mu}}\|\nabla
u\|_{L^2}^\frac{3}{2}\|\nabla\theta\|_{L^2}\|\sqrt{\rho}\u\|_{L^2}^\frac{1}{2}+C\|\nabla
u\|_{L^2}\|\nabla\theta\|_{L^2}\|\rho\theta\|_{L^2}^\frac{1}{2}\|\sqrt{\rho}\u\|_{L^2}^\frac{1}{2}
\\&+C\|\nabla
u\|_{L^2}\|\nabla\theta\|_{L^2}^2\|\rho\|_{L^6}.
\end{split}
\eeq
 Substituting (\ref{3d-III 1}) and (\ref{3d-III_2, 3}) into (\ref{3d-dt rho
theta}), and using the Cauchy inequality, we have \beq\label{3d-dt
rho theta+1}\begin{split}
&\kappa\int_{\mathbb{R}^3}|\nabla\theta|^2+\frac{1}{2}\frac{d}{dt}\int_{\mathbb{R}^3}\rho|\theta|^2\\
\le& \frac{C(2\mu+\lambda)}{\sqrt{\mu}}\|\nabla
u\|_{L^2}^\frac{3}{2}\|\nabla\theta\|_{L^2}\|\sqrt{\rho}\u\|_{L^2}^\frac{1}{2}+C\|\nabla
u\|_{L^2}\|\nabla\theta\|_{L^2}\|\rho\theta\|_{L^2}^\frac{1}{2}\|\sqrt{\rho}\u\|_{L^2}^\frac{1}{2}
\\&+C\|\nabla
u\|_{L^2}\|\nabla\theta\|_{L^2}^2\|\rho\|_{L^6}\\
\le&
\frac{\kappa}{2}\int_{\mathbb{R}^3}|\nabla\theta|^2+\frac{C(2\mu+\lambda)^2}{\kappa\mu}\|\nabla
u\|_{L^2}^3\|\sqrt{\rho}\u\|_{L^2}+\frac{C}{\kappa}\|\nabla
u\|_{L^2}^2\|\rho\theta\|_{L^2}\|\sqrt{\rho}\u\|_{L^2}
\\&+C\|\nabla
u\|_{L^2}\|\nabla\theta\|_{L^2}^2\|\rho\|_{L^6}.
\end{split}
\eeq Integrating (\ref{3d-dt rho theta+1}) over $[0,t]$, and using
(\ref{a priori assumption}), the H\"older inequality, (\ref{3d-mass
conservation}), and (\ref{3d-energy1}), we have \bex\begin{split}
B(t) \le&
\int_{\mathbb{R}^3}\rho_0|\theta_0|^2+\frac{C(2\mu+\lambda)^2}{\kappa\mu}\int_0^t\|\nabla
u\|_{L^2}^3\|\sqrt{\rho}\u\|_{L^2}\\&+\frac{C}{\kappa}\int_0^t\|\nabla
u\|_{L^2}^2\|\rho\theta\|_{L^2}\|\sqrt{\rho}\u\|_{L^2}
+C\int_0^t\|\nabla u\|_{L^2}\|\nabla\theta\|_{L^2}^2\|\rho\|_{L^6}\\
\le&
K_2+\frac{C(2\mu+\lambda)^2A(T)^\frac{3}{2}}{\kappa\sqrt{\mu}}\|\nabla
u\|_{L^2([0,t];L^2)}+\frac{C\sqrt{K_2}A(T)\sqrt{\mu}}{\kappa}\|\nabla
u\|_{L^2([0,t];L^2)}
\\&+C\sqrt{A(T)}m_0^\frac{1}{6}\frac{K_2}{\kappa}\\ \le&K_2+
\frac{C(2\mu+\lambda)^2\tilde{E}^\frac{3}{2}}{\kappa\sqrt{\mu}}(K_2m_0^\frac{1}{4})+\frac{C\tilde{E}\sqrt{\mu}}{\kappa}(K_2m_0^\frac{1}{4})
+C\sqrt{\tilde{E}}m_0^\frac{1}{6}\frac{K_2}{\kappa}.\end{split} \eex
Thus, \bex\begin{split} B(t) \le K_2+\frac{K_2}{2}=\frac{3K_2}{2},
\end{split}
\eex provided that
$$
m_0\le \min\left\{C_2,
\frac{\kappa^4\mu^2}{6^4C^4(2\mu+\lambda)^8\tilde{E}^6},\
\frac{\kappa^4}{6^4C^4\tilde{E}^4\mu^2},\
\frac{\kappa^6}{6^6C^6\tilde{E}^3}\right\}\triangleq C_3.
$$
\endpf

\begin{lemma}\label{3d-le:rho}Under the conditions of Proposition \ref{prop 3.1}, it holds that
\be\label{3d-upper bound of rho}0\le\rho\le
\frac{3\bar{\rho}}{2},\ee for any $(x,t)\in\mathbb{R}^3\times[0,T]$,
provided that {\small\bex \begin{split}m_0\le&\min\left\{C_3,\ \frac{(2\mu+\lambda)^{6}(\log\frac{3}2)^{6}}{(4C)^6\tilde{E}^3},\
\frac{(2\mu+\lambda)^\frac{12}{7}(\log\frac{3}{2})^\frac{12}{7}}{(4
C)^\frac{12}{7}},\
\frac{\mu^{12}(2\mu+\lambda)^{12}(\log\frac{3}{2})^{12}}{(4C)^{12}\tilde{E}^{12}},\ \frac{(2\mu+\lambda)^{36}\kappa^{12}(\log\frac{3}{2})^{12}}{(4C)^{12}}\right\}\\
\triangleq& \varepsilon_0.\end{split} \eex }
\end{lemma}
\pf The first inequality of (\ref{3d-upper bound of rho}) is
obvious. In fact, this has been obtained in \cite{cho-Kim: perfect
gas} for any
$(x,t)\in\mathbb{R}^3\times[0,T]\subset\mathbb{R}^3\times[0,T^*)$.
We only need to prove the second inequality of (\ref{3d-upper bound
of rho}).

Let us mention that the Zlotnik inequality (see Appendix A) used in
\cite{Huang-Li-Xin} seems not working here. The main ingredient for
handling such the difficulty is an equation obtained from (\ref{full
N-S+1})$_1$ involving $\log\rho$. It was introduced by P.L. Lions
(\cite{Lions2}) to prove global existence of weak solutions of the
compressible isentropic Navier-Stokes equations, and was later used
by B. Desjardins (\cite{Desjardin}) et al to study the regularity of weak
solutions of the compressible isentropic Navier-Stokes equations for
small time under periodic boundary conditions.

More precisely,  for any given $(x, t)\in\mathbb{R}^3\times[0,T]$.
Denote
$$
\rho^\delta(y,s)=\rho(y,s)+\delta\exp\{-\int_0^s\mathrm{div}u\left(X(\tau;x,t),\tau\right)\,d\tau\}>0,
$$ where $X(s;x,t)$ is given by
    $$ \left\{
      \begin{array}{l}
       \displaystyle\frac{d}{ds}X(s;x,t)=u\left(X(s;x,t), s\right),\ 0\le s<t,\\
        X(t; x,t)=x.\\
      \end{array}
      \right.
    $$
It is easy to verify that
$$
\frac{d}{d s}\rho^\delta\left(X(s;x,t),
s\right)+\rho^\delta\left(X(s;x,t),
s\right)\mathrm{div}u\left(X(s;x,t),s\right),
$$
due to (\ref{full N-S+1})$_1$. This gives \be\label{rho
identity}\begin{split} Y^\prime(s)=g(s)+b^\prime(s),
\end{split}\ee where
$$
Y(s)=\log\rho^\delta\left(X(s;x,t),s\right),\
g(s)=-\frac{P\left(X(s;x,t),s\right)}{2\mu+\lambda},\
b(s)=-\frac{1}{2\mu+\lambda}\int_0^sG\left(X(\tau;x,t),\tau\right)\,d\tau,
$$ and $G=(2\mu+\lambda)\mathrm{div}u-P$=$(2\mu+\lambda)\mathrm{div}u-\rho\theta$.

By (\ref{equation of G}) and (\ref{full N-S+1})$_1$, we have
\bex\begin{split} G\left(X(t;x,\tau),\tau\right)=&-(-\Delta)^{-1}
\mathrm{div}\big[(\rho u)_\tau+ \mathrm{div}(\rho u\otimes
u)\big]=-[(-\Delta)^{-1} \mathrm{div}(\rho u)]_\tau  \\&- (-\Delta)^{-1}
\mathrm{div}\mathrm{div}(\rho u\otimes u)\\=&-[(-\Delta)^{-1}
\mathrm{div}(\rho u)]_\tau - u\cdot\nabla(-\Delta)^{-1} \mathrm{div}(\rho
u) + u\cdot\nabla(-\Delta)^{-1} \mathrm{div}(\rho u) \\&- (-\Delta)^{-1}
\mathrm{div}\mathrm{div}(\rho u\otimes
u)\\=&-\frac{d}{d\tau}[(-\Delta)^{-1} \mathrm{div}(\rho
u)] + u\cdot\nabla(-\Delta)^{-1} \mathrm{div}(\rho u) - (-\Delta)^{-1}
\mathrm{div}\mathrm{div}(\rho u\otimes
u)\\=&-\frac{d}{d\tau}[(-\Delta)^{-1} \mathrm{div}(\rho
u)]+[u_i,R_{ij}](\rho u_j),
\end{split}\eex
where $[u_i,R_{ij}]=u_iR_{ij}-R_{ij}u_i$ and
$R_{ij}=\partial_i(-\Delta)^{-1}\partial_j$. This deduces
 \bex\begin{split}
b(t)-b(0)=&\frac{1}{2\mu+\lambda}\int_0^t\left[\frac{d}{d\tau}[(-\Delta)^{-1}
\mathrm{div}(\rho u)]-[u_i,R_{ij}](\rho
u_j)\right]\,d\tau\\=&\frac{1}{2\mu+\lambda}(-\Delta)^{-1}
\mathrm{div}(\rho u)-\frac{1}{2\mu+\lambda}(-\Delta)^{-1}
\mathrm{div}(\rho_0
u_0)-\frac{1}{2\mu+\lambda}\int_0^t[u_i,R_{ij}](\rho u_j)\,d\tau\\
\le& \frac{1}{2\mu+\lambda}\|(-\Delta)^{-1} \mathrm{div}(\rho
u)\|_{L^\infty}+\frac{1}{2\mu+\lambda}\|(-\Delta)^{-1}
\mathrm{div}(\rho_0
u_0)\|_{L^\infty}\\&+\frac{1}{2\mu+\lambda}\int_0^t\|[u_i,R_{ij}](\rho
u_j)\|_{L^\infty}\,d\tau=\sum\limits_{i=1}^3IV_i.
\end{split}\eex For $IV_1$, using the Gagliardo-Nirenberg inequality,
the Sobolev inequality, the Calderon-Zygmund inequality, the
H\"older inequality, (\ref{a priori assumption}), and (\ref{3d-mass
conservation}), we have \be\label{3d-IV 1}\begin{split}
IV_1\le&\frac{C}{2\mu+\lambda}\|(-\Delta)^{-1} \mathrm{div}(\rho
u)\|_{L^6}^\frac{1}{3}\|\nabla(-\Delta)^{-1}
\mathrm{div}(\rho u)\|_{L^4}^\frac{2}{3}\\
\le&\frac{C}{2\mu+\lambda}\|\rho u\|_{L^2}^\frac{1}{3}\|\rho
u\|_{L^4}^\frac{2}{3}
\le\frac{C}{2\mu+\lambda}\|\rho\|_{L^3}^\frac{1}{3} \|
u\|_{L^6}^\frac{1}{3}\|\rho\|_{L^{12}}^\frac{2}{3}\|
u\|_{L^6}^\frac{2}{3}\\
\le&\frac{Cm_0^\frac{1}{6}}{2\mu+\lambda}\|\nabla
u\|_{L^2}\le\frac{Cm_0^\frac{1}{6}\sqrt{\tilde{E}}}{2\mu+\lambda}.
\end{split}\ee
Similarly, for $IV_2$, we have \be\label{3d-IV 2}\begin{split}
IV_2\le\frac{Cm_0^\frac{1}{6}\sqrt{\tilde{E}}}{2\mu+\lambda}.
\end{split}\ee
 Since $u(\cdot,t)\in
W^{1,6}(\mathbb{R}^3)$, $\rho u(\cdot,t)\in
L^{12}(\mathbb{R}^3)$ and $\frac{1}{4}=\frac{1}{6}+\frac{1}{12}$, in view of the conclusions by Desjardins ((33),\cite{Desjardin}) or by Choe-Jin (Section 4, \cite{choe-jin}) and references therein, it holds that
 \bex
\|[u_i,R_{ij}](\rho u_j)\|_{W^{1,4}}\le C\|u\|_{W^{1,6}}\|\rho
u\|_{L^{12}}. \eex This, combined with (\ref{a priori assumption}),
(\ref{3d-mass conservation}), the Sobolev inequality, the
Calderon-Zygmund inequality similar to (\ref{Calderon-Zygmund}),
gives \be\label{RR1}\begin{split}
&\|[u_i,R_{ij}](\rho u_j)\|_{W^{1,4}}\le Cm_0^\frac{1}{12}\|u\|_{W^{1,6}}\|u\|_{L^\infty}\le Cm_0^\frac{1}{12}\|u\|_{W^{1,6}}^2\\
\le&Cm_0^\frac{1}{12}\left(\|\nabla
u\|_{L^2}^2+\|\mathrm{curl}u\|_{L^6}^2+\|\mathrm{div}u\|_{L^6}^2\right)\\
\le&Cm_0^\frac{1}{12}\left(\|\nabla
u\|_{L^2}^2+\|\nabla\mathrm{curl}u\|_{L^2}^2+\frac{1}{(2\mu+\lambda)^2}\|G\|_{L^6}^2+\frac{1}{(2\mu+\lambda)^2}\|\rho\theta\|_{L^6}^2\right)
\\
\le&Cm_0^\frac{1}{12}\left(\|\nabla
u\|_{L^2}^2+\|\nabla\mathrm{curl}u\|_{L^2}^2+\frac{1}{(2\mu+\lambda)^2}\|\nabla
G\|_{L^2}^2+\frac{1}{(2\mu+\lambda)^2}\|\nabla\theta\|_{L^2}^2\right).
\end{split}
\ee (\ref{RR1}) combined with the Sobolev inequality, (\ref{3d-H 1
of curl u}) and(\ref{3d-H 1 of G}) deduces
 \be\label{RR2}\begin{split}
&\|[u_i,R_{ij}](\rho
u_j)\|_{L^\infty}\\
\le&Cm_0^\frac{1}{12}\left(\|\nabla u\|_{L^2}^2+\frac{1}{\mu^2}
 \|\sqrt{\rho}\u\|_{L^{2}}^2+\frac{1}{(2\mu+\lambda)^2}\|\nabla\theta\|_{L^2}^2\right).
\end{split}
\ee
Substituting (\ref{RR2}) into $IV_3$, we have \be\label{3d-IV
3}\begin{split}
IV_3\le&\frac{Cm_0^\frac{1}{12}}{2\mu+\lambda}\int_0^t\left(\|\nabla
u\|_{L^2}^2+\frac{1}{\mu^2}
 \|\sqrt{\rho}\u\|_{L^{2}}^2+\frac{1}{(2\mu+\lambda)^2}\|\nabla\theta\|_{L^2}^2\right)\,d\tau\\
\le&
\frac{Cm_0^\frac{1}{12}}{2\mu+\lambda}\left(m_0^\frac{1}{2}+\frac{\tilde{E}}{\mu}
 +\frac{1}{(2\mu+\lambda)^2\kappa}\right).
\end{split}
\ee
 By (\ref{3d-IV 1}), (\ref{3d-IV 2}) and (\ref{3d-IV 3}), we
have \be\label{b(t)-b(0)}\begin{split}
b(t)-b(0)\le&\frac{Cm_0^\frac{1}{6}\sqrt{\tilde{E}}}{2\mu+\lambda}+\frac{Cm_0^\frac{1}{12}}{2\mu+\lambda}\left(m_0^\frac{1}{2}+\frac{\tilde{E}}{\mu}
 +\frac{1}{(2\mu+\lambda)^2\kappa}\right)\\
\le&\frac{Cm_0^\frac{1}{6}\sqrt{\tilde{E}}}{2\mu+\lambda}+\frac{Cm_0^\frac{7}{12}}{2\mu+\lambda}
+\frac{Cm_0^\frac{1}{12}\tilde{E}}{\mu(2\mu+\lambda)}
 +\frac{Cm_0^\frac{1}{12}}{(2\mu+\lambda)^3\kappa}\\
\le&\log\frac{3}{2},
\end{split}\ee provided that
{\small\bex
\begin{split}m_0\le&\min\left\{C_3,\ \frac{(2\mu+\lambda)^{6}(\log\frac{3}2)^{6}}{(4C)^6\tilde{E}^3},\
\frac{(2\mu+\lambda)^\frac{12}{7}(\log\frac{3}{2})^\frac{12}{7}}{(4
C)^\frac{12}{7}},\
\frac{\mu^{12}(2\mu+\lambda)^{12}(\log\frac{3}{2})^{12}}{(4C)^{12}\tilde{E}^{12}},\ \frac{(2\mu+\lambda)^{36}\kappa^{12}(\log\frac{3}{2})^{12}}{(4C)^{12}}\right\}\\
\triangleq& \varepsilon_0.\end{split} \eex }Integrating (\ref{rho
identity}) w.r.t. $s$ over $[0,t]$, we get \bex\begin{split}
\log\rho^\delta(x, t)=&
\log\left[\rho_0\left(X(t;x,0)\right)+\delta\right]+\int_0^tg(\tau)\,d\tau+b(t)-b(0)\\
\le&\log\left(\bar{\rho}+\delta\right)+\log\frac{3}{2},
\end{split}
\eex provided that $m_0\le \varepsilon_0$. Let
$\delta\rightarrow0^+$, we have
$$
\rho\le\frac{3\bar{\rho}}{2}.
$$

\endpf
\section{Asymptotic behavior in time}\label{3d-sec 5}

 \setcounter{equation}{0} \setcounter{theorem}{0}

In this section, we denote generic constants by $\bar{C}$ depending
on the initial data, coefficients of viscosity and heat conduction
and some other known constants but independent of $t$. Theorem
\ref{3d-th:1.2} will be proved in Sections \ref{3d-sec 5.1} and
\ref{3d-sec 5.2}.

\subsection{Large-time behavior}\label{3d-sec 5.1}

The main result in Section \ref{3d-sec 5.1} is stated as follows.
\begin{proposition}\label{3d-prop large time}
Under the conditions of Theorem \ref{3d-th:1.2}, it holds that
\be\label{3d-longtime u theta}
\int_{\mathbb{R}^3}\left(\rho|\theta|^2+|\nabla
u|^2+|\nabla\theta|^2\right)\rightarrow0, \ee as
$t\rightarrow\infty$.
\end{proposition}

\bigskip

To prove Proposition \ref{3d-prop large time}, we need some
estimates uniform for $t$. In fact, the lower order estimates of the
solutions have been made uniformly for $t$ in Section \ref{3d-sec
4}. More precisely,
 \begin{lemma}\label{3d-le:low uniform}
  Under the conditions of Theorem \ref{3d-th:1.2}, it holds that
 \be\label{3d-low uniform 1}0\le\rho\le\bar{C},\ee and
\be\label{3d-low uniform 2}
\int_{\mathbb{R}^3}\left(\rho+\rho|\theta|^2+|\nabla
u|^2\right)+\int_0^t\int_{\mathbb{R}^3}\left(|\nabla
u|^2+|\nabla\theta|^2+\rho|\u|^2\right)\le\bar{C}, \ee for any
$(x,t)\in\mathbb{R}^3\times[0,\infty)$.
 \end{lemma}

 With Lemma \ref{3d-le:low uniform}, one can follow the
 proofs of  Lemma \ref{blowup-le: int rho u t}, Corollary \ref{blowup-cor
 3.5} and Lemma \ref{blowup-le: H 2 of theta} step
by step, and easily get the following higher order estimates uniform
for $t$, respectively.
\begin{lemma}\label{3d-le: H 2 of u and H 1of theta}Under the conditions of Theorem \ref{3d-th:1.2}, it holds that
\be\label{3d-le: H 2 of u}
\int_{\mathbb{R}^3}\left(\rho|\u|^2+|\nabla
\theta|^2\right)+\int_0^t\int_{\mathbb{R}^3}\left(|\nabla\u|^2+\rho
|\dot{\theta}|^2\right)\le \bar{C}, \ee for any $t\in[0,\infty)$.
\end{lemma}

\begin{corollary}\label{corollary-long time}
Under the conditions of Theorem \ref{3d-th:1.2}, it holds that
\be\label{3d-u infty} \|u\|_{L^\infty}+\|\nabla
u\|_{L^r}+\int_0^t\int_{\mathbb{R}^3}\left(\rho|\theta_t|^2+|\nabla^2\theta|^2\right)\le\bar{C},
\ee for any $r\in[2,6]$, and any $t\in[0,\infty)$.
\end{corollary}

\begin{lemma}\label{3d-le: H 2 of theta} Under the conditions of Theorem \ref{3d-th:1.2}, it holds that
\be\label{3d-H 2 of theta}
\int_{\mathbb{R}^3}\rho|\theta_t|^2+\int_0^t\int_{\mathbb{R}^3}|\nabla\theta_t|^2\le
\bar{C}, \ee for any $t\in[0,\infty)$.
\end{lemma}

\bigskip

{\noindent\bf \underline{Proof of Proposition \ref{3d-prop large time}��}}\\

Denote $$ F(t)=\int_{\mathbb{R}^3}\left(\mu|\mathrm{curl}
u|^2+\frac{G^2}{2\mu+\lambda}\right).
$$ By (\ref{3d-low uniform 1}) and (\ref{3d-low uniform 2}), we have
\be\label{3d-F in L 1}F\in L^1(0,\infty). \ee Moreover, by
(\ref{3d-dt nabla u 2-1}), (\ref{3d-II 3+II 4}), (\ref{3d-low
uniform 1}) and (\ref{3d-low uniform 2}), we have
\beq\label{3d-longtime u}\begin{split}
|\frac{d}{dt}F(t)|\le&\bar{C}\int_{\mathbb{R}^3}\rho |\u|^2+
\bar{C}\|\rho\|_{L^6}^2\|\nabla\theta\|_{L^2}^2\|\nabla
u\|_{L^2}^2+\bar{C}\|\nabla u\|_{L^2}^6\\&+\bar{C}\|\nabla
u\|_{L^2}^4\|\rho\theta\|_{L^2}^2+\bar{C}\|\nabla\theta\|_{L^2}^2\\
\le&\bar{C}\int_{\mathbb{R}^3}\rho |\u|^2+
\bar{C}\|\nabla\theta\|_{L^2}^2+\bar{C}\|\nabla u\|_{L^2}^2,
\end{split}
\eeq where we have used $\Delta
u=\nabla\mathrm{div}u-\nabla\times(\mathrm{curl}u)$ such that
\be\label{3d-nabla curl div} \int_{\mathbb{R}^3}|\nabla
u|^2=\int_{\mathbb{R}^3}(|\mathrm{div}u|^2+|\mathrm{curl}u|^2). \ee
By (\ref{3d-longtime u}), (\ref{3d-low uniform 2}) and (\ref{3d-F in
L 1}), we conclude that
$$
F\in W^{1,1}(0,\infty),
$$ which deduces
\be\label{3d-longtime nabla u}
\int_{\mathbb{R}^3}\left(\mu|\mathrm{curl}
u|^2+\frac{G^2}{2\mu+\lambda}\right)(t)=F(t)\rightarrow0, \ee as
$t\rightarrow\infty$.

It follows from (\ref{3d-low uniform 2}) and (\ref{3d-H 2 of theta})
that \bex \|\nabla\theta\|_{L^2}^2(\cdot)\in W^{1,1}(0,\infty), \eex
which deduces \be\label{3d-longtime theta}
\|\nabla\theta\|_{L^2}^2(t)\rightarrow0, \ee as
$t\rightarrow\infty$.

 By (\ref{3d-longtime nabla u}), (\ref{3d-longtime theta}), (\ref{3d-nabla curl div}), (\ref{3d-low uniform 1}) and (\ref{3d-low uniform 2}), we
get (\ref{3d-longtime u theta}).\\

 \subsection{Decay estimates}\label{3d-sec 5.2}

\begin{proposition}
Under the conditions of Theorem \ref{3d-th:1.2}, we get
\be\label{3d-decay rate}
\int_{\mathbb{R}^3}\left(\rho\theta^2+|\nabla u|^2\right)\le
\bar{C}\exp\{-\bar{C}_1t\}, \ee for any $t\in[0,\infty)$, provided that
$$
m_0\le\min\{\varepsilon_0,\tilde{\varepsilon}_0\},
$$ for some $\tilde{\varepsilon}_0>0$ depending on $\mu,\lambda,\kappa, K_1, K_2,
\bar{\rho}$, and some other known constants but independent of $t$.
\end{proposition}
\begin{remark}
The decay rate of $\|\nabla\theta\|_{L^2}$ is still unknown.
\end{remark}
\pf
 By (\ref{3d-dt
nabla u 2}), (\ref{3d-nabla curl div}), Corollary \ref{3d-cor 3.2}
and (\ref{3d-mass conservation}), we have
  \be\label{3d-asyu}\begin{split} &\int_{\mathbb{R}^3}\rho
|\u|^2+\frac{d}{dt}\int_{\mathbb{R}^3}\left(\mu|\mathrm{curl}
u|^2+\frac{G^2}{2\mu+\lambda}\right)\\
\le&
\frac{C}{(2\mu+\lambda)^2}\|\rho\|_{L^6}^2\|\nabla\theta\|_{L^2}^2\|\nabla
u\|_{L^2}^2+\frac{C}{\mu^2}\|\nabla
u\|_{L^2}^6+\frac{C}{(2\mu+\lambda)^4}\|\nabla
u\|_{L^2}^4\|\rho\theta\|_{L^2}^2 \\&+\frac{2\bar{\rho}\kappa^2
}{(2\mu+\lambda)^2}\|\nabla\theta\|_{L^2}^2\\
\le&
\frac{C\tilde{E}m_0^\frac{1}{3}}{(2\mu+\lambda)^2}\|\nabla\theta\|_{L^2}^2+\frac{C\tilde{E}^2}{\mu^2}\|\nabla
u\|_{L^2}^2+\frac{C\tilde{E}}{(2\mu+\lambda)^4}\|\nabla u\|_{L^2}^2
+\frac{C\kappa^2
}{(2\mu+\lambda)^2}\|\nabla\theta\|_{L^2}^2\\
=& M_1\|\nabla\theta\|_{L^2}^2+M_2\|\nabla u\|_{L^2}^2,
\end{split}
\ee where
$M_1=\left(\frac{C\tilde{E}m_0^\frac{1}{3}}{(2\mu+\lambda)^2}+\frac{C\kappa^2
}{(2\mu+\lambda)^2}\right)$, and
$M_2=\left(\frac{C\tilde{E}^2}{\mu^2}+\frac{C\tilde{E}}{(2\mu+\lambda)^4}\right)$.

By (\ref{3d-dt rho theta+1}), we have
\be\label{3d-asytheta}\begin{split}
\kappa\int_{\mathbb{R}^3}|\nabla\theta|^2+\frac{d}{dt}\int_{\mathbb{R}^3}\rho|\theta|^2
\le& \frac{C(2\mu+\lambda)^2}{\kappa\mu}\|\nabla
u\|_{L^2}^3\|\sqrt{\rho}\u\|_{L^2}+\frac{C}{\kappa}\|\nabla
u\|_{L^2}^2\|\rho\theta\|_{L^2}\|\sqrt{\rho}\u\|_{L^2}
\\&+C\|\nabla
u\|_{L^2}\|\nabla\theta\|_{L^2}^2\|\rho\|_{L^6}.
\end{split}
\ee Multiplying (\ref{3d-asytheta}) by
$\displaystyle\frac{2M_1}{\kappa}$, and using Corollary \ref{3d-cor
3.2}, (\ref{3d-mass conservation}), and Cauchy inequality, we have
 \be\label{3d-asytheta+1}\begin{split}
&2M_1\int_{\mathbb{R}^3}|\nabla\theta|^2+\frac{2M_1}{\kappa}\frac{d}{dt}\int_{\mathbb{R}^3}\rho|\theta|^2\\
\le& \frac{C(2\mu+\lambda)^2\tilde{E}M_1}{\kappa^2\mu}\|\nabla
u\|_{L^2}\|\sqrt{\rho}\u\|_{L^2}+\frac{CM_1\sqrt{\tilde{E}}}{\kappa^2}\|\nabla
u\|_{L^2}\|\sqrt{\rho}\u\|_{L^2}
\\&+\frac{Cm_0^\frac{1}{6}M_1\sqrt{\tilde{E}}}{\kappa}\|\nabla\theta\|_{L^2}^2\\
\le&
\frac{1}{2}\|\sqrt{\rho}\u\|_{L^2}^2+\left(\frac{C(2\mu+\lambda)^4\tilde{E}^2M_1^2}{\kappa^4\mu^2}+\frac{CM_1^2\tilde{E}}{\kappa^4}\right)\|\nabla
u\|_{L^2}^2
\\&+\frac{Cm_0^\frac{1}{6}M_1\sqrt{\tilde{E}}}{\kappa}\|\nabla\theta\|_{L^2}^2.
\end{split}
\ee Adding (\ref{3d-asytheta+1}) into (\ref{3d-asyu}), we have
\bex\begin{split}
&M_1\int_{\mathbb{R}^3}|\nabla\theta|^2+\frac{2M_1}{\kappa}\frac{d}{dt}\int_{\mathbb{R}^3}\rho|\theta|^2+\frac{1}{2}\int_{\mathbb{R}^3}\rho
|\u|^2+\frac{d}{dt}\int_{\mathbb{R}^3}\left(\mu|\mathrm{curl}
u|^2+\frac{G^2}{2\mu+\lambda}\right)\\
\le& M_2\|\nabla
u\|_{L^2}^2+\left(\frac{C(2\mu+\lambda)^4\tilde{E}^2M_1^2}{\kappa^4\mu^2}+\frac{CM_1^2\tilde{E}}{\kappa^4}\right)\|\nabla
u\|_{L^2}^2
+\frac{Cm_0^\frac{1}{6}M_1\sqrt{\tilde{E}}}{\kappa}\|\nabla\theta\|_{L^2}^2\\
\le& M_3\|\nabla u\|_{L^2}^2 +\frac{M_1}{2}\|\nabla\theta\|_{L^2}^2,
\end{split}
\eex provided that
$$
m_0\le\min\left\{\varepsilon_0,
\frac{\kappa^6}{2^6C^6\tilde{E}^3}\right\}.
$$Here
$\displaystyle
M_3=M_2+\frac{C(2\mu+\lambda)^4\tilde{E}^2M_1^2}{\kappa^4\mu^2}+\frac{CM_1^2\tilde{E}}{\kappa^4}$.

Thus, \be\label{3d-asytheta u}\begin{split}
\int_{\mathbb{R}^3}\left(\frac{M_1}{2}|\nabla\theta|^2+\frac{1}{2}\rho
|\u|^2\right)+\frac{d}{dt}\int_{\mathbb{R}^3}\left(\frac{2M_1}{\kappa}\rho|\theta|^2+\mu|\mathrm{curl}
u|^2+\frac{G^2}{2\mu+\lambda}\right) \le M_3\|\nabla u\|_{L^2}^2.
\end{split}
\ee Multiplying (\ref{3d-dt rho u}) by $\displaystyle
\frac{2M_3}{\mu}$, and adding the resulting inequality into
(\ref{3d-asytheta u}), we have \be\label{3d-asytheta
u+1}\begin{split} &\int_{\mathbb{R}^3}\left(M_3|\nabla
u|^2+\frac{M_1}{2}|\nabla\theta|^2+\frac{1}{2}\rho
|\u|^2\right)\\&+\frac{d}{dt}\int_{\mathbb{R}^3}\left(\frac{M_3}{\mu}\rho|u|^2+\frac{2M_1}{\kappa}\rho|\theta|^2+\mu|\mathrm{curl}
u|^2+\frac{G^2}{2\mu+\lambda}\right)\\ \le&
\frac{Cm_0^\frac{2}{3}M_3}{\mu(\mu+\lambda)}\|\nabla\theta\|_{L^2}^2\le
\frac{M_1}{4}\|\nabla\theta\|_{L^2}^2,
\end{split}
\ee provided that
$$
m_0\le\min\left\{\varepsilon_0, \tilde{\varepsilon}_0\right\},
$$ where $$\tilde{\varepsilon}_0=\min\left\{\frac{\kappa^6}{2^6C^6\tilde{E}^3},\
\frac{\mu^\frac{3}{2}(\mu+\lambda)^\frac{3}{2}M_1^\frac{3}{2}}{8C^\frac{3}{2}M_3^\frac{3}{2}}\right\}.$$
By (\ref{3d-asytheta u+1}), together with the facts \bex
\int_{\mathbb{R}^3}\rho|u|^2\le
\|\rho\|_{L^\frac{3}{2}}\|u\|_{L^6}^2\le \bar{C}\|\nabla
u\|_{L^2}^2,\eex  \bex \int_{\mathbb{R}^3}\rho|\theta|^2\le
\bar{C}\|\nabla \theta\|_{L^2}^2,\eex  \bex
\int_{\mathbb{R}^3}|\mathrm{curl}u|^2\le\bar{C}\|\nabla\mathrm{curl}u\|_{L^\frac{6}{5}}^2\le\bar{C}\|\rho\u\|_{L^\frac{6}{5}}^2
\le\bar{C}\|\sqrt{\rho}\u\|_{L^2}^2\|\sqrt{\rho}\|_{L^3}^2\le\bar{C}\|\sqrt{\rho}\u\|_{L^2}^2,\eex
and \bex
\int_{\mathbb{R}^3}|G|^2\le\bar{C}\|\sqrt{\rho}\u\|_{L^2}^2,\eex
 we get
\be\label{3d-asytheta u+2}\begin{split}
&\bar{C}_1\int_{\mathbb{R}^3}\left(\frac{M_3}{\mu}\rho|u|^2+\frac{2M_1}{\kappa}\rho|\theta|^2+\mu|\mathrm{curl}
u|^2+\frac{G^2}{2\mu+\lambda}\right)\\&+\frac{d}{dt}\int_{\mathbb{R}^3}\left(\frac{M_3}{\mu}\rho|u|^2+\frac{2M_1}{\kappa}\rho|\theta|^2+\mu|\mathrm{curl}
u|^2+\frac{G^2}{2\mu+\lambda}\right)\le0,
\end{split}
\ee for some constant $\bar{C}_1>0$ depending on $\mu,\lambda,
\kappa, M_1, M_3$ and other known constants but independent of $t$.
(\ref{3d-asytheta u+2}) deduces
\be\label{3d-asyresult}\int_{\mathbb{R}^3}\left(\frac{M_3}{\mu}\rho|u|^2+\frac{2M_1}{\kappa}\rho|\theta|^2+\mu|\mathrm{curl}
u|^2+\frac{G^2}{2\mu+\lambda}\right)\le A\exp\{-\bar{C}_1t\},\ee
where $$
A=\int_{\mathbb{R}^3}\left(\frac{M_3}{\mu}\rho_0|u_0|^2+\frac{2M_1}{\kappa}\rho_0|\theta_0|^2+\mu|\mathrm{curl}
u_0|^2+\frac{G_0^2}{2\mu+\lambda}\right).
$$
By (\ref{3d-nabla curl div}) and (\ref{3d-asyresult}), we get
(\ref{3d-decay rate}).
\endpf

\section{Global classical solution}
\label{3d-sec 7}

\setcounter{equation}{0} \setcounter{theorem}{0}

The proof of local existence and uniqueness of the classical
solution as in Theorem \ref{3d-th1.3} can be found in Section
\ref{3d-sec 8} (see Appendix B below). Let $T^*_1>0$ be the maximal existence
time of the classical solution. Our aim is to
prove $T^*_1=\infty$. To do this, we use contradiction arguments
similar to Section \ref{3d-sec 4}.

More precisely, we assume that $0<T_1^*<\infty$. In this section, we
denote generic constants by $\check{C}$ depending on the initial
data, $\mu,\lambda,\kappa$, $T_1^*$ and some other known constants
but independent of $t\in[0,T_1^*)$. In this case, we shall prove
\be\label{3d-classicalaim1} \|\rho(\cdot,t)\|_{H^2\cap W^{2,q}}\le
\check{C}, \ee for any $t\in[0,T_1^*)$, and
\be\label{3d-classicalaim2}
 \|\nabla
 u(\cdot,t)\|_{H^1}+\|\nabla\theta(\cdot,t)\|_{H^1}\le\check{C},\ee for any $t\in[0,T_1^*)$, and
\be\label{3d-classicalaim3} \|\sqrt{\rho}
u_t(\cdot,t)\|_{L^2}+\|\sqrt{\rho}
\theta_t(\cdot,t)\|_{L^2}\le\check{C},\ee for a.e. $t\in[0,T_1^*)$.
With (\ref{3d-classicalaim1}), (\ref{3d-classicalaim2}) and
(\ref{3d-classicalaim3}), we can define a new initial data at
$T_1^*$ \be \label{3d-redd.intial} \big(\rho(\cdot, T_1^*), u(\cdot,
T_1^*), \theta(\cdot, T_1^*)\big)=\lim\limits_{t\nearrow
T_1^*}\big(\rho(\cdot, t), u(\cdot, t), \theta(\cdot, t)\big), \ee
which satisfies the conditions of Appendix B. This means that the
life span of the classical solution beyond $T_1^*$, which is the
desired contradiction.

  To get (\ref{3d-classicalaim1}), (\ref{3d-classicalaim2}) and
(\ref{3d-classicalaim3}), we begin with the following lemma which is
essentially obtained in Section \ref{3d-sec 4}.

\begin{lemma}\label{3d-le:6.1} Under the conditions of Theorem \ref{3d-th1.3}, it holds that
\be\label{3d-rho}
\|\rho(\cdot,t)\|_{L^\infty}+\|\rho(\cdot,t)\|_{H^1\cap
W^{1,q}}+\|\rho_t(\cdot,t)\|_{L^2\cap L^q}\le\check{C}, \ee
\be\label{3d-u theta}
 \|u(\cdot,t)\|_{L^\infty}+\|\nabla
 u(\cdot,t)\|_{H^1}+\|\nabla\theta(\cdot,t)\|_{H^1}+\|\nabla^2u\|_{L^2([0,t];L^q)}+\|\nabla^2\theta\|_{L^2([0,t];L^q)}\le\check{C},\ee
\be\label{3d-u t theta t}\|\nabla u_t\|_{L^2([0,t];L^2)}+\|\nabla
\theta_t\|_{L^2([0,t];L^2)}+ \|\sqrt{\rho}
u_t(\cdot,t)\|_{L^2}+\|\sqrt{\rho}
\theta_t(\cdot,t)\|_{L^2}\le\check{C},\ee for a.e. $t\in[0,T_1^*)$.
\end{lemma}

From Lemma \ref{3d-le:6.1}, (\ref{3d-classicalaim2}) and
(\ref{3d-classicalaim3}) have been obtained. What we need to do is
to get (\ref{3d-classicalaim1}).
\begin{lemma}\label{3d-le: H 2 of rho}
Under the conditions of Theorem \ref{3d-th1.3}, it holds that
 \be\label{3d-H 2 of rho}
\|\rho(\cdot,t)\|_{H^2}+\int_0^t\|\nabla
u(\cdot,t)\|_{H^2}^2\,dt\le\check{C},
 \ee for any $t\in[0,T_1^*)$.
\end{lemma}
\pf Taking $\nabla^2$ on both sides of (\ref{full N-S+1})$_1$, we
have \be\label{3d-equation of H 2 rho}\begin{split}
&\nabla^2\rho_t+2\nabla
u^j\otimes\nabla\nabla_j\rho+u^j\nabla^2\nabla_j\rho+\nabla^2u^j\nabla_j\rho+\nabla^2\rho\mathrm{div}u+2\nabla\rho\otimes\nabla\mathrm{div}u\\&+\rho\nabla^2\mathrm{div}u=0.
\end{split}
\ee Multiplying (\ref{3d-equation of H 2 rho}) by $\nabla^2\rho$,
integrating by parts over $\mathbb{R}^3$, and using the Cauchy
inequality, the H\"older inequality, the Sobolev inequality and
(\ref{3d-rho}), we have \be\label{3d-dt rho 2}\begin{split}
\frac{d}{dt}\int_{\mathbb{R}^3}|\nabla^2\rho|^2\le&\check{C}\|\nabla
u\|_{L^\infty}\int_{\mathbb{R}^3}|\nabla^2\rho|^{2}+\check{C}\int_{\mathbb{R}^3}|\nabla^2\rho||\nabla^2u||\nabla\rho|
+\check{C}\int_{\mathbb{R}^3}|\nabla^2\rho||\nabla^3u|\\
\le&\check{C}(\|\nabla
u\|_{L^\infty}+1)\int_{\mathbb{R}^3}|\nabla^2\rho|^{2}
+\check{C}\|\nabla^2u\|_{L^6}^2\|\nabla\rho\|_{L^3}^2+\check{C}\int_{\mathbb{R}^3}|\nabla^3u|^2\\
\le&\check{C}(\|\nabla
u\|_{L^\infty}+1)\int_{\mathbb{R}^3}|\nabla^2\rho|^{2}
+\check{C}\int_{\mathbb{R}^3}|\nabla^3u|^2.
\end{split}
\ee By (\ref{full N-S+1})$_2$ and the standard elliptic estimates,
we have \be\label{3d-nabla 3 u}\begin{split}
\|\nabla^3u\|_{L^2}^2\le&\check{C}\|\nabla\rho\|_{L^3}^2\|u_t\|_{L^6}^2+\check{C}\|\nabla
u_t\|_{L^2}^2+\check{C}\|\nabla\rho\|_{L^3}^2\|\nabla
u\|_{L^6}^2+\check{C}\|\nabla
u\|_{L^4}^4\\&+\check{C}\|\nabla^2u\|_{L^2}^2+\check{C}\|\nabla^2\rho\|_{L^2}^2
+\check{C}\|\nabla\rho\|_{L^3}^2\|\nabla\theta\|_{L^6}^2+\check{C}\|\nabla^2\theta\|_{L^2}^2\\
\le&\check{C}\|\nabla
u_t\|_{L^2}^2+\check{C}\|\nabla^2\rho\|_{L^2}^2+\check{C},
\end{split}
\ee where we have used the H\"older inequality, the Sobolev
inequality, (\ref{3d-rho}) and (\ref{3d-u theta}).

Substituting (\ref{3d-nabla 3 u}) into (\ref{3d-dt rho 2}), and
using (\ref{3d-rho}), (\ref{3d-u theta}) and the Sobolev inequality
again, together with (\ref{3d-u t theta t}) and the Gronwall
inequality, we have \be \label{3d-nabla 2 rho}\begin{split}
\|\rho(\cdot,t)\|_{H^2} \le\check{C},
\end{split}
\ee for any $t\in[0,T_1^*)$. By (\ref{3d-u theta}), (\ref{3d-u t
theta t}), (\ref{3d-nabla 3 u}), (\ref{3d-nabla 2 rho}), we get
$$
\int_0^t\|\nabla u(\cdot,t)\|_{H^2}^2dt\le \check{C},
$$ for any $t\in[0,T_1^*)$.
This completes the proof of Lemma \ref{3d-le: H 2 of rho}.
\endpf
\begin{corollary}\label{3d-cor rho tt}
Under the conditions of Theorem \ref{3d-th1.3}, it holds that
\be\label{3d-rhott}
\|\rho_t(\cdot,t)\|_{H^1}+\int_0^t\|\rho_{ss}\|_{L^2}^2\le\check{C},
\ee for a.e. $t\in[0,T_1^*)$.
\end{corollary}
\pf Taking $\nabla$ on both sides of (\ref{full N-S+1})$_1$, we have
$$
\nabla\rho_t=-\nabla(\rho\mathrm{div}u+u\cdot\nabla\rho)=-\nabla\rho\mathrm{div}u-\rho\nabla\mathrm{div}u-\nabla
u\cdot\nabla\rho-u\cdot\nabla\nabla\rho.
$$
This, together with (\ref{3d-rho}), the H\"older inequality, the
Sobolev inequality, (\ref{3d-u theta}) and (\ref{3d-H 2 of rho}),
deduces \be\label{3d-rho t H 1}\begin{split}
\|\rho_t(\cdot,t)\|_{H^1}=&\|\rho_t(\cdot,t)\|_{L^2}+\|\nabla\rho_t(\cdot,t)\|_{L^2}\\
\le&\check{C}+\check{C}\|\nabla\rho(\cdot,t)\|_{L^3}\|\mathrm{div}u(\cdot,t)\|_{L^6}+\check{C}\|\rho(\cdot,t)\|_{L^\infty}\|\nabla\mathrm{div}u(\cdot,t)\|_{L^2}
\\&+\check{C}\|\nabla
u(\cdot,t)\|_{L^6}\|\nabla\rho(\cdot,t)\|_{L^3}+\check{C}\|u(\cdot,t)\|_{L^\infty}\|\nabla^2\rho(\cdot,t)\|_{L^2}\\
\le&\check{C}.
\end{split}\ee
Using (\ref{full N-S+1})$_1$ again, similar to (\ref{3d-rho t H 1}),
we have \be\label{3d-rho tt}\begin{split}
\|\rho_{tt}\|_{L^2}=&\|(\rho\mathrm{div}u+u\cdot\nabla\rho)_t\|_{L^2}=\|\rho_t\mathrm{div}u+\rho\mathrm{div}u_t+u_t\cdot\nabla\rho+u\cdot\nabla\rho_t\|_{L^2}
\\ \le&\|\rho_t\|_{L^6}\|\mathrm{div}u\|_{L^3}+\|\rho\|_{L^\infty}\|\mathrm{div}u_t\|_{L^2}+\|u_t\|_{L^6}\|\nabla\rho\|_{L^3}+\|u\|_{L^\infty}\|\nabla\rho_t\|_{L^2}
\\ \le&\check{C}\|\nabla u_t\|_{L^2}+\check{C}.
\end{split}
\ee By (\ref{3d-rho tt}) and (\ref{3d-u t theta t}), we get
$$
\int_0^t\|\rho_{tt}\|_{L^2}^2\le\check{C}.
$$ for a.e. $t\in[0,T_1^*)$.
 The proof of Corollary \ref{3d-cor rho tt} is
complete.
\endpf
\begin{lemma}\label{3d-le:6.4}
Under the conditions of Theorem \ref{3d-th1.3}, it holds that \be
 \label{3d-na ut} \int_{\mathbb{R}^3}t|\nabla
u_t|^2+\int_0^t\int_{\mathbb{R}^3}s\rho |u_{ss}|^2\le \check{C}, \ee
for a.e. $t\in[0,T_1^*)$.
\end{lemma}
\pf Differentiating (\ref{full N-S+1})$_2$ w.r.t. $t$, we have
\be\label{rho u tt+.=} \rho u_{tt}+\rho_tu_t+\rho_t u\cdot\nabla
u+\rho u_t\cdot\nabla u+\rho u\cdot\nabla u_t+\nabla P_t=\mu\Delta
u_t+(\mu+\lambda)\nabla\mathrm{div}u_t. \ee Multiplying (\ref{rho u
tt+.=}) by $u_{tt}$, integrating by parts over $\mathbb{R}^3$, we
have \be\label{3d-dt nab ut}\begin{split}
&\int_{\mathbb{R}^3}\rho|u_{tt}|^2+\frac{1}{2}\frac{d}{dt}\int_{\mathbb{R}^3}\left(\mu|\nabla
u_t|^2+(\mu+\lambda)|\mathrm{div}u_t|^2\right)\\
&=-\int_{\mathbb{R}^3}(\rho_tu_t+\rho_t u\cdot\nabla u)\cdot
u_{tt}-\int_{\mathbb{R}^3}(\rho u_t\cdot\nabla u+\rho u\cdot\nabla
u_t)\cdot u_{tt}-\int_{\mathbb{R}^3}\nabla P_t\cdot u_{tt}\\
&=\sum\limits_{i=1}^3VI_i.
\end{split}
\ee For $VI_1$, using (\ref{full N-S+1})$_1$, integration by parts,
the H\"older inequality, the Sobolev inequality, (\ref{3d-rho}),
(\ref{3d-u theta}), (\ref{3d-u t theta t}) and (\ref{3d-rhott}), we
have \be\label{3d-VI 1}\begin{split}
VI_1=&-\frac{d}{dt}\int_{\mathbb{R}^3}\left(\frac{1}{2}\rho_t|u_t|^2+\rho_t
(u\cdot\nabla) u\cdot
u_{t}\right)+\frac{1}{2}\int_{\mathbb{R}^3}\rho_{tt}|u_t|^2\\&+\int_{\mathbb{R}^3}(\rho_{tt}
u\cdot\nabla u+\rho_t u_t\cdot\nabla u+\rho_t u\cdot\nabla u_t)\cdot
u_{t}\\
\le&-\frac{d}{dt}\int_{\mathbb{R}^3}\left(\frac{1}{2}\rho_t|u_t|^2+\rho_t
(u\cdot\nabla) u\cdot u_{t}\right)+\int_{\mathbb{R}^3}(\rho_{t}
u+\rho u_t)\cdot\nabla u_t\cdot u_t\\&+\|\rho_{tt}\|_{L^2}
\|u\|_{L^\infty}\|\nabla
u\|_{L^3}\|u_{t}\|_{L^6}+\|\rho_t\|_{L^2}\|u_t\|_{L^6}^2\|\nabla
u\|_{L^6}\\&+\|\rho_t\|_{L^6}\|u\|_{L^6}\|\nabla
u_t\|_{L^2}\|u_{t}\|_{L^6}\\
\le&-\frac{d}{dt}\int_{\mathbb{R}^3}\left(\frac{1}{2}\rho_t|u_t|^2+\rho_t
(u\cdot\nabla) u\cdot
u_{t}\right)+\|\rho_{t}\|_{L^3}\|u\|_{L^\infty}\|\nabla
u_t\|_{L^2}\|u_t\|_{L^6}\\&+
\|\sqrt{\rho}\|_{L^\infty}\|u_t\|_{L^\infty}\|\nabla
u_t\|_{L^2}\|\sqrt{\rho} u_t\|_{L^2}+\check{C}\|\nabla
u_t\|_{L^2}^2+\check{C}\|\rho_{tt}\|_{L^2}^2\\
\le&-\frac{d}{dt}\int_{\mathbb{R}^3}\left(\frac{1}{2}\rho_t|u_t|^2+\rho_t
(u\cdot\nabla) u\cdot u_{t}\right)+\check{C}(\|\nabla
u_t\|_{L^2}+\|\nabla^2u_t\|_{L^2})\|\nabla
u_t\|_{L^2}\\&+\check{C}\|\rho_{tt}\|_{L^2}^2.
\end{split}
\ee For $VI_2$, using the H\"older inequality, the Sobolev
inequality, (\ref{3d-rho}) and (\ref{3d-u theta}) again, we have
\be\label{3d-VI 2}\begin{split}
VI_2\le&\frac{1}{4}\int_{\mathbb{R}^3}\rho|u_{tt}|^2+\check{C}\|\rho\|_{L^\infty}
\|u_t\|_{L^6}^2\|\nabla
u\|_{L^3}^2+\check{C}\|\rho\|_{L^\infty}\|u\|_{L^\infty}^2\|\nabla
u_t\|_{L^2}^2\\
\le&\frac{1}{4}\int_{\mathbb{R}^3}\rho|u_{tt}|^2+\check{C} \|\nabla
u_t\|_{L^2}^2.
\end{split}
\ee For $VI_3$, we have \be\label{3d-VI 3}\begin{split}
VI_3=&\frac{d}{dt}\int_{\mathbb{R}^3} P_t\mathrm{div}
u_{t}-\int_{\mathbb{R}^3}
P_{tt}\mathrm{div}u_{t}\\=&\frac{d}{dt}\int_{\mathbb{R}^3}
P_t\mathrm{div} u_{t}-\frac{1}{2\mu+\lambda}\int_{\mathbb{R}^3}
P_{tt}G_t-\frac{1}{2(2\mu+\lambda)}\frac{d}{dt}\int_{\mathbb{R}^3}
P_{t}^2.
\end{split}
\ee By (\ref{full N-S+1})$_1$ and (\ref{full N-S+1})$_3$, we have
\bex\begin{split} P_{tt}=&-\mathrm{div}(\rho\theta
u)_t-(\rho\theta\mathrm{div}u)_t+\mu\left(\nabla u+(\nabla
u)^\prime\right):\left(\nabla u_t+(\nabla
u_t)^\prime\right)\\&+2\lambda(\mathrm{div}u)\mathrm{div}u_t+\kappa\Delta\theta_t.
\end{split}\eex
Substituting this equality into the second term of the right side of
(\ref{3d-VI 3}), and using integration by parts, the H\"older
inequality, the Sobolev inequality, (\ref{3d-rho}), (\ref{3d-u
theta}) and (\ref{3d-u t theta t}), we have
\be\label{3d-PttGt}\begin{split}-\frac{1}{2\mu+\lambda}\int_{\mathbb{R}^3}
P_{tt}G_t=&-\frac{1}{2\mu+\lambda}\int_{\mathbb{R}^3}(\rho\theta
u)_t\cdot\nabla G_t+\frac{1}{2\mu+\lambda}\int_{\mathbb{R}^3}
G_t(\rho\theta\mathrm{div}u)_t\\&-\frac{\mu}{2\mu+\lambda}\int_{\mathbb{R}^3}
G_t\left(\nabla u+(\nabla u)^\prime\right):\left(\nabla u_t+(\nabla
u_t)^\prime\right)\\&-\frac{2\lambda}{2\mu+\lambda}\int_{\mathbb{R}^3}
G_t(\mathrm{div}u)\mathrm{div}u_t+\frac{\kappa}{2\mu+\lambda}\int_{\mathbb{R}^3}\nabla
G_t\cdot\nabla\theta_t\\
\le&\check{C}\|(\rho\theta u)_t\|_{L^2}\|\nabla
G_t\|_{L^2}+\check{C}\|G_t\|_{L^2}\|\rho\theta\mathrm{div}u_t\|_{L^2}
\\&+\check{C}\|G_t\|_{L^6}\|\rho\theta_t\mathrm{div}u\|_{L^\frac{6}{5}}
+\check{C}\|G_t\|_{L^6}\|\rho_t\theta\mathrm{div}u\|_{L^\frac{6}{5}}\\&+\check{C}\|G_t\|_{L^6}\|\nabla
u\|_{L^3}\|\nabla u_t\|_{L^2}+\check{C}\|\nabla
G_t\|_{L^2}\|\nabla\theta_t\|_{L^2}\\
\le&\check{C}\|\nabla G_t\|_{L^2}(\|\nabla
u_t\|_{L^2}+\|\nabla\theta_t\|_{L^2}+1)+\check{C}\|\nabla
u_t\|_{L^2}^2+\check{C}.
\end{split}\ee
It follows from (\ref{equation of G}) that
\be\label{3d-nablaGt}\begin{split} \|\nabla
G_t\|_{L^2}\le&\check{C}\|(\rho u_t+\rho u\cdot\nabla
u)_t\|_{L^2}\\
\le&\check{C}\|\rho_t\|_{L^3}\|u_t\|_{L^6}+\check{C}\|\sqrt{\rho}\|_{L^\infty}\|\sqrt{\rho}
u_{tt}\|_{L^2}+\check{C}\|\rho_t\|_{L^3}\|u\|_{L^\infty}\|\nabla
u\|_{L^6}\\&+\check{C}\|\rho\|_{L^\infty}\|u_t\|_{L^6}\|\nabla
u\|_{L^3}+\check{C}\|\rho\|_{L^\infty}\|u\|_{L^\infty}\|\nabla
u_t\|_{L^2}\\
\le&\check{C}\|\nabla u_t\|_{L^2}+\check{C}\|\sqrt{\rho}
u_{tt}\|_{L^2}+\check{C},
\end{split}
\ee where we have used the H\"older inequality, the Sobolev
inequality, (\ref{3d-rho}) and (\ref{3d-u theta}).

 Substituting (\ref{3d-nablaGt}) into (\ref{3d-PttGt}), and using
the Cauchy inequality, we have
\be\label{3d-PttGt+1}\begin{split}-\frac{1}{2\mu+\lambda}\int_{\mathbb{R}^3}
P_{tt}G_t
\le&\frac{1}{4}\int_{\mathbb{R}^3}\rho|u_{tt}|^2+\check{C}\|\nabla
u_t\|_{L^2}^2+\check{C}\|\nabla\theta_t\|_{L^2}^2+\check{C}.
\end{split}\ee
Putting (\ref{3d-VI 1}), (\ref{3d-VI 2}), (\ref{3d-VI 3}) and
(\ref{3d-PttGt+1}) into (\ref{3d-dt nab ut}), and using the Cauchy
inequality, we have \be\label{3d-dt nab ut+1}\begin{split}
&\frac{1}{2}\int_{\mathbb{R}^3}\rho|u_{tt}|^2+\frac{1}{2}\frac{d}{dt}\int_{\mathbb{R}^3}\left(\mu|\nabla
u_t|^2+(\mu+\lambda)|\mathrm{div}u_t|^2\right)\\
\le&-\frac{d}{dt}\int_{\mathbb{R}^3}\left(\frac{1}{2}\rho_t|u_t|^2+\rho_t
(u\cdot\nabla) u\cdot u_{t}\right)+\frac{d}{dt}\int_{\mathbb{R}^3}
P_t\mathrm{div}
u_{t}-\frac{1}{2(2\mu+\lambda)}\frac{d}{dt}\int_{\mathbb{R}^3}
P_{t}^2\\&+\epsilon\|\nabla^2u_t\|_{L^2}^2+\check{C}_\epsilon\|\nabla
u_t\|_{L^2}^2+\check{C}\|\rho_{tt}\|_{L^2}^2+\check{C}\|\nabla\theta_t\|_{L^2}^2+\check{C},
\end{split}
\ee
 for $\epsilon>0$ to be
decided later.

  By (\ref{rho u tt+.=}) and the elliptic estimates, together with the H\"older
inequality, the Sobolev inequality, (\ref{3d-rho}), (\ref{3d-u
theta}) and (\ref{3d-rhott}), we have \be\label{3d-H 2 of u
t}\begin{split} \|\nabla^2u_t(\cdot,t)\|_{L^2}\le&\check{C}\|\rho
u_{tt}\|_{L^2}+\check{C}\|\rho_tu_t\|_{L^2}+\check{C}\|\rho_t
u\cdot\nabla u\|_{L^2}+\check{C}\|\rho u_t\cdot\nabla
u\|_{L^2}\\&+\check{C}\|\rho u\cdot\nabla
u_t\|_{L^2}+\check{C}\|\nabla
\rho_t\theta+\rho_t\nabla\theta\|_{L^2}+\check{C}\|\nabla
\rho\theta_t+\rho\nabla\theta_t\|_{L^2}\\
\le&\check{C}\|\sqrt{\rho}
u_{tt}\|_{L^2}+\check{C}\|\rho_t\|_{L^3}\|u_t\|_{L^6}+\check{C}\|\rho_t\|_{L^3}
\|u\|_{L^\infty}\|\nabla u\|_{L^6}\\&+\check{C}\|\rho\|_{L^\infty}
\|u_t\|_{L^6}\|\nabla u\|_{L^3}+\check{C}\|\rho\|_{L^\infty}
\|u\|_{L^\infty}\|\nabla u_t\|_{L^2}+\check{C}\|\nabla
\rho_t\|_{L^2}\|\theta\|_{L^\infty}\\&+\check{C}\|\rho_t\|_{L^3}\|\nabla\theta\|_{L^6}
+\check{C}\|\nabla\rho\|_{L^3}\|\theta_t\|_{L^6}+\check{C}\|\rho\|_{L^\infty}\|\nabla\theta_t\|_{L^2}\\
\le&\check{C}\|\sqrt{\rho} u_{tt}\|_{L^2}+\check{C}\|\nabla
u_t\|_{L^2}+\check{C}\|\nabla\theta_t\|_{L^2}+\check{C}.
\end{split}\ee
Substituting (\ref{3d-H 2 of u t}) into (\ref{3d-dt nab ut+1}),
taking $\epsilon$ sufficiently small, and then multiplying the
result by $t$, we have \be\label{3d-dt nab ut+2}\begin{split}
&\frac{1}{4}\int_{\mathbb{R}^3}t\rho|u_{tt}|^2+\frac{1}{2}\frac{d}{dt}\int_{\mathbb{R}^3}t\left(\mu|\nabla
u_t|^2+(\mu+\lambda)|\mathrm{div}u_t|^2\right)\\
\le&\frac{1}{2}\int_{\mathbb{R}^3}\left(\mu|\nabla
u_t|^2+(\mu+\lambda)|\mathrm{div}u_t|^2\right)-\frac{d}{dt}\int_{\mathbb{R}^3}t\left(\frac{1}{2}\rho_t|u_t|^2+\rho_t
(u\cdot\nabla) u\cdot
u_{t}\right)\\&+\int_{\mathbb{R}^3}\left(\frac{1}{2}\rho_t|u_t|^2+\rho_t
(u\cdot\nabla) u\cdot u_{t}\right)+\frac{d}{dt}\int_{\mathbb{R}^3}t
P_t\mathrm{div} u_{t}-\int_{\mathbb{R}^3} P_t\mathrm{div}
u_{t}\\&-\frac{1}{2(2\mu+\lambda)}\frac{d}{dt}\int_{\mathbb{R}^3}t
P_{t}^2+\frac{1}{2(2\mu+\lambda)}\int_{\mathbb{R}^3}
P_{t}^2+\check{C}t\|\nabla
u_t\|_{L^2}^2+\check{C}t\|\rho_{tt}\|_{L^2}^2\\&+\check{C}t\|\nabla\theta_t\|_{L^2}^2+\check{C}.
\end{split}
\ee Integrating (\ref{3d-dt nab ut+2}) over $[0,t]$ for
$t\in[0,T_1^*)$, and using (\ref{full N-S+1})$_1$, integration by
parts, (\ref{3d-rho}), (\ref{3d-u theta}), (\ref{3d-u t theta t}),
(\ref{3d-rhott}), the Cauchy inequality and the H\"older inequality,
we have \bex\begin{split}
&\frac{1}{4}\int_0^t\int_{\mathbb{R}^3}s\rho|u_{ss}|^2+\frac{1}{2}\int_{\mathbb{R}^3}t\left(\mu|\nabla
u_t|^2+(\mu+\lambda)|\mathrm{div}u_t|^2\right)\\
\le&-\int_{\mathbb{R}^3}t\Big(\rho u\cdot \nabla u_t\cdot u_t+\rho_t
(u\cdot\nabla) u\cdot
u_{t}\Big)+\int_0^t\int_{\mathbb{R}^3}\left(\rho u\cdot\nabla
u_s\cdot u_s+\rho_s (u\cdot\nabla) u\cdot
u_{s}\right)\\&+\int_{\mathbb{R}^3}t P_t\mathrm{div}
u_{t}-\int_0^t\int_{\mathbb{R}^3} P_s\mathrm{div}
u_{s}+\frac{1}{2(2\mu+\lambda)}\int_0^t\int_{\mathbb{R}^3}
P_{s}^2+\check{C}\\ \le& \frac{\mu t}{8}\int_{\mathbb{R}^3}|\nabla
u_t|^2+\check{C}t\|\rho_t\|_{L^2}\|u\|_{L^\infty}\|\nabla
u\|_{L^3}\|u_{t}\|_{L^6}+\check{C}\\ \le& \frac{\mu
t}{4}\int_{\mathbb{R}^3}|\nabla u_t|^2+\check{C}.
\end{split}
\eex This gives (\ref{3d-na ut}).
\endpf

\begin{corollary}\label{3d-cor:6.5}
Under the conditions of Theorem \ref{3d-th1.3}, it holds that \be
\label{3d-uH3}t\|\nabla
u(\cdot,t)\|_{H^2}^2+\int_0^t\int_{\mathbb{R}^3}s|\nabla^2u_s|^2\le
\check{C}, \ee for a.e. $t\in[0,T_1^*)$.
\end{corollary}
\pf By (\ref{3d-nabla 3 u}), (\ref{3d-nabla 2 rho}), (\ref{3d-na
ut}) and (\ref{3d-u theta}), we have \be \label{3d-H 3 of u}
t\|\nabla u(\cdot,t)\|_{H^2}^2\le \check{C}, \ee for a.e
$t\in[0,T_1^*)$.

By (\ref{3d-H 2 of u t}), (\ref{3d-u t theta t}) and (\ref{3d-na
ut}), we get
$$
\int_0^t\int_{\mathbb{R}^3}s|\nabla^2u_s|^2\le \check{C},
$$  for a.e. $t\in[0,T_1^*)$.
\endpf
\begin{lemma}\label{3d-le: H 3 of rho}
Under the conditions of Theorem \ref{3d-th1.3}, it holds that
 \be\label{3d-H 3 of rho}
\|\nabla^2\rho(\cdot,t)\|_{L^q}+\int_0^t\|\nabla
\theta(\cdot,s)\|_{H^2}^2ds\le\check{C},
 \ee for any $t\in[0,T_1^*)$.
\end{lemma}
\pf  Multiplying (\ref{3d-equation of H 2 rho}) by
$q|\nabla^2\rho|^{q-2}\nabla^2\rho$, integrating by parts over
$\mathbb{R}^3$, and using the H\"older inequality, the Sobolev
inequality and (\ref{3d-rho}), we have \be\label{3d-dtrho
2q}\begin{split}
\frac{d}{dt}\int_{\mathbb{R}^3}|\nabla^2\rho|^q\le&\check{C}\|\nabla
u\|_{L^\infty}\int_{\mathbb{R}^3}|\nabla^2\rho|^{q}+\check{C}\|\nabla^2\rho\|_{L^q}^{q-1}\|\nabla^2u\|_{W^{1,q}}\|\nabla\rho\|_{L^q}
\\&
+\check{C}\|\nabla^2\rho\|_{L^q}^{q-1}\|\nabla^3u\|_{L^q}\\
\le&\check{C}\|\nabla
u\|_{H^2}\int_{\mathbb{R}^3}|\nabla^2\rho|^{q}+\check{C}\|\nabla^2\rho\|_{L^q}^{q-1}\|\nabla^2u\|_{W^{1,q}}.
\end{split}
\ee By (\ref{full N-S+1})$_2$, (\ref{full N-S+1})$_3$ and the
elliptic estimates, together with the Sobolev inequality,
(\ref{3d-rho}), (\ref{3d-u theta}), (\ref{3d-u t theta t}), the
H\"older inequality and the Gagliardo-Nirenberg inequality, we have
\be\label{3d-u w 2q}\begin{split}
\|\nabla^2u\|_{W^{1,q}}\le&\check{C}\|\rho
u_t\|_{W^{1,q}}+\check{C}\|\rho u\cdot\nabla
u\|_{W^{1,q}}+\check{C}\|\nabla(\rho \theta)\|_{W^{1,q}}\\
\le&\check{C}\|\nabla(\rho u_t)\|_{L^{2}}+\check{C}\|\nabla(\rho
u_t)\|_{L^{q}}+\check{C}\|\nabla(\rho u\cdot\nabla
u)\|_{L^{q}}+\check{C}\|\nabla^2(\rho \theta)\|_{L^{q}}+\check{C}\\
\le&\check{C}\|\nabla u_t\|_{L^{2}}+\check{C}\|\nabla
u_t\|_{L^{q}}+\check{C}\|\nabla\rho\|_{L^\infty}\| \nabla
u_t\|_{L^{2}}+\check{C}\|\nabla^3
u\|_{L^{2}}\\&+\check{C}\|\nabla\rho\|_{L^\infty}+\check{C}\| \nabla
u\|_{L^\infty}+\check{C}\|\nabla^2\rho\|_{L^{q}}+\check{C}\|\nabla^3\theta\|_{L^{2}}+\check{C}\\
\le&\check{C}\|\nabla u_t\|_{L^{2}}^{\frac{6-q}{2q}}\|\nabla^2
u_t\|_{L^{2}}^{\frac{3q-6}{2q}}+\check{C}(\|\nabla^2\rho\|_{L^q}+1)\|
\nabla u_t\|_{L^{2}}+\check{C}\|\nabla^3
u\|_{L^{2}}\\&+\check{C}\|\nabla^2\rho\|_{L^{q}}+\check{C}\|\nabla^3\theta\|_{L^{2}}+\check{C},
\end{split}
\ee and
 \be \label{3d-H 3 of theta}
\begin{split}
\|\nabla^3\theta\|_{L^2}\le&\check{C}\|\nabla\rho
\theta_t+\rho\nabla\theta_t\|_{L^2}+\check{C}\|\nabla\rho
u\cdot\nabla\theta\|_{L^2}+\check{C}\|\rho \nabla
u\cdot\nabla\theta\|_{L^2}+\check{C}\|\rho
u\cdot\nabla\nabla\theta\|_{L^2}\\&+\check{C}\|\nabla\rho\theta\mathrm{div}u\|_{L^2}+\check{C}\|\rho\nabla\theta\mathrm{div}u\|_{L^2}
+\check{C}\|\rho\theta\nabla\mathrm{div}u\|_{L^2}+\check{C}\|\nabla
u\nabla^2u\|_{L^2}\\
\le&\check{C}\|\nabla\theta_t\|_{L^2}+\check{C}\|\nabla^3u\|_{L^2}+\check{C}.
\end{split}
\ee Substituting (\ref{3d-u w 2q}) and (\ref{3d-H 3 of theta}) into
(\ref{3d-dtrho 2q}), and using the Young inequality, we
have\be\label{3d-dtrho 2q+1}\begin{split}
\frac{d}{dt}\int_{\mathbb{R}^3}|\nabla^2\rho|^q\le&\check{C}(\|\nabla^3
u\|_{L^2}+\|\nabla
u_t\|_{L^{2}}+\|\nabla\theta_t\|_{L^2}+1)\left(\int_{\mathbb{R}^3}|\nabla^2\rho|^{q}+1\right)\\&+\check{C}\|\nabla
u_t\|_{L^{2}}^{\frac{6-q}{2q}}\|\nabla^2
u_t\|_{L^{2}}^{\frac{3q-6}{2q}}\left(\int_{\mathbb{R}^3}|\nabla^2\rho|^{q}+1\right).
\end{split}
\ee It is easy to see \be\label{3d-gronwall}
\begin{split}
\|\nabla u_t\|_{L^{2}}^{\frac{6-q}{2q}}\|\nabla^2
u_t\|_{L^{2}}^{\frac{3q-6}{2q}}=&\left(t\|\nabla
u_t\|_{L^{2}}^2\right)^{\frac{6-q}{4q}}t^{-\frac{6-q}{4q}}\left(t\|\nabla^2
u_t\|_{L^{2}}^2\right)^{\frac{3q-6}{4q}}t^{-\frac{3q-6}{4q}}\\ \le&
\check{C}t^{-\frac{1}{2}}\left(t\|\nabla^2
u_t\|_{L^{2}}^2\right)^{\frac{3q-6}{4q}}\\ \le&
\check{C}t^{-\frac{2q}{q+6}}+\check{C}t\|\nabla^2 u_t\|_{L^{2}}^2\in
L^1([0,T]),
\end{split}
\ee where we have used (\ref{3d-na ut}), (\ref{3d-uH3}), the Young
inequality and $q<6$.

By (\ref{3d-u t theta t}), (\ref{3d-H 2 of rho}),
(\ref{3d-gronwall}) and the Gronwall inequality, we have
 \be \label{3d- uw 2q+1}
\begin{split}
\|\nabla^2\rho(\cdot,t)\|_{L^q}\le\check{C},
\end{split}
\ee for any $t\in[0,T_1^*)$. By (\ref{3d-H 3 of theta}), (\ref{3d-u
t theta t}) and (\ref{3d-H 2 of rho}), we get
$$
\int_0^t\|\nabla \theta(\cdot,s)\|_{H^2}^2ds\le\check{C},
$$
for any $t\in[0,T_1^*)$.
\endpf

\section{Appendix}\setcounter{equation}{0} \setcounter{theorem}{0}
\label{3d-sec 8}

{\noindent\bf Appendix A} (Zlotnik inequality)\\

{\it\noindent Let the function $y$ satisfy
$$
y^\prime(t)=g(y)+b^\prime(t)\ \mathrm{on}\ [0,T],\ y(0)=y^0,
$$ with $g\in C(\mathbb{R})$ and $y,b\in W^{1,1}(0,T)$. If
$g(\infty)=-\infty$ and
$$
b(t_2)-b(t_1)\le N_0+N_1(t_2-t_1)
$$ for all $0\le t_1<t_2\le T$ with some $N_0\ge0$ and $N_1\ge0$,
then
$$
y(t)\le\max\{y^0,\bar{\zeta}\}+N_0<\infty\ \mathrm{on}\ [0,T],
$$ where $\bar{\zeta}$ is a constant such that
$$
g(\zeta)\le-N_1,\ \mathrm{for}\ \zeta\ge\bar{\zeta}.
$$ }

{\noindent\bf Appendix B} (Local classical solution) {\it Assume
that the initial data $(\rho_0,u_0,\theta_0)$ satisfies
\be\label{local-initial data}\rho_0\geq0,\ \theta_0\ge0,\ in\
\mathbb{R}^3,\ \rho_0\in H^2\cap W^{2,q},\ u_0\in D^2\cap D_0^1,\
\theta_0\in D^2\cap D_0^1, \ee for some $q\in(3,6)$, and the
compatibility conditions \beq\label{local-compatibility}
\begin{cases}
\mu\Delta u_0+(\mu+\lambda)\nabla\mathrm{div} u_0-\nabla P(\rho_0, \theta_0)=\sqrt{\rho_0}g_1,\\
\kappa\Delta\theta_0+\frac{\mu}{2}\left|\nabla u_0+(\nabla
u_0)^\prime\right|^2+\lambda(\mathrm{div}u_0)^2=\sqrt{\rho_0}g_2,\
x\in\mathbb{R}^3,
\end{cases}
\eeq for some $g_i\in L^2$, $i=1,2$. Then there exist a positive
constant $T_0>0$ and a unique classical solution $(\rho,u,\theta)$
in $\mathbb{R}^3\times[0,T_0]$ such that
\be\label{apo-regularities}\begin{split} &\rho\in C([0,T_0]; H^2\cap
W^{2,q}),\ \rho_t\in C([0,T_0]; H^1),\ \rho\ge0,\ \theta\ge0\
\mathrm{in}\ \mathbb{R}^3\times[0,T_0],&\\& (u,\theta)\in
C([0,T_0];D^2\cap D^1_0)\cap L^2([0,T_0];D^3),\ (u_t,\theta_t)\in
L^2([0,T_0]; D_0^1),&\\&
  (\sqrt{\rho} u_t, \sqrt{\rho} \theta_t)\in
L^\infty([0,T_0];L^2),\ \sqrt{t}\sqrt{\rho} u_{tt}\in L^2([0,T_0];
L^2),\ t\sqrt{\rho} u_{tt}\in L^\infty([0,T_0];L^2),&\\&
\sqrt{t}u\in L^\infty([0,T_0];D^3),\ \sqrt{t} u_t\in
L^\infty([0,T_0];D_0^1)\cap L^2([0,T_0]; D^2),&\\& tu\in
L^\infty([0,T_0];D^{3,q}),\ tu_t\in L^\infty([0,T_0]; D^2),\
tu_{tt}\in L^2([0,T_0];D_0^1),&\\& t\theta\in
L^\infty([0,T_0];D^3)\cap L^2([0,T_0];D^4),\ t\theta_t\in
L^\infty([0,T_0]; D_0^1)\cap L^2([0,T_0];D^2),&\\&
t^\frac{3}{2}\theta\in L^\infty([0,T_0]; D^4),\
t^\frac{3}{2}\theta_t\in L^\infty([0,T_0]; D^2),\
t\sqrt{\rho}\theta_{tt}\in
L^2([0,T_0];L^2),&\\&t^\frac{3}{2}\sqrt{\rho}\theta_{tt}\in
L^\infty([0,T_0];L^2),\  t^\frac{3}{2}\theta_{tt}\in
L^2([0,T_0];D_0^1).
\end{split}
\ee }

{\noindent\bf\underline{Proof of Appendix B:}}\\

Using some arguments similar to \cite{Cho-Kim, cho-Kim: perfect
gas}, we can construct a sequence of approximate classical solutions
$(\rho^k,u^k,\theta^k)$ to (\ref{full N-S+1})-(\ref{3d-boundary})
satisfying \be\label{apo-full N-S}
\begin{cases}
\rho^k_t+\nabla \cdot (\rho^k u^{k-1})=0, \\
\rho^k  u^k_t+\rho^{k} u^{k-1}\cdot\nabla u^{k}+\nabla
P(\rho^k,\theta^k)=\mu\Delta
u^k+(\mu+\lambda)\nabla\mathrm{div}u^k,\\[1.8mm]
\rho^k \theta^k_t+\rho^k
u^{k-1}\cdot\nabla\theta^{k}+\rho^k\theta^{k}\mathrm{div}u^{k-1}\\=\displaystyle\frac{\mu}{2}\left|\nabla
u^{k-1}+(\nabla
u^{k-1})^\prime\right|^2+\lambda(\mathrm{div}u^{k-1})^2+\kappa\Delta\theta^k,
\end{cases}
\ee with initial conditions \be\label{apo-initial} (\rho^k, u^k,
\theta^k)|_{t=0}=(\rho_0^\delta, u_0, \theta_0),\ x\in \mathbb{R}^3,
\ee and \be\label{apo-boundary} (\rho^k,
u^k,\theta^k)\rightarrow(\delta, 0, 0)\ \mathrm{as}\
|x|\rightarrow\infty,\ \mathrm{for}\ t\ge0,\ee where
$\rho_0^\delta=\rho_0+\delta$ for $\delta\in(0,1)$, and $k\ge1$.
Here we take $(u^0,\theta^0)=(u_0,\theta_0)$.

From the compatibility condition (\ref{local-compatibility}), we
obtain \be\label{apo-compa}\begin{cases} \mu\Delta
u_0+(\mu+\lambda)\nabla\mathrm{div}u_0-\nabla
P(\rho_0^\delta,\theta_0)=\sqrt{\rho_0^\delta} g_1^\delta,\\
\kappa\Delta\theta_0+\frac{\mu}{2}\left|\nabla u_0+(\nabla
u_0)^\prime\right|^2+\lambda(\mathrm{div}u_0)^2=\sqrt{\rho_0^\delta}g_2^\delta,
\end{cases}\ee where
$$
g_1^\delta=\left(\frac{\rho_0}{\rho_0^\delta}\right)^\frac{1}{2}g_1-\frac{\delta}{\sqrt{\rho_0^\delta}}\nabla\theta_0,\
\mathrm{and}\
g_2^\delta=\left(\frac{\rho_0}{\rho_0^\delta}\right)^\frac{1}{2}g_2.
$$
It is easy to verify  \be\label{apo-g1 g2} \|g_1^\delta\|_{L^2}\le
\|g_1\|_{L^2}+\sqrt{\delta}\|\nabla\theta_0\|_{L^2},\ \mathrm{and}\
\|g_2^\delta\|_{L^2}\le\|g_2\|_{L^2}. \ee

{\noindent\bf Step 1: Some estimates.}\\

From \cite{cho-Kim: perfect gas} together with (\ref{apo-full
N-S})-(\ref{apo-g1 g2}), we get the following lemma.

\begin{lemma}\label{apo-le: 7.1} Under the condition of (\ref{local-initial data}), (\ref{apo-compa}) and
(\ref{apo-g1 g2}), there exists a constant $T_0\in(0,1)$ independent
of $k$ and $\delta$, such that \be\label{3d-rho k} \rho^k>0,\
\|\rho^k(\cdot,t)\|_{L^\infty}+\|\rho^k(\cdot,t)-\delta\|_{H^1\cap
W^{1,q}}+\|\rho^k_t(\cdot,t)\|_{L^2\cap L^q}\le\check{C}, \ee
\be\label{3d-u theta k}\begin{split}
 &\|u^k(\cdot,t)\|_{L^\infty}+\|\nabla
 u^k(\cdot,t)\|_{H^1}+\|\nabla\theta^k(\cdot,t)\|_{H^1}+\|\nabla^2u^k\|_{L^2([0,t];L^q)}
 \\&+\|\nabla^2\theta^k\|_{L^2([0,t];L^q)}\le\check{C},\end{split}\ee
\be\label{3d-u t theta t k}\|\nabla
u_t^k\|_{L^2([0,t];L^2)}+\|\nabla \theta_t^k\|_{L^2([0,t];L^2)}+
\|\sqrt{\rho^k} u_t^k(\cdot,t)\|_{L^2}+\|\sqrt{\rho^k}
\theta_t^k(\cdot,t)\|_{L^2}\le\check{C},\ee for any $k\ge1$ and a.e.
$t\in[0,T_0]$, where $\check{C}$ is independent of $k$, $\delta$ and
$t$. Furthermore, ${\it \theta^k\ge0}$.\footnote{This can be
obtained by using
 (\ref{apo-full N-S})$_3$ and the maximal principle for the parabolic
 equation.}
\end{lemma}
Based on {\it Lemma \ref{apo-le: 7.1}}\footnote{It is similar to
Lemma \ref{3d-le:6.1}.}, we derive the next lemma by using some
arguments similar to Lemmas \ref{3d-le: H 2 of rho}, \ref{3d-le:6.4}
and \ref{3d-le: H 3 of rho}, and Corollaries \ref{3d-cor rho tt} and
\ref{3d-cor:6.5}.
\begin{lemma}\label{apo-le:7.2}
 Under the condition of (\ref{local-initial data}), (\ref{apo-compa}) and
(\ref{apo-g1 g2}), it holds that
 \be\label{apo-total}\begin{cases}
\sqrt{t}\|\nabla u^k(\cdot,t)\|_{H^2}+\sqrt{t}\|\nabla
u_t^k(\cdot,t)\|_{L^2}+\|\rho_t^k(\cdot,t)\|_{H^1}+\|\rho^k(\cdot,t)-\delta\|_{H^2\cap W^{2,q}}\le\check{C},\\[3mm]
\displaystyle\int_0^{T_0}\left(\|\rho^k_{tt}\|_{L^2}^2+\|\nabla
u^k\|_{H^2}^2+t\|\sqrt{\rho^k}
u_{tt}^k\|_{L^2}^2+t\|\nabla^2u_t^k\|_{L^2}^2+\|\nabla
\theta^k\|_{H^2}^2
 \right)(\cdot,t)dt\le\check{C},\end{cases}\ee for any $k\ge1$ and a.e. $t\in[0,T_0]$.
\end{lemma}

We need some higher order estimates for $(\rho^k,u^k,\theta^k)$
which are included in the following lemmas.
\begin{lemma}\label{apo-le:7.3}
Under the condition of (\ref{local-initial data}), (\ref{apo-compa})
and (\ref{apo-g1 g2}), it holds that
 \be \label{apo-le:utheta total}
 \begin{split}
\int_{\mathbb{R}^3}t^2\left(|\nabla\theta_t^k|^2+\rho^{k}
|u_{tt}^{k}|^2\right)
+\int_0^{T_0}\int_{\mathbb{R}^3}t^2\left(\rho^k
|\theta^k_{tt}|^2+|\nabla u_{tt}^{k}|^2\right)\le\check{C},
 \end{split}\ee for any $k\ge1$ and a.e. $t\in[0,T_0]$.
\end{lemma}
\pf Differentiating (\ref{apo-full N-S})$_2$ w.r.t. $t$ two times,
we have \be\label{rho u tt k+.=}\begin{split} &\rho^k
u_{ttt}^k+2\rho_t^k u_{tt}^k+\rho_{tt}^ku_t^k+\rho_{tt}^k
u^{k-1}\cdot\nabla u^k+2\rho_t^k u_t^{k-1}\cdot\nabla u^k+2\rho_t^k
u^{k-1}\cdot\nabla u_t^k+\\&\rho^k u_{tt}^{k-1}\cdot\nabla
u^k+2\rho^k u_t^{k-1}\cdot\nabla u_t^k+\rho^k u^{k-1}\cdot\nabla
u_{tt}^k+\nabla P_{tt}^k=\mu\Delta
u_{tt}^k+(\mu+\lambda)\nabla\mathrm{div}u_{tt}^k.\end{split} \ee
Multiplying (\ref{rho u tt k+.=}) by $u_{tt}^{k}$, integrating by
parts over $\mathbb{R}^3$, and using (\ref{apo-full N-S})$_1$, we
have \be
 \label{apo-dt rhoutt}
 \begin{split}
&\frac{1}{2}\frac{d}{dt}\int_{\mathbb{R}^3}\rho^{k}
|u_{tt}^{k}|^2+\int_{\mathbb{R}^3}\left(\mu|\nabla
u_{tt}^{k}|^2+(\mu+\lambda)|\mathrm{div}u_{tt}^{k}|^2\right)\\
=&-2\int_{\mathbb{R}^3}\rho_t^{k}
|u_{tt}^{k}|^2-\int_{\mathbb{R}^3}\rho_{tt}^{k}u_t^{k}\cdot
u_{tt}^{k}-\int_{\mathbb{R}^3}\rho_{tt}^{k} u^{k-1}\cdot\nabla
u^{k}\cdot u_{tt}^{k}-2\int_{\mathbb{R}^3}\rho_t^{k}
u_t^{k-1}\cdot\nabla u^{k}\cdot
u_{tt}^{k}\\&-2\int_{\mathbb{R}^3}\rho_t^{k} u^{k-1}\cdot\nabla
u_t^{k}\cdot u_{tt}^{k}-\int_{\mathbb{R}^3}\rho^{k}
u_{tt}^{k-1}\cdot\nabla u^{k}\cdot
u_{tt}^{k}-2\int_{\mathbb{R}^3}\rho^{k} u_t^{k-1}\cdot\nabla
u_t^{k}\cdot u_{tt}^{k}\\&+\int_{\mathbb{R}^3}
P_{tt}^{k}\mathrm{div}u_{tt}^{k}=\sum\limits_{i=1}^8VII_i.
 \end{split}
\ee For $VII_1$, using (\ref{apo-full N-S})$_1$ and integration by
parts again, together with the Cauchy inequality, (\ref{3d-rho k})
and (\ref{3d-u theta k}), we have \be\label{apo-VII 1}
\begin{split}
VII_1=&-4\int_{\mathbb{R}^3}\rho^{k}u^{k-1}\cdot \nabla
u_{tt}^{k}\cdot u_{tt}^k\le\frac{\mu}{16}\int_{\mathbb{R}^3}|\nabla
u_{tt}^{k}|^2+\check{C}\int_{\mathbb{R}^3}\rho^{k}|u_{tt}^{k}|^2.
\end{split}
\ee For $VII_2$, we have \be\label{apo-VII 2}
\begin{split}VII_2=&-\int_{\mathbb{R}^3}(\rho^{k}u^{k-1})_t\cdot\nabla(u_t^{k}\cdot
u_{tt}^{k})\\
\le&\|\rho^{k}_t\|_{L^3}\|u^{k-1}\|_{L^\infty}\|\nabla
u_t^{k}\|_{L^2}\|u_{tt}^{k}\|_{L^6}+\|\rho^{k}_t\|_{L^3}\|u^{k-1}\|_{L^\infty}\|u_t^{k}\|_{L^6}\|\nabla
u_{tt}^{k}\|_{L^2}\\&+\|\sqrt{\rho^{k}}\|_{L^\infty}\|u^{k-1}_t\|_{L^6}\|\nabla
u_t^{k}\|_{L^3}\|\sqrt{\rho^k}
u_{tt}^{k}\|_{L^2}+\|\rho^{k}u_t^{k}\|_{L^3}\|u^{k-1}_t\|_{L^6}\|\nabla
u_{tt}^{k}\|_{L^2}\\
\le&\frac{\mu}{16}\int_{\mathbb{R}^3}|\nabla
u_{tt}^{k}|^2+\check{C}\|\nabla u_t^{k}\|_{L^2}^2+\check{C}\|\nabla
u^{k-1}_t\|_{L^2}^2\|\nabla u_t^{k}\|_{L^2}\|\nabla
u_t^{k}\|_{L^6}\\&+\check{C}\|\sqrt{\rho^k}
u_{tt}^{k}\|_{L^2}^2+\check{C}\|\rho^{k}u_t^{k}\|_{L^2}\|\rho^{k}u_t^{k}\|_{L^6}\|\nabla
u^{k-1}_t\|_{L^2}^2\\
\le&\frac{\mu}{16}\int_{\mathbb{R}^3}|\nabla
u_{tt}^{k}|^2+\check{C}\|\nabla u_t^{k}\|_{L^2}^4+\check{C}\|\nabla
u^{k-1}_t\|_{L^2}^4+\check{C}\|\nabla u_t^{k}\|_{L^2}^2\|\nabla^2
u_t^{k}\|_{L^2}^2\\&+\check{C}\|\sqrt{\rho^k}u_{tt}^{k}\|_{L^2}^2+\check{C},
\end{split}
\ee where we have used (\ref{apo-full N-S})$_1$, integration by
parts, the H\"older inequality, the Cauchy inequality, the Sobolev
inequality, the Gagliardo-Nirenberg inequality, (\ref{3d-rho k}),
(\ref{3d-u theta k}) and (\ref{3d-u t theta t k}).

Similarly, we have \be\label{apo-VII 3-5}
\begin{split}
&VII_3+VII_4+VII_5\\
\le&\check{C}\|\rho_{tt}^{k}\|_{L^2}\|u^{k-1}\|_{L^\infty}\|\nabla
u^{k}\|_{L^3}\|u_{tt}^{k}\|_{L^6}+\check{C}\|\rho_t^{k}\|_{L^2}
\|u_t^{k-1}\|_{L^6}\|\nabla
u^{k}\|_{L^6}\|u_{tt}^{k}\|_{L^6}\\&+\check{C}\|\rho_t^{k}\|_{L^2}
\|u^{k-1}\|_{L^6}\|\nabla u_t^{k}\|_{L^6}\|u_{tt}^{k}\|_{L^6}
\\
\le&\frac{3\mu}{16}\int_{\mathbb{R}^3}|\nabla
u_{tt}^{k}|^2+\check{C}\|\rho_{tt}^{k}\|_{L^2}^2+\check{C} \|\nabla
u_t^{k-1}\|_{L^2}^2+\check{C} \|\nabla^2 u_t^{k}\|_{L^2}^2,
\end{split}
\ee  and\be\label{apo-VII 6-7}
\begin{split}
VII_6+VII_7 \le&\check{C}
\|\rho^{k}u_{tt}^{k}\|_{L^2}\|u_{tt}^{k-1}\|_{L^6}\|\nabla
u^{k}\|_{L^3}+\check{C}\|\rho^{k}u_{tt}^{k}\|_{L^2}
\|u_t^{k-1}\|_{L^6}\|\nabla u_t^{k}\|_{L^3}\\
\le&\frac{\mu}{8}\int_{\mathbb{R}^3}|\nabla
u_{tt}^{k-1}|^2+\check{C}
\|\sqrt{\rho^{k}}u_{tt}^{k}\|_{L^2}^2+\check{C} \|\nabla
u_t^{k-1}\|_{L^2}^2\|\nabla u_t^{k}\|_{L^2}\|\nabla^2
u_t^{k}\|_{L^2}\\
\le&\frac{\mu}{8}\int_{\mathbb{R}^3}|\nabla
u_{tt}^{k-1}|^2+\check{C}
\|\sqrt{\rho^{k}}u_{tt}^{k}\|_{L^2}^2+\check{C} \|\nabla
u_t^{k-1}\|_{L^2}^4+\check{C}\|\nabla u_t^{k}\|_{L^2}^2\|\nabla^2
u_t^{k}\|_{L^2}^2,
\end{split}
\ee and \be\label{apo-VII 8}\begin{split} VII_8\le&
\frac{\mu}{16}\int_{\mathbb{R}^3}|\nabla
u_{tt}^{k}|^2+\check{C}\int_{\mathbb{R}^3}|\rho^k_{tt}\theta^k+2\rho^k_t\theta^k_t+\rho^k\theta_{tt}^k|^2\\
\le&\frac{\mu}{16}\int_{\mathbb{R}^3}|\nabla
u_{tt}^{k}|^2+\check{C}\int_{\mathbb{R}^3}|\rho^k_{tt}|^2+\check{C}\|\rho^k_t\|_{L^3}^2\|\theta^k_t\|_{L^6}^2+\check{C}\int_{\mathbb{R}^3}\rho^k|\theta_{tt}^k|^2
\\
\le&\frac{\mu}{16}\int_{\mathbb{R}^3}|\nabla
u_{tt}^{k}|^2+\check{C}\int_{\mathbb{R}^3}|\rho^k_{tt}|^2+\check{C}\int_{\mathbb{R}^3}|\nabla\theta^k_t|^2+\check{C}\int_{\mathbb{R}^3}\rho^k|\theta_{tt}^k|^2.
\end{split}
\ee Substituting (\ref{apo-VII 1}), (\ref{apo-VII 2}), (\ref{apo-VII
3-5}), (\ref{apo-VII 6-7}) and (\ref{apo-VII 8}) into (\ref{apo-dt
rhoutt}), and multiplying the result by $t^2$, we have\be
 \label{apo-dt rhoutt+1}
 \begin{split}
&\frac{1}{2}\frac{d}{dt}\int_{\mathbb{R}^3}t^2\rho^{k}
|u_{tt}^{k}|^2+\frac{5}{8}\int_{\mathbb{R}^3}t^2\left(\mu|\nabla
u_{tt}^{k}|^2+(\mu+\lambda)|\mathrm{div}u_{tt}^{k}|^2\right)\\
\le&\int_{\mathbb{R}^3}t\rho^{k}
|u_{tt}^{k}|^2+\frac{\mu}{8}\int_{\mathbb{R}^3}t^2|\nabla
u_{tt}^{k-1}|^2+\check{C}\int_{\mathbb{R}^3}t^2\rho^{k}|u_{tt}^{k}|^2+\check{C}t^2\|\nabla
u_t^{k}\|_{L^2}^4\\&+\check{C}t^2\|\nabla
u^{k-1}_t\|_{L^2}^4+\check{C}t^2\|\rho_{tt}^{k}\|_{L^2}^2+\check{C}
t^2\|\nabla^2 u_t^{k}\|_{L^2}^2+\check{C}t^2\|\nabla
u_t^{k}\|_{L^2}^2\|\nabla^2
u_t^{k}\|_{L^2}^2\\&+\check{C}\int_{\mathbb{R}^3}t^2|\nabla\theta^k_t|^2+\check{C}\int_{\mathbb{R}^3}t^2\rho^k|\theta_{tt}^k|^2+\check{C}.
 \end{split}
\ee Integrating (\ref{apo-dt rhoutt+1}) over $[0,t]$ for
$t\in[0,T_0]$, and using (\ref{3d-u t theta t k}) and
(\ref{apo-total}), for any given $N\in \mathbb{Z}_+$, we have \be
\label{apo-int rhoutt+1}
 \begin{split}
&\max\limits_{1\le k\le N}\int_{\mathbb{R}^3}t^2\rho^{k}
|u_{tt}^{k}|^2+\max\limits_{1\le k\le
N}\mu\int_0^t\int_{\mathbb{R}^3}s^2|\nabla u_{ss}^{k}|^2
\le\check{C}\max\limits_{1\le k\le
N}\int_0^t\int_{\mathbb{R}^3}s^2\rho^k|\theta_{ss}^k|^2+\check{C}.
 \end{split}
\ee Differentiating (\ref{apo-full N-S})$_3$ w.r.t. $t$, we have
\be\label{apo-rho theta tt=}
\begin{split}
&\rho^k \theta^k_{tt}+\rho^k_t \theta^k_t+\rho^k_t
u^{k-1}\cdot\nabla\theta^{k}+\rho^k
u^{k-1}_t\cdot\nabla\theta^{k}+\rho^k
u^{k-1}\cdot\nabla\theta^{k}_t+\rho^k_t\theta^{k}\mathrm{div}u^{k-1}\\&+\rho^k\theta^{k}_t\mathrm{div}u^{k-1}+\rho^k\theta^{k}\mathrm{div}u^{k-1}_t
=\mu\left(\nabla u^{k-1}+(\nabla u^{k-1})^\prime\right):\left(\nabla
u^{k-1}_t+(\nabla
u^{k-1}_t)^\prime\right)\\&+2\lambda\mathrm{div}u^{k-1}\mathrm{div}u^{k-1}_t+\kappa\Delta\theta^k_t.
\end{split}
\ee Multiplying (\ref{apo-rho theta tt=}) by $\theta_{tt}^k$, and
integrating by parts over $\mathbb{R}^3$, we have \be
\label{apo-dtnabla thetat}
\begin{split}
&\int_{\mathbb{R}^3}\rho^k
|\theta^k_{tt}|^2+\frac{\kappa}{2}\frac{d}{dt}\int_{\mathbb{R}^3}|\nabla\theta_t^k|^2\\
=&-\int_{\mathbb{R}^3}\rho^k_t
\theta^k_t\theta^k_{tt}-\int_{\mathbb{R}^3}\rho^k_t
u^{k-1}\cdot\nabla\theta^{k}\theta^k_{tt}-\int_{\mathbb{R}^3}\rho^k
u^{k-1}_t\cdot\nabla\theta^{k}\theta_{tt}^k-\int_{\mathbb{R}^3}\rho^k
u^{k-1}\cdot\nabla\theta^{k}_t\theta_{tt}^k\\&-\int_{\mathbb{R}^3}\rho^k_t\theta^{k}\mathrm{div}u^{k-1}\theta_{tt}^k
-\int_{\mathbb{R}^3}\rho^k\theta^{k}_t\mathrm{div}u^{k-1}\theta_{tt}^k-\int_{\mathbb{R}^3}\rho^k\theta^{k}\mathrm{div}u^{k-1}_t\theta_{tt}^k
\\&+\mu\int_{\mathbb{R}^3}\left(\nabla u^{k-1}+(\nabla
u^{k-1})^\prime\right):\left(\nabla u^{k-1}_t+(\nabla
u^{k-1}_t)^\prime\right)\theta_{tt}^k+2\lambda\int_{\mathbb{R}^3}\mathrm{div}u^{k-1}\mathrm{div}u^{k-1}_t\theta_{tt}^k\\
=&\sum\limits_{i=1}^9VIII_{i}.
\end{split}
\ee  For $VIII_1$ and $VIII_2$, similar to (\ref{3d-VI 1}), we have
\be \label{apo-VIII1}
\begin{split}
VIII_1 =&-\frac{1}{2}\frac{d}{dt}\int_{\mathbb{R}^3}\rho^k_t
|\theta^k_t|^2+\int_{\mathbb{R}^3}(\rho^ku^{k-1})_{t}\cdot
\nabla\theta^k_t\theta_t^k\\
\le&-\frac{1}{2}\frac{d}{dt}\int_{\mathbb{R}^3}\rho^k_t
|\theta^k_t|^2+\|\rho^k_{t}\|_{L^3}\|u^{k-1}\|_{L^\infty}\|\nabla\theta^k_t\|_{L^2}
\|\theta_t^k\|_{L^6}\\&+\|\sqrt{\rho^k}\|_{L^\infty}\|u^{k-1}_{t}\|_{L^\infty}\|\nabla\theta^k_t\|_{L^2}
\|\sqrt{\rho^k}\theta_t^k\|_{L^2}\\
\le&-\frac{1}{2}\frac{d}{dt}\int_{\mathbb{R}^3}\rho^k_t
|\theta^k_t|^2+\check{C}\|\nabla\theta^k_t\|_{L^2}^2+\check{C}\|\nabla
u^{k-1}_t\|_{L^2}^2+\check{C}\|\nabla^2 u^{k-1}_t\|_{L^2}^2,
\end{split}
\ee and \be \label{apo-VIII2}
\begin{split}
VIII_2=&-\frac{d}{dt}\int_{\mathbb{R}^3}\rho^k_t
u^{k-1}\cdot\nabla\theta^{k}\theta^k_{t}+\int_{\mathbb{R}^3}\rho^k_{tt}
u^{k-1}\cdot\nabla\theta^{k}\theta^k_{t}+\int_{\mathbb{R}^3}\rho^k_t
u^{k-1}_t\cdot\nabla\theta^{k}\theta^k_{t}\\&+\int_{\mathbb{R}^3}\rho^k_t
u^{k-1}\cdot\nabla\theta^{k}_t\theta^k_{t}\\
\le&-\frac{d}{dt}\int_{\mathbb{R}^3}\rho^k_t
u^{k-1}\cdot\nabla\theta^{k}\theta^k_{t}+\|\rho^k_{tt}\|_{L^2}
\|u^{k-1}\|_{L^\infty}\|\nabla\theta^{k}\|_{L^3}\|\theta^k_{t}\|_{L^6}\\&+\|\rho^k_t\|_{L^2}
\|u^{k-1}_t\|_{L^6}\|\nabla\theta^{k}\|_{L^6}\|\theta^k_{t}\|_{L^6}+\|\rho^k_t\|_{L^6}
\|u^{k-1}\|_{L^6}\|\nabla\theta^{k}_t\|_{L^2}\|\theta^k_{t}\|_{L^6}\\
\le&-\frac{d}{dt}\int_{\mathbb{R}^3}\rho^k_t
u^{k-1}\cdot\nabla\theta^{k}\theta^k_{t}+\check{C}\|\rho^k_{tt}\|_{L^2}^2+\check{C}\|\nabla\theta^{k}_t\|_{L^2}^2+\check{C}
\|\nabla u^{k-1}_t\|_{L^2}^2.
\end{split}
\ee Similarly, for the rest terms of the right side of
(\ref{apo-dtnabla thetat}), we have \be\label{apo-VIII3467}
\begin{split}&VIII_3+VIII_4+VIII_6+VIII_7\\
\le&\|\sqrt{\rho^k}\|_{L^\infty}\|\sqrt{\rho^k}\theta_{tt}^k\|_{L^2}
\left(\|u^{k-1}_t\|_{L^6}\|\nabla\theta^{k}\|_{L^3}+\|u^{k-1}\|_{L^\infty}\|\nabla\theta^{k}_t\|_{L^2}\right)
\\&
+\|\sqrt{\rho^k}\|_{L^\infty}\|\sqrt{\rho^k}\theta_{tt}^k\|_{L^2}\left(\|\theta^{k}_t\|_{L^6}\|\mathrm{div}u^{k-1}\|_{L^3}
+\|\theta^{k}\|_{L^\infty}\|\mathrm{div}u^{k-1}_t\|_{L^2}\right)\\
\le&\frac{1}{2}\|\sqrt{\rho^k}\theta_{tt}^k\|_{L^2}^2+
\check{C}\|\nabla
u^{k-1}_t\|_{L^2}^2+\check{C}\|\nabla\theta^{k}_t\|_{L^2}^2,
\end{split}
\ee and \be \label{apo-VIII5}
\begin{split}
VIII_5=&-\frac{d}{dt}\int_{\mathbb{R}^3}\rho^k_t\theta^{k}\mathrm{div}u^{k-1}\theta_{t}^k+
\int_{\mathbb{R}^3}\rho^k_{tt}\theta^{k}\mathrm{div}u^{k-1}\theta_{t}^k+\int_{\mathbb{R}^3}\rho^k_t|\theta^{k}_t|^2\mathrm{div}u^{k-1}\\&+
\int_{\mathbb{R}^3}\rho^k_t\theta^{k}\mathrm{div}u^{k-1}_t\theta_{t}^k\\
\le&-\frac{d}{dt}\int_{\mathbb{R}^3}\rho^k_t\theta^{k}\mathrm{div}u^{k-1}\theta_{t}^k+
\|\rho^k_{tt}\|_{L^2}\|\theta^{k}\|_{L^\infty}\|\mathrm{div}u^{k-1}\|_{L^3}\|\theta_{t}^k\|_{L^6}\\&+\|\rho^k_t\|_{L^2}
\|\theta^{k}_t\|_{L^6}^2\|\mathrm{div}u^{k-1}\|_{L^6}+
\|\rho^k_t\|_{L^3}\|\theta^{k}\|_{L^\infty}\|\mathrm{div}u^{k-1}_t\|_{L^2}\|\theta_{t}^k\|_{L^6}\\
\le&-\frac{d}{dt}\int_{\mathbb{R}^3}\rho^k_t\theta^{k}\mathrm{div}u^{k-1}\theta_{t}^k+
\check{C}\|\rho^k_{tt}\|_{L^2}^2+\check{C}
\|\nabla\theta^{k}_t\|_{L^2}^2+\|\nabla u^{k-1}_t\|_{L^2}^2,
\end{split}
\ee and \be \label{apo-VIII8}
\begin{split}
VIII_8 \le&\mu\frac{d}{dt}\int_{\mathbb{R}^3}\left(\nabla
u^{k-1}+(\nabla u^{k-1})^\prime\right):\left(\nabla
u^{k-1}_t+(\nabla
u^{k-1}_t)^\prime\right)\theta_{t}^k\\&+\check{C}\|\nabla
u^{k-1}\|_{L^3}\|\nabla
u^{k-1}_{tt}\|_{L^2}\|\theta_{t}^k\|_{L^6}+\check{C}\|\nabla
u^{k-1}_{t}\|_{L^2}\|\nabla
u^{k-1}_{t}\|_{L^3}\|\theta_{t}^k\|_{L^6}\\
\le&\mu\frac{d}{dt}\int_{\mathbb{R}^3}\left(\nabla u^{k-1}+(\nabla
u^{k-1})^\prime\right):\left(\nabla u^{k-1}_t+(\nabla
u^{k-1}_t)^\prime\right)\theta_{t}^k\\&+\check{C}\|\nabla
u^{k-1}_{tt}\|_{L^2}\|\nabla\theta_{t}^k\|_{L^2}+\check{C}\|\nabla
u^{k-1}_{t}\|_{L^2}\left(\|\nabla u^{k-1}_{t}\|_{L^2}+\|\nabla^2
u^{k-1}_{t}\|_{L^2}\right)\|\nabla\theta_{t}^k\|_{L^2},
\end{split}
\ee and \be\label{apo-VIII9}
\begin{split}
VIII_9
\le&2\lambda\frac{d}{dt}\int_{\mathbb{R}^3}\mathrm{div}u^{k-1}\mathrm{div}u^{k-1}_t\theta_{t}^k+\check{C}\|\nabla
u^{k-1}_{tt}\|_{L^2}\|\nabla\theta_{t}^k\|_{L^2}\\&+\check{C}\|\nabla
u^{k-1}_{t}\|_{L^2}\left(\|\nabla u^{k-1}_{t}\|_{L^2}+\|\nabla^2
u^{k-1}_{t}\|_{L^2}\right)\|\nabla\theta_{t}^k\|_{L^2}.
\end{split}
\ee Substituting (\ref{apo-VIII1}), (\ref{apo-VIII2}),
(\ref{apo-VIII3467}), (\ref{apo-VIII5}), (\ref{apo-VIII8}) and
(\ref{apo-VIII9}) into (\ref{apo-dtnabla thetat}), and multiplying
the result by $t^2$, we have \be \label{apo-dtnabla thetat+1}
\begin{split}
&\frac{1}{2}\int_{\mathbb{R}^3}t^2\rho^k
|\theta^k_{tt}|^2+\frac{\kappa}{2}\frac{d}{dt}\int_{\mathbb{R}^3}t^2|\nabla\theta_t^k|^2\\
\le&\kappa\int_{\mathbb{R}^3}t|\nabla\theta_t^k|^2-\frac{1}{2}\frac{d}{dt}\int_{\mathbb{R}^3}t^2\rho^k_t
|\theta^k_t|^2+\int_{\mathbb{R}^3}t\rho^k_t
|\theta^k_t|^2-\frac{d}{dt}\int_{\mathbb{R}^3}t^2\rho^k_t
u^{k-1}\cdot\nabla\theta^{k}\theta^k_{t}\\&+2\int_{\mathbb{R}^3}t\rho^k_t
u^{k-1}\cdot\nabla\theta^{k}\theta^k_{t}-\frac{d}{dt}\int_{\mathbb{R}^3}t^2\rho^k_t\theta^{k}\mathrm{div}u^{k-1}\theta_{t}^k+
2\int_{\mathbb{R}^3}t\rho^k_t\theta^{k}\mathrm{div}u^{k-1}\theta_{t}^k\\&+
\mu\frac{d}{dt}\int_{\mathbb{R}^3}t^2\left(\nabla u^{k-1}+(\nabla
u^{k-1})^\prime\right):\left(\nabla u^{k-1}_t+(\nabla
u^{k-1}_t)^\prime\right)\theta_{t}^k\\&-2\mu\int_{\mathbb{R}^3}t\left(\nabla
u^{k-1}+(\nabla u^{k-1})^\prime\right):\left(\nabla
u^{k-1}_t+(\nabla
u^{k-1}_t)^\prime\right)\theta_{t}^k\\
&+2\lambda\frac{d}{dt}\int_{\mathbb{R}^3}t^2\mathrm{div}u^{k-1}\mathrm{div}u^{k-1}_t\theta_{t}^k-
4\lambda\int_{\mathbb{R}^3}t\mathrm{div}u^{k-1}\mathrm{div}u^{k-1}_t\theta_{t}^k\\&
+\check{C}\|\rho^k_{tt}\|_{L^2}^2+\check{C}t^2\|\nabla
u^{k-1}_{tt}\|_{L^2}\|\nabla\theta_{t}^k\|_{L^2}+\check{C}t\|\nabla^2
u^{k-1}_{t}\|_{L^2}^2+\check{C},
\end{split}
\ee where we have used (\ref{apo-total}), (\ref{apo-le:utheta
total}) and the Cauchy inequality.

Integrating (\ref{apo-dtnabla thetat+1}) over $[0,t]$ for
$t\in[0,T_0]$, and using (\ref{apo-full N-S})$_1$, integration by
parts, the H\"older inequality, (\ref{3d-rho k}), (\ref{3d-u theta
k}), (\ref{3d-u t theta t k}), (\ref{apo-total}) and the Cauchy
inequality, we have \bex
\begin{split}
&\frac{1}{2}\int_0^t\int_{\mathbb{R}^3}s^2\rho^k
|\theta^k_{ss}|^2+\frac{\kappa}{2}\int_{\mathbb{R}^3}t^2|\nabla\theta_t^k|^2\\
\le&t^2\|\sqrt{\rho^k}\|_{L^\infty}\|\sqrt{\rho^k}\theta^k_t\|_{L^2}\|u^{k-1}\|_{L^\infty}\|\nabla\theta^k_t\|_{L^2}
+t^2\|\rho^k_t\|_{L^2}
\|u^{k-1}\|_{L^\infty}\|\nabla\theta^{k}\|_{L^3}\|\theta^k_{t}\|_{L^6}\\&+t^2\|\rho^k_t\|_{L^2}\|\theta^{k}\|_{L^\infty}
\|\mathrm{div}u^{k-1}\|_{L^3}\|\theta_{t}^k\|_{L^6}+4\mu t^2\|\nabla
u^{k-1}\|_{L^3}\|\nabla
u^{k-1}_t\|_{L^2}\|\theta_{t}^k\|_{L^6}\\&+2\lambda
t^2\|\mathrm{div}u^{k-1}\|_{L^3}\|\mathrm{div}u^{k-1}_t\|_{L^2}\|\theta_{t}^k\|_{L^6}
+\check{C}\int_0^t s^2\|\nabla
u^{k-1}_{ss}\|_{L^2}\|\nabla\theta_{s}^k\|_{L^2}+\check{C}\\
\le&\frac{\kappa}{4}\int_{\mathbb{R}^3}t^2|\nabla\theta_t^k|^2
+\frac{\epsilon}{4}\int_0^t s^2\|\nabla
u^{k-1}_{ss}\|_{L^2}^2+\check{C}_\epsilon,
\end{split}
\eex for $\epsilon>0$ to be decided later. This gives \be
\label{apo-intnabla thetat}
\begin{split}
\max\limits_{1\le k\le N}\int_0^t\int_{\mathbb{R}^3}s^2\rho^k
|\theta^k_{ss}|^2+\max\limits_{1\le k\le
N}\kappa\int_{\mathbb{R}^3}t^2|\nabla\theta_t^k|^2
\le\epsilon\max\limits_{1\le k\le N}\int_0^t s^2\|\nabla
u^{k}_{ss}\|_{L^2}^2+\check{C}_\epsilon,
\end{split}
\ee for any given $N\in\mathbb{Z}_+$.

Multiplying (\ref{apo-intnabla thetat}) by $2\check{C}$, putting the
result into (\ref{apo-int rhoutt+1}), and taking $\epsilon>0$
sufficiently small, we get

 \be \label{apo-utheta total}
 \begin{split}
\int_0^t\int_{\mathbb{R}^3}s^2\left(\rho^k |\theta^k_{ss}|^2+|\nabla
u_{ss}^{k}|^2\right)+\int_{\mathbb{R}^3}t^2\left(|\nabla\theta_t^k|^2+\rho^{k}
|u_{tt}^{k}|^2\right) \le\check{C},
 \end{split}\ee for any $k\in[1,N]$ and a.e. $t\in[0,T_0]$. Since $N$ is
 arbitrary, (\ref{apo-utheta total}) implies (\ref{apo-le:utheta
 total}).
\endpf
\begin{corollary}\label{apo-cor7.4}
Under the condition of (\ref{local-initial data}), (\ref{apo-compa})
and (\ref{apo-g1 g2}), it holds that
\be\label{apo-cor7.4k}\begin{cases} \displaystyle
t\|\nabla^2u_t^k(\cdot,t)\|_{L^2}+t\|\nabla^3\theta^k(\cdot,t)\|_{L^2}+t\|\nabla^3u^k(\cdot,t)\|_{L^q}\le\check{C},\\[2mm]
 \displaystyle\int_0^{T_0}t^2\left(
\|\nabla^2\theta^k_t\|_{L^2}^2+\|\nabla^4\theta^k\|_{L^2}^2\right)\le\check{C},
\end{cases}\ee for any $k\ge1$ and a.e. $t\in[0,T_0]$.
\end{corollary}
\pf Similar to (\ref{3d-H 2 of u t}), (\ref{3d-u w 2q}), (\ref{3d-H
3 of theta}) and (\ref{3d-gronwall}), together with
(\ref{apo-total}) and (\ref{apo-le:utheta total}), we have
\be\label{3d-H 2 of u t k}\begin{split}
t\|\nabla^2u_t^k(\cdot,t)\|_{L^2} \le&\check{C}t\|\sqrt{\rho^k}
u_{tt}^k\|_{L^2}+\check{C}t\|\nabla u_t^k\|_{L^2}+\check{C}t\|\nabla
u_t^{k-1}\|_{L^2}+\check{C}t\|\nabla\theta_t^k\|_{L^2}+\check{C}\\
\le&\check{C},
\end{split}\ee
 and\be \label{3d-H 3 of theta k}
\begin{split}
t\|\nabla^3\theta^k(\cdot,t)\|_{L^2}
\le\check{C}t\|\nabla\theta_t^k\|_{L^2}+\check{C}t\|\nabla^3u^{k-1}\|_{L^2}+\check{C}\le
\check{C},
\end{split}
\ee and\be\label{3d-u w 2q k}\begin{split}
t\|\nabla^3u^k(\cdot,t)\|_{L^{q}}\le&\check{C}t^{1-\frac{2q}{q+6}}+\check{C}t^2\|\nabla^2
u_t^k\|_{L^{2}}^2+\check{C}t\| \nabla
u_t^k\|_{L^{2}}+\check{C}t\|\nabla^3
u^k\|_{L^{2}}\\&+\check{C}t\|\nabla^3\theta^k\|_{L^{2}}+\check{C}
\le\check{C}.
\end{split}
\ee By (\ref{apo-rho theta tt=}) and the $H^2$-estimates for the
elliptic equation, together with the H\"older inequality, the
Sobolev inequality, (\ref{3d-rho k}), (\ref{3d-u theta k}) and
(\ref{3d-u t theta t k}), we have
 \be\label{apo-nabla 2theta t}
\begin{split}
\|\nabla^2\theta^k_t(\cdot,t)\|_{L^2}\le&\check{C}\|\rho^k
\theta^k_{tt}\|_{L^2}+\check{C}\|\rho^k_t\|_{L^3}
\|\theta^k_t\|_{L^6}+\check{C}\|\rho^k_t\|_{L^3}
\|u^{k-1}\|_{L^\infty}\|\nabla\theta^{k}\|_{L^6}\\&+\check{C}\|\rho^k\|_{L^\infty}
\|u^{k-1}_t\|_{L^6}\|\nabla\theta^{k}\|_{L^3}+\check{C}\|\rho^k\|_{L^\infty}
\|u^{k-1}\|_{L^\infty}\|\nabla\theta^{k}_t\|_{L^2}\\&+\check{C}\|\rho^k_t\|_{L^3}\|\theta^{k}\|_{L^\infty}\|\nabla
u^{k-1}\|_{L^6}+\check{C}\|\rho^k\|_{L^\infty}\|\theta^{k}_t\|_{L^6}\|\nabla
u^{k-1}\|_{L^3}\\&+\check{C}
\|\rho^k\|_{L^\infty}\|\theta^{k}\|_{L^\infty}\|\nabla
u^{k-1}_t\|_{L^2}+\check{C}\|\nabla u^{k-1}\|_{L^3}\|\nabla
u^{k-1}_t\|_{L^6}\\
\le&\check{C}\|\sqrt{\rho^k} \theta^k_{tt}\|_{L^2}+\check{C}
\|\nabla\theta^k_t\|_{L^2}+\check{C} \|\nabla
u^{k-1}_t\|_{L^2}+\check{C}\|\nabla^2 u^{k-1}_t\|_{L^2}+\check{C},
\end{split}
\ee which together with (\ref{apo-total}) and (\ref{apo-le:utheta
total}) deduces \be\label{apo-na thetatk H2}
\int_0^{T_0}t^2\|\nabla^2\theta^k_t\|_{L^2}^2\le\check{C}. \ee Using
(\ref{apo-full N-S})$_3$ and the elliptic estimates, together with
the H\"older inequality, the Sobolev inequality, (\ref{3d-rho k}),
(\ref{3d-u theta k}), (\ref{3d-u t theta t k}) and
(\ref{apo-total}), we have \be\label{apo-H 4 theta}
\begin{split}\|\nabla^4\theta^k\|_{L^2}\le&\check{C}\|\nabla^2\rho^k
\theta^k_t\|_{L^2}+\check{C}\|\nabla\rho^k
\nabla\theta^k_t\|_{L^2}+\check{C}\|\rho^k\nabla^2\theta^k_t\|_{L^2}+\check{C}\|\nabla\nabla\rho^k
u^{k-1}\cdot\nabla\theta^{k}\|_{L^2}\\&+\check{C}\|\nabla\rho^k
\nabla u^{k-1}\cdot\nabla\theta^{k}\|_{L^2}+\check{C}\|\nabla\rho^k
u^{k-1}\cdot\nabla\nabla\theta^{k}\|_{L^2}+\check{C}\|\rho^k
\nabla\nabla
u^{k-1}\cdot\nabla\theta^{k}\|_{L^2}\\&+\check{C}\|\rho^k \nabla
u^{k-1}\cdot\nabla\nabla\theta^{k}\|_{L^2}+\check{C}\|\rho^k
u^{k-1}\cdot\nabla\nabla^2\theta^{k}\|_{L^2}+\check{C}\|\nabla^2\rho^k\theta^{k}\mathrm{div}u^{k-1}\|_{L^2}
\\&+\check{C}\|\nabla\rho^k\nabla\theta^{k}\mathrm{div}u^{k-1}\|_{L^2}+\check{C}\|\nabla\rho^k\theta^{k}\nabla\mathrm{div}u^{k-1}\|_{L^2}
+\check{C}\|\rho^k\nabla^2\theta^{k}\mathrm{div}u^{k-1}\|_{L^2}\\&+\check{C}\|\rho^k\nabla\theta^{k}\nabla\mathrm{div}u^{k-1}\|_{L^2}
+\check{C}\|\rho^k\theta^{k}\nabla^2\mathrm{div}u^{k-1}\|_{L^2}
+\check{C}\|\nabla^2 u^{k-1}\nabla^2
u^{k-1}\|_{L^2}\\&+\check{C}\|\nabla u^{k-1}\nabla^3
u^{k-1}\|_{L^2}\\
\le&\check{C}
\|\nabla\theta^k_t\|_{L^2}+\check{C}\|\nabla^2\theta^k_t\|_{L^2}
+\check{C} \|\nabla^3\theta^{k}\|_{L^2} +\check{C}\|\nabla
u^{k-1}\|_{H^2}^2+\check{C},
\end{split}
\ee which together with (\ref{apo-total}), (\ref{apo-le:utheta
total}), (\ref{3d-H 3 of theta k}), (\ref{apo-na thetatk H2})
deduces
$$
\int_0^{T_0}t^2\|\nabla^4\theta^k\|_{L^2}^2\le\check{C}.
$$
The proof of Corollary \ref{apo-cor7.4} is complete.
\endpf
\begin{lemma}\label{apo-le:rho theta tt k}
Under the condition of (\ref{local-initial data}), (\ref{apo-compa})
and (\ref{apo-g1 g2}), it holds that \be\label{apo-rho theta ttk}
\int_{\mathbb{R}^3}t^3\rho^k
|\theta^k_{tt}|^2+\int_0^{T_0}\int_{\mathbb{R}^3}t^3|\nabla\theta^k_{tt}|^2
\le\check{C}, \ee  for any $k\ge1$ and a.e. $t\in[0,T_0]$.
\end{lemma}
\pf Differentiating (\ref{apo-rho theta tt=}) w.r.t. $t$,
multiplying the result by $\theta_{tt}^k$, and integrating by parts
over $\mathbb{R}^3$, we have \be\label{apo-dt rho theta ttk}
\begin{split}
&\frac{1}{2}\frac{d}{dt}\int_{\mathbb{R}^3}\rho^k
|\theta^k_{tt}|^2+\kappa\int_{\mathbb{R}^3}|\nabla\theta^k_{tt}|^2\\&=-2\int_{\mathbb{R}^3}\rho^k_t
|\theta^k_{tt}|^2-\int_{\mathbb{R}^3}\rho^k_{tt}
\theta^k_t\theta^k_{tt}-\int_{\mathbb{R}^3}\rho^k_{tt}
u^{k-1}\cdot\nabla\theta^{k}\theta^k_{tt}-2\int_{\mathbb{R}^3}\rho^k_t
u^{k-1}_t\cdot\nabla\theta^{k}\theta_{tt}^k\\&-2\int_{\mathbb{R}^3}\rho^k_t
u^{k-1}\cdot\nabla\theta^{k}_t\theta^k_{tt}-\int_{\mathbb{R}^3}\rho^k
u^{k-1}_{tt}\cdot\nabla\theta^{k}\theta^k_{tt}-2\int_{\mathbb{R}^3}\rho^k
u^{k-1}_t\cdot\nabla\theta^{k}_t\theta_{tt}^k\\&-\int_{\mathbb{R}^3}\rho^k_{tt}\theta^{k}\mathrm{div}u^{k-1}\theta^k_{tt}
-2\int_{\mathbb{R}^3}\rho^k_t\theta^{k}_t\mathrm{div}u^{k-1}\theta_{tt}^k-2\int_{\mathbb{R}^3}\rho^k_t\theta^{k}\mathrm{div}u^{k-1}_t\theta_{tt}^k
\\&-\int_{\mathbb{R}^3}\rho^k\theta^{k}_{tt}\mathrm{div}u^{k-1}\theta_{tt}^k-2\int_{\mathbb{R}^3}\rho^k\theta^{k}_t\mathrm{div}u^{k-1}_t\theta_{tt}^k
-\int_{\mathbb{R}^3}\rho^k\theta^{k}\mathrm{div}u^{k-1}_{tt}\theta_{tt}^k\\&+\mu\int_{\mathbb{R}^3}\left(\nabla
u^{k-1}+(\nabla u^{k-1})^\prime\right):\left(\nabla
u^{k-1}_{tt}+(\nabla
u^{k-1}_{tt})^\prime\right)\theta_{tt}^k\\&+\mu\int_{\mathbb{R}^3}\left|\nabla
u^{k-1}_t+(\nabla
u^{k-1}_t)^\prime\right|^2\theta_{tt}^k+2\lambda\int_{\mathbb{R}^3}
\mathrm{div}u^{k-1}\mathrm{div}u^{k-1}_{tt}\theta_{tt}^k+2\lambda\int_{\mathbb{R}^3}|\mathrm{div}u^{k-1}_{t}|^2\theta_{tt}^k\\
&=\sum\limits_{i=1}^{17}VIV_i.
\end{split}
\ee From (\ref{apo-VII 1}), (\ref{apo-VII 2}), (\ref{apo-VII 3-5})
and (\ref{apo-VII 6-7}) with $u^k$ replaced by $\theta^k$, we obtain
\be\label{apo-VIV 1}
\begin{split}
VIV_1\le\frac{\kappa}{12}\int_{\mathbb{R}^3}|\nabla
\theta_{tt}^{k}|^2+\check{C}\int_{\mathbb{R}^3}\rho^{k}|\theta_{tt}^{k}|^2,
\end{split}
\ee and \be\label{apo-VIV 2}
\begin{split}VIV_2
\le&\frac{\kappa}{12}\int_{\mathbb{R}^3}|\nabla
\theta_{tt}^{k}|^2+\check{C}\|\nabla
\theta_t^{k}\|_{L^2}^2+\check{C}\|\nabla u^{k-1}_t\|_{L^2}^2\|\nabla
\theta_t^{k}\|_{L^2}\|\nabla^2
\theta_t^{k}\|_{L^2}\\&+\check{C}\|\sqrt{\rho^k}\theta_{tt}^{k}\|_{L^2}^2+\check{C}\|\nabla
u^{k-1}_t\|_{L^2}^4,
\end{split}
\ee and \be\label{apo-VIV 3-5}
\begin{split}
VIV_3+VIV_4+VIV_5 \le\frac{\kappa}{12}\int_{\mathbb{R}^3}|\nabla
\theta_{tt}^{k}|^2+\check{C}\|\rho_{tt}^{k}\|_{L^2}^2+\check{C}
\|\nabla u_t^{k-1}\|_{L^2}^2+\check{C} \|\nabla^2
\theta_t^{k}\|_{L^2}^2,
\end{split}
\ee and \be\label{apo-VIV 6-7}
\begin{split}
VIV_6+VIV_7 \le&\check{C}\int_{\mathbb{R}^3}|\nabla
u_{tt}^{k-1}|^2+\check{C}
\|\sqrt{\rho^{k}}\theta_{tt}^{k}\|_{L^2}^2+\check{C} \|\nabla
u_t^{k-1}\|_{L^2}^2\|\nabla \theta_t^{k}\|_{L^2}\|\nabla^2
\theta_t^{k}\|_{L^2}.
\end{split}
\ee Similarly, we have \be\label{apo-VIV 8-10}\begin{split}
&VIV_8+VIV_9+VIV_{10}\\
\le&\|\rho^k_{tt}\|_{L^2}\|\theta^{k}\|_{L^\infty}\|\mathrm{div}u^{k-1}\|_{L^3}\|\theta^k_{tt}\|_{L^6}+
2\|\rho^k_t\|_{L^2}\|\theta^{k}_t\|_{L^6}\|\mathrm{div}u^{k-1}\|_{L^6}\|\theta_{tt}^k\|_{L^6}\\&+2\|\rho^k_t\|_{L^2}\|\theta^{k}\|_{L^6}
\|\mathrm{div}u^{k-1}_t\|_{L^6}\|\theta_{tt}^k\|_{L^6}\\ \le&
\frac{\kappa}{12}\int_{\mathbb{R}^3}|\nabla
\theta_{tt}^{k}|^2+\check{C}\|\rho^k_{tt}\|_{L^2}^2+\check{C}\|\nabla\theta^{k}_t\|_{L^2}^2+\check{C}
\|\nabla^2u^{k-1}_t\|_{L^2}^2,
\end{split}
\ee and\be\label{apo-VIV 11-13}
\begin{split}
&VIV_{11}+VIV_{12}+VIV_{13}\\
\le&
\|\sqrt{\rho^k}\|_{L^\infty}\|\sqrt{\rho^k}\theta^{k}_{tt}\|_{L^2}\|\mathrm{div}u^{k-1}\|_{L^3}\|\theta_{tt}^k\|_{L^6}
+2\|\rho^k\theta^{k}_{t}\|_{L^3}\|\mathrm{div}u^{k-1}_t\|_{L^2}\|\theta^{k}_{tt}\|_{L^6}
\\&+\|\sqrt{\rho^k}\|_{L^\infty}\|\sqrt{\rho^k}\theta^{k}_{tt}\|_{L^2}\|\theta^{k}\|_{L^\infty}\|\mathrm{div}u^{k-1}_{tt}\|_{L^2}
\\
\le&
\check{C}\|\sqrt{\rho^k}\theta^{k}_{tt}\|_{L^2}\|\nabla\theta_{tt}^k\|_{L^2}
+\check{C}\|\rho^k\theta^{k}_{t}\|_{L^2}^\frac{1}{2}\|\rho^k\theta^{k}_{t}\|_{L^6}^\frac{1}{2}\|\mathrm{div}u^{k-1}_t\|_{L^2}\|\nabla\theta^{k}_{tt}\|_{L^2}
\\&+\check{C}\|\sqrt{\rho^k}\theta^{k}_{tt}\|_{L^2}\|\mathrm{div}u^{k-1}_{tt}\|_{L^2}\\
\le&\frac{\kappa}{12}\int_{\mathbb{R}^3}|\nabla
\theta_{tt}^{k}|^2+\check{C}\|\sqrt{\rho^k}\theta^{k}_{tt}\|_{L^2}^2+\check{C}\|\nabla\theta^{k}_{t}\|_{L^2}\|\nabla
u^{k-1}_t\|_{L^2}^2+\check{C}\|\nabla u^{k-1}_{tt}\|_{L^2}^2,
\end{split}
\ee and \be\label{apo-VIV 14-17}
\begin{split}
&VIV_{14}+VIV_{15}+VIV_{16}+VIV_{17}\\
\le&\check{C}\|\nabla u^{k-1}\|_{L^3}\|\nabla
u^{k-1}_{tt}\|_{L^2}\|\theta_{tt}^k\|_{L^6}+\check{C}\|\nabla
u^{k-1}_{t}\|_{L^2}\|\nabla
u^{k-1}_{t}\|_{L^3}\|\theta_{tt}^k\|_{L^6}\\
\le&\frac{\kappa}{12}\int_{\mathbb{R}^3}|\nabla
\theta_{tt}^{k}|^2+\check{C}\|\nabla
u^{k-1}_{tt}\|_{L^2}^2+\check{C}\|\nabla
u^{k-1}_{t}\|_{L^2}^2\|\nabla u^{k-1}_{t}\|_{L^2}\|\nabla^2
u^{k-1}_{t}\|_{L^2}.
\end{split} \ee
Substituting (\ref{apo-VIV 1}), (\ref{apo-VIV 2}), (\ref{apo-VIV
3-5}), (\ref{apo-VIV 6-7}), (\ref{apo-VIV 8-10}), (\ref{apo-VIV
11-13}) and (\ref{apo-VIV 14-17}) into (\ref{apo-dt rho theta ttk}),
and using the Cauchy inequality, we have
 \be\label{apo-dt rho theta ttk+1}
\begin{split}
&\frac{1}{2}\frac{d}{dt}\int_{\mathbb{R}^3}\rho^k
|\theta^k_{tt}|^2+\frac{\kappa}{2}\int_{\mathbb{R}^3}|\nabla\theta^k_{tt}|^2\\
\le&\check{C}\|\nabla \theta_t^{k}\|_{L^2}^2+\check{C}\|\nabla
u^{k-1}_t\|_{L^2}^4+\check{C}\|\nabla u^{k-1}_t\|_{L^2}^2\|\nabla
\theta_t^{k}\|_{L^2}\|\nabla^2
\theta_t^{k}\|_{L^2}\\&+\check{C}\|\rho_{tt}^{k}\|_{L^2}^2+\check{C}
\|\nabla^2 \theta_t^{k}\|_{L^2}^2+\check{C}
\|\sqrt{\rho^{k}}\theta_{tt}^{k}\|_{L^2}^2+\check{C}\|\nabla
u^{k-1}_{tt}\|_{L^2}^2\\&+\check{C}\|\nabla^2
u^{k-1}_{t}\|_{L^2}^4+\check{C}.
\end{split}
\ee Multiplying (\ref{apo-dt rho theta ttk+1}) by $t^3$, integrating
the result over $[0,t]$ for $t\in[0,T_0]$, and using
(\ref{apo-total}), (\ref{apo-le:utheta total}), (\ref{apo-cor7.4k}),
we have \bex
\begin{split}
\int_{\mathbb{R}^3}t^3\rho^k
|\theta^k_{tt}|^2+\int_0^{T_0}\int_{\mathbb{R}^3}t^3|\nabla\theta^k_{tt}|^2
\le\check{C}.
\end{split}
\eex The proof of Lemma \ref{apo-le:rho theta tt k} is complete.
\endpf
\begin{corollary}\label{apo-cor7.6k}
Under the condition of (\ref{local-initial data}), (\ref{apo-compa})
and (\ref{apo-g1 g2}), it holds that \be\label{apo-H4 thetak}
t^3\|\nabla^2\theta^k_{t}(\cdot,t)\|_{L^2}^2+t^3\|\nabla^4\theta^k(\cdot,t)\|_{L^2}^2
\le\check{C}, \ee  for any $k\ge1$ and a.e. $t\in[0,T_0]$.
\end{corollary}
\pf It follows from (\ref{apo-nabla 2theta t}), (\ref{apo-total}),
(\ref{apo-le:utheta total}), (\ref{apo-cor7.4k}), (\ref{apo-rho
theta ttk}) that\be\label{apo-tnabla 2theta t}
\begin{split}
t^\frac{3}{2}\|\nabla^2\theta^k_t(\cdot,t)\|_{L^2}\le&\check{C}t^\frac{3}{2}\|\sqrt{\rho^k}
\theta^k_{tt}(\cdot,t)\|_{L^2}+\check{C}t^\frac{3}{2}
\|\nabla\theta^k_t(\cdot,t)\|_{L^2}\\&+\check{C}
t^\frac{3}{2}\|\nabla
u^{k-1}_t(\cdot,t)\|_{L^2}+\check{C}t^\frac{3}{2}\|\nabla^2
u^{k-1}_t(\cdot,t)\|_{L^2}+\check{C}\\ \le&\check{C}.
\end{split}
\ee By (\ref{apo-H 4 theta}), (\ref{apo-total}), (\ref{apo-le:utheta
total}), (\ref{apo-cor7.4k}) and (\ref{apo-tnabla 2theta t}), we
have\be\label{apo-tH 4 theta}
\begin{split}t^\frac{3}{2}\|\nabla^4\theta^k(\cdot,t)\|_{L^2}\le&\check{C}
t^\frac{3}{2}\|\nabla\theta^k_t(\cdot,t)\|_{L^2}+\check{C}t^\frac{3}{2}\|\nabla^2\theta^k_t(\cdot,t)\|_{L^2}
\\&+\check{C}
t^\frac{3}{2}\|\nabla^3\theta^{k}(\cdot,t)\|_{L^2}
+\check{C}t^\frac{3}{2}\|\nabla u^{k-1}(\cdot,t)\|_{H^2}^2+\check{C}\\
\le&\check{C}.
\end{split}
\ee
\endpf

{\noindent\bf Step 2: Completion of proof of Appendix B.}\\

Using some arguments similar to \cite{cho-Kim: perfect gas}, we
obtain that the full sequence $(\rho^k,u^k,\theta^k)$ converges to a
limit $(\rho^\delta,u^\delta,\theta^\delta)$ for any given
$\delta\in(0,1)$ in the following strong sense: \bex\begin{cases}
\rho^k\rightarrow\rho^\delta\ \mathrm{in}\ L^\infty([0,T_0];
L^2),\ \mathrm{as}\ k\rightarrow\infty,\\
(u^k,\theta^k)\rightarrow(u^\delta,\theta^\delta)\ \mathrm{in}\
L^2([0,T_0]; D_0^1),\ \mathrm{as}\ k\rightarrow\infty,
\end{cases}
\eex and $(\rho^\delta,u^\delta,\theta^\delta)$ is the unique
solution to (\ref{full N-S+1})-(\ref{3d-initial}) with initial data
replaced by $(\rho_0^\delta, u_0, \theta_0)$, where $\rho^\delta>0$
and $\theta^\delta\ge0$. With Lemmas \ref{apo-le: 7.1},
\ref{apo-le:7.2}, \ref{apo-le:7.3} and \ref{apo-le:rho theta tt k},
and Corollaries \ref{apo-cor7.4} and \ref{apo-cor7.6k}, and the
lower semi-continuity of the norms, we have
 \be\label{3d-u theta delta}\begin{split}
 &\|\rho^\delta(\cdot,t)-\delta\|_{H^2\cap
W^{2,q}}+\|\rho^\delta_t(\cdot,t)\|_{H^1}+\|\nabla
 u^\delta(\cdot,t)\|_{H^1}+\\&\|\nabla\theta^\delta(\cdot,t)\|_{H^1}+\|\nabla^3u^\delta\|_{L^2([0,t];L^2)}
 +\|\nabla^3\theta^\delta\|_{L^2([0,t];L^2)}\le\check{C},\end{split}\ee
\be\label{3d-u t theta t delta}\begin{split}&\|\sqrt{\rho^\delta}
u_t^\delta(\cdot,t)\|_{L^2}+\|\sqrt{\rho^\delta}
\theta_t^\delta(\cdot,t)\|_{L^2}+\|\nabla
u_t^\delta\|_{L^2([0,t];L^2)}+\\&\|\nabla
\theta_t^\delta\|_{L^2([0,t];L^2)}\le\check{C}, \end{split}\ee
 \be\label{apo-total delta}\begin{split}
\sqrt{t}\|\nabla u^\delta(\cdot,t)\|_{H^2}+\sqrt{t}\|\nabla
u_t^\delta(\cdot,t)\|_{L^2}+\int_0^{T_0}\left(t\|\sqrt{\rho^\delta}
u_{tt}^\delta\|_{L^2}^2+t\|\nabla^2u_t^\delta\|_{L^2}^2
 \right)dt\le\check{C},\end{split}\ee
\be\label{apo-cor7.4k delta}\begin{split}
&t\|\nabla^2u_t^\delta(\cdot,t)\|_{L^2}+t\|\nabla^3\theta^\delta(\cdot,t)\|_{L^2}+t\|\nabla^3u^\delta(\cdot,t)\|_{L^q}
+t\|\nabla\theta_t^\delta(\cdot,t)\|_{L^2}+t\|\sqrt{\rho^{\delta}}
u_{tt}^{\delta}(\cdot,t)\|_{L^2}\\& +\int_0^{T_0}t^2\left(
\|\nabla^2\theta^\delta_t\|_{L^2}^2+\|\nabla^4\theta^\delta\|_{L^2}^2+\|\sqrt{\rho^\delta}
\theta^\delta_{tt}\|_{L^2}^2+\|\nabla
u_{tt}^{\delta}\|_{L^2}^2\right)dt\le\check{C},
\end{split}\ee  and
 \be\label{apo-rho theta ttk delta}
t^3\|\nabla^2\theta^\delta_{t}(\cdot,t)\|_{L^2}^2+t^3\|\nabla^4\theta^\delta(\cdot,t)\|_{L^2}^2+t^3\|\sqrt{\rho^\delta}
\theta^\delta_{tt}\|_{L^2}^2+\int_0^{T_0}\int_{\mathbb{R}^3}t^3|\nabla\theta^\delta_{tt}|^2\,dx\,dt
\le\check{C}. \ee By (\ref{3d-u theta delta}), (\ref{3d-u t theta t
delta}), (\ref{apo-total delta}), (\ref{apo-cor7.4k delta}),
(\ref{apo-rho theta ttk delta}), we pass
$(\rho^\delta,u^\delta,\theta^\delta)$ to a limit $(\rho,u,\theta)$
(take subsequence if necessary) which is the unique solution to
(\ref{full N-S+1})-(\ref{3d-boundary}). By the lower semi-continuity
of the norms and some arguments which are concerned with the
time-continuity of the solutions as in \cite{Cho-Kim, cho-Kim:
perfect gas} and references therein, we get
(\ref{apo-regularities}).

\section*{Acknowledgements}
The authors would like to thank the anonymous referees for their
helpful suggestions. Wen was supported by the National Natural
Science Foundation of China $\#$11301205 and $\#$11671150 and by
the Fundamental Research Funds for the Central Universities
$\#$D2154560. Zhu was supported by the National Natural Science
Foundation of China $\#$11331005, the Program for Changjiang
Scholars and Innovative Research Team in University $\#$IRT13066,
and the Special Fund Basic Scientific Research of Central Colleges
$\#$CCNU12C01001.

\end{document}